\documentclass[12pt]{amsart}
\usepackage{amssymb,amsmath,amsthm,slashed}
\newcommand{\leg}[2]{\genfrac{(}{)}{}{}{#1}{#2}}

\newtheorem{theorem}{Theorem}
\newtheorem{lemma}[theorem]{Lemma}

\newtheorem*{conjecture}{\bf Conjecture}
\newtheorem{proposition}[theorem]{Proposition}
\newtheorem{definition}[theorem]{Definition}
\newtheorem{fact}[theorem]{Fact}

\theoremstyle{remark}

\newtheorem*{remark}{Remark}
\newtheorem*{example}{Example}
\numberwithin{theorem}{section} \numberwithin{equation}{section}

\newcommand{\R}{\mathbb{R}}
\newcommand{\C}{\mathbb{C}}
\newcommand{\D}{{\text {\bf D}}}

\newcommand{\Z}{\mathbb{Z}}

\newcommand{\sgn}{\operatorname{sgn}}

\newcommand{\tr}{{\text {\rm tr}}}

\newcommand{\re}{\textnormal{Re}}
\newcommand{\im}{\textnormal{Im}}
\def\H{\mathbb{H}}
\begin{document}
\title[$\mathrm{SO}(3)$-Donaldson invariants of $\mathbb{C}\mathrm{P}^2$ and Mock Theta Functions]
{$\mathrm{SO}(3)$-Donaldson invariants of $\mathbb{C}\mathrm{P}^2$ and Mock Theta Functions}
\author{Andreas Malmendier and Ken Ono}
\address{Department of Mathematics, Colby College,
Waterville, Maine 04901}
\email{andreas.malmendier@colby.edu}
\address{Department of Mathematics and Computer Science, Emory University,
Atlanta, Georgia 30322} \email{ono@mathcs.emory.edu}
\thanks{The second author thanks
the generous support of the National Science Foundation, the
Manasse family, and the Hilldale Foundation, and the Candler Fund.}
\begin{abstract}
We compute the Moore-Witten regularized $u$-plane integral
on $\mathbb{C}\mathrm{P}^2$,
and we confirm the conjecture
that it is the generating function for the $\mathrm{SO}(3)$-Donaldson
invariants of $\mathbb{C}\mathrm{P}^2$.
We also derive generating functions for the $\mathrm{SO}(3)$-Donaldson invariants with
$2 N_f$ massless monopoles using the geometry of certain rational elliptic surfaces ($N_f \in \{0,2,3,4\}$),
and we
show that the partition function for $N_f=4$
is nearly modular.
Our results rely heavily on the theory of mock theta functions and harmonic Maass
forms (for example, see \cite{OnoCDM}).
\end{abstract}
\maketitle
\section{Introduction and Statement of Results}
In his plenary address at the Ramanujan Centenary Conference at the
University of Illinois (Urbana-Champaign) in 1987, Freeman Dyson
proclaimed his hope \cite{DysonGarden} that the ``theory" of
Ramanujan's mock theta functions would someday play a role
in mathematical physics.
\bigskip

\noindent ``{\it The mock theta-functions give us tantalizing hints
of a grand synthesis still to be discovered. Somehow it should be
possible to build them into a coherent group-theoretical structure,
analogous to the structure of modular forms which Hecke built around
the old theta-functions of Jacobi. This remains a challenge for the
future. My dream is that I will live to see the day when our young
physicists, struggling to bring the predictions of superstring
theory into correspondence with the facts of nature, will be led to
enlarge their analytic machinery to include mock
theta-functions...But before this can happen, the purely
mathematical exploration of the mock-modular forms and their
mock-symmetries must be carried a great deal further.''}

\smallskip \ \ \hskip4.in Freeman Dyson
\medskip

The 2002 Ph.D. thesis of Sander Zwegers \cite{Z2} has provided
Dyson's desired ``coherent group-theoretical structure", and it
turns out that Ramanujan's mock theta functions are {\it holomorphic
parts} of weight 1/2 {\it harmonic weak Maass forms} (see
Section~\ref{HarmonicMaassForms} for definitions).
Zwegers's thesis has sparked a flurry of recent activity involving
such Maass forms. Indeed, harmonic Maass forms are now known to play
a central role in the study of Ramanujan's mock theta functions, as
well as other important mathematical topics: Borcherds products,
derivatives of modular $L$-functions, Gross-Zagier formulas and
Faltings heights of CM cycles, partitions, and traces of singular
moduli (see \cite{BO1, BO2, BOR, Br, BF, BrO, BrYang, OnoCDM,ZagierBourbaki, Z2}). 
Here we give an application, which combined
with earlier important works of G\"ottsche and Zagier \cite{Goettsche, GoettscheZagier},
realizes Dyson's original hope that
mock theta functions would one day play a role in mathematical
physics.
We shall relate mock theta functions to Donaldson invariants of a smooth simply connected four-dimensional manifold \cite{DonaldsonKronheimer}.
There are two families of Donaldson invariants, corresponding to the $\mathrm{SU}(2)$-gauge
theory and the $\mathrm{SO}(3)$-gauge theory with non-trivial Stiefel-Whitney class. In each case, the Donaldson invariants are graded homogeneous
polynomials on the homology $H_0(\mathbb{C}\mathrm{P}^2) \oplus H_2(\mathbb{C}\mathrm{P}^2)$, where $H_i(\mathbb{C}\mathrm{P}^2)$ is considered to
have degree $(4-i)/2$, defined using the fundamental homology classes of the corresponding moduli spaces of anti-selfdual instantons arising in gauge
theory. The Donaldson invariants depend by definition on the choice of a Riemannian metric. In the case $b_2^+ > 1$, the Donaldson invariants
are independent of the metric as long as it is generic.  In the case $b_2^+=1$, Kotschick and Morgan \cite{KotschickMorgan} showed that the Donaldson
invariants depended on the metric but only through the chamber of its period point in the positive cone in $H^2(X,\mathbb{R})$. Crossing through a wall
in $H^2(X,\mathbb{R})$ defined by a class $\xi \in H^2(X,\mathbb{Z})$ adds a certain wall-crossing term $\delta^X_{\xi,k}$ to the Donaldson invariants.
Kotschick and Morgan made the following conjecture:
$\delta^X_{\xi,k}$ is a polynomial in the multiplication by $\xi$ and the quadratic form $Q_X$ on $H_2(X,\mathbb{Z})$ whose coefficients
depend only on $\xi^2$, the instanton number $k$, and the homotopy type on $X$.
We restrict ourselves to the simplest manifold to which Donaldson's definition applies, the complex projective plane
$\mathbb{C}\mathrm{P}^2$ with the Fubini-Study metric.
In earlier work, G\"ottsche and Zagier \cite{GoettscheZagier} gave a formula for the Donaldson
invariants of rational surfaces in terms of theta functions of indefinite lattices.
These are Jacobi forms for special choices of the polarization on the
boundary of the positive cone. As an application, G\"ottsche \cite{Goettsche} derived
closed expressions for the two families of the Donaldson invariants of $\mathbb{C}\mathrm{P}^2$
assuming the Kotschick-Morgan conjecture. In the recent work of G\"ottsche et al. \cite{GoettscheNakajimaYoshioka},
it was shown that the assumption of the Kotschick-Morgan conjecture was \emph{not} necessary for the formula
for the Donaldson invariants.
The present paper concerns deep conjectural relations between these works and
constructions in theoretical physics which are presently the focus of much study.
From the viewpoint of theoretical physics \cite{Witten1}, these
two families of Donaldson invariants and the related Seiberg-Witten invariants
are the correlation functions of a supersymmetric topological
gauge theory with gauge group $\mathrm{SU}(2)$ and $\mathrm{SO}(3)$ respectively.
Using physical considerations, Witten \cite{Witten2} argued that one should be able
to compute these correlation
functions in a so called {\it low
energy effective field theory} instead. This theory has the advantage of being an
{\em abelian} $\mathcal{N}=2$ supersymmetric topological gauge theory.
The data required to define the theory
only involves line bundles of even (resp. odd) first Chern class on
$\mathbb{C}\mathrm{P}^2$ if the gauge group is $\mathrm{SU}(2)$ (resp. $\mathrm{SO}(3)$).
The vacua of the low energy effective field theory are parametrized by the $u$-plane which is a certain analytically marked
rational elliptic surface. Seiberg and Witten \cite{SeibergWitten1} argued
further that the rational elliptic surface is the
modular curve $\H/\Gamma_0(4)$,
together with a meromorphic one-form. Moreover, Moore and Witten \cite{MooreWitten} obtained the correlation functions as a regularized integral over the
$u$-plane. The integrand is a modular invariant function which is determined by the marked elliptic surface and the gauge group. The regularization procedure
depends on interpreting the
integrand as a total derivative, combined with
constant term contributions from cusps.
This integration by parts naturally
introduces non-holomorphic modular forms of weight $3/2$ (resp. $1/2$)
for the gauge group $\mathrm{SU}(2)$ (resp. $\mathrm{SO}(3)$).
Thus, the regularized $u$-plane integral defines a way of extracting
certain contributions for each boundary
component near the cusps at $\tau=0,2,\infty$ of the rational elliptic surface. Moore and Witten observed \cite{MooreWitten} that the
cuspidal
contributions
at $\tau=0,2$ vanish trivially, which
coincides with the mathematical statement that
all Seiberg-Witten invariants on $\mathbb{C}\mathrm{P}^2$
vanish due to the presence of a
Fubini-Study metric of positive scalar curvature \cite{Witten3}.
Moore and Witten went further and made the following
comprehensive conjecture:
\begin{conjecture}[Moore and Witten \cite{MooreWitten}]
The contribution from the cusp at $\tau=\infty$ to the regularized $u$-plane integral is the generating function for the
Donaldson invariants of $\mathbb{C}\mathrm{P}^2$.
\end{conjecture}
As evidence for this conjecture,
in the case of the gauge group $\mathrm{SU}(2)$, Moore and Witten \cite{MooreWitten} computed the first 40 invariants and found them to be in agreement
with the results of Ellingsrud and G\"ottsche \cite{EllingsrudGoettsche}.
However, the conjecture remains open.
\begin{remark}
Moore and Witten \cite{MooreWitten} also showed that the
general equality for $\mathrm{SU}(2)$ implies some interesting
relations involving the classical Hurwitz class numbers which
arise in number theory. It is in this way that harmonic
Maass forms make their first appearance. The holomorphic
part of Zagier's weight 3/2 Maass-Eisenstein series,
which first arose \cite{HirzebruchZagier} in connection with
intersection theory for certain Hilbert modular
surfaces, is the generating function for Hurwitz class numbers.
\end{remark}
Our main result concerns the case of the $\mathrm{SO}(3)$-gauge theory
where little was known about the conjecture.
We prove the following theorem:
\begin{theorem}\label{MainTheorem}
The conjecture of Moore and Witten in the case of the $\mathrm{SO}(3)$-gauge theory on $\mathbb{C}\mathrm{P}^2$ is true.
\end{theorem}
Donaldson theory can be generalized by introducing sections of
a $\operatorname{spin}_{\mathbb{C}}$-structure on $\mathbb{C}\mathrm{P}^2$ coupled to the $\mathrm{SO}(3)$-gauge bundle.
In the physics literature, these additional fields are called {\it massless monopoles},
and their number is denoted by $2 \, N_f$, where $0\leq N_f\leq 4$.
The case where $N_f=0$ corresponds to the
original $\mathrm{SO}(3)$-Donaldson invariants\footnote{The number of sections is even and $N_f \le 4$;
otherwise the quantum theory is inconsistent \cite{MooreWitten}.}.
The corresponding moduli spaces and invariants
when $N_f>0$ have not been given much consideration in mathematics.
However in physics,
Seiberg and Witten \cite{SeibergWitten2} argued that including the monopoles
changes the rational elliptic surface. For $N_f>0$ the modular elliptic surface for $\Gamma_0(4)$ is replaced by rational elliptic
surfaces determined by $N_f$ vectors in the $E_8$ root lattice.
Rather nicely, the integrand of the $u$-plane integral remains a modular invariant
when the gauge group is $\mathrm{SO}(3)$.
We classify all rational elliptic surfaces which appear as Seiberg-Witten curves for $0 \le N_f \le 4$ (see Lemma~\ref{sw-curves}).
For $N_f=2, 3$, it turns out that the elliptic surfaces for the Donaldson theory with massless monopoles are still modular elliptic surfaces related
to the elliptic surface for $N_f=0$ by two-isogeny or a twist of the elliptic fibration.
However, in the case where $N_f=1$ the surface is not modular elliptic (see Lemma \ref{presentation} and \cite{Nahm} for a physical explanation).
In the case $N_f=4$, elliptic surfaces for the Donaldson theory with massive monopoles can be obtained from twisting the fibration of the modular elliptic surface in the case $N_f=2$. We give explicit Weierstrass representations (see Lemmas \ref{presentation} and
\ref{presentation2}) for the cases $N_f=2,3,4$, and we compute the generating functions (see
Theorems \ref{generating_function} and \ref{partition_function_tau_a}) using the regularized $u$-plane integral.
We conjecture:
\begin{conjecture}
The generating function in (\ref{generating_function}) is the generating function (\ref{generatingfunction2}) for
the $\mathrm{SO}(3)$-Donaldson invariants of $\mathbb{C}\mathrm{P}^2$ with $2 \, N_f$ monopoles for $N_f=2,3$.
\end{conjecture}
\noindent
Using the Atiyah-Singer index theorem
and Theorem \ref{MainTheorem}, we can also show that the conjecture
passes a non-trivial test proving that the highest order contributions to
each coefficient in the partition functions agree for $N_f=2, 3$.
In the $N_f=4$ case, instead of asymptotic freedom one has conformal invariance of the physical system. Thus, the partition
function is a function of the complex coupling constant $\tau$. We explicitly compute the partition
function. To this end, we shall define a power series in $q=\exp(2\pi i \tau)$ (see Sections 4 and 7)
\begin{equation*}
 Q^+\left(q\right) = \frac{1}{q^{\frac{1}{8}}}\left( 1 + 28 \, q^{\frac{1}{2}} + 39 \, q + 196 \, q^{\frac{3}{2}} + 161 \, q^2 + \dots \right) \;.
\end{equation*}
Adding a non-holomorphic part $Q^-(q)$ we obtain the Maass form $Q(q) = Q^+(q) + Q^-(q)$ (see Section 4).
In terms of the function $\pmb{\mathrm{Z}}(\tau)=\eta^3(\tau) \, Q(q)$, we prove the following theorem:
\begin{theorem}
The regularized partition function of the massless $N_f=4$ low energy effective field theory on $\mathbb{C}\mathrm{P}^2$
is well defined, and is given by
\begin{equation}
\begin{split}\label{partition_function_Nf4_1st}
  \pmb{\mathrm{Z}}^{4}_{\mathrm{UP}}  =  \left[ \frac{1}{2}\, \frac{q}{\eta^4(\tau)} \, \frac{d}{dq} \left( \frac{q}{\eta^4(\tau)}  \, \frac{d}{dq} \right)
+ g(\tau)  \right] \, \frac{\pmb{\mathrm{Z}}(\tau)}{\eta^4(\tau)}\;,
\end{split}
\end{equation}
where
\begin{equation*}
  g(\tau) =
- \dfrac{1}{2^2 \, 3^2 } \, \left[ \left(\frac{\vartheta_2(\tau)}{\eta(\tau)}\right)^8 - \left(\frac{\vartheta_2(\tau)}{\eta(\tau)}\right)^4
\left(\frac{\vartheta_3(\tau)}{\eta(\tau)}\right)^4 + \left(\frac{\vartheta_3(\tau)}{\eta(\tau)}\right)^8 \right]\;.
\end{equation*}
\end{theorem}
Here we discuss this result in the context of the theory of modular forms.
Recall that the Donaldson polynomials arise as correlation functions of the $\mathcal{N}=2$ supersymmetric topological gauge theory
with gauge group $\mathrm{SU}(2)$ or $\mathrm{SO}(3)$. There are also quantum field theories on $\mathbb{C}\mathrm{P}^2$ with more
supersymmetry. Roughly speaking, there are three choices for a $\mathcal{N}=4$ supersymmetric topological gauge theory, determined by
what is referred to as {\it the topological twist}. One particular choice, called the {\it Vafa-Witten twist}, has been studied
extensively in the physics literature \cite{VafaWitten}, and its partition function has been computed in many examples.
It turns out that the coefficient of $q^k$ in the partition function
equals the Euler characteristic of the instanton moduli space of instanton number $k$. The
{\it S-duality} sends $\tau \mapsto - 1/\tau$  and exchanges electric and magnetic fields. In the $\mathcal{N}=4$ theory, the S-duality
is a symmetry of the physical theory whence the partition function is in fact a weakly holomorphic weight $0$ modular form.
Indeed, Vafa and Witten even point out the strong relation with the non-holomorphic
weight $3/2$ Eisenstein series of Zagier which is perhaps the first
half-integral weight harmonic Maass form which was seriously investigated.
In the $\mathcal{N}=2$ theory, the S-duality is not a symmetry. Accordingly, we see that the partition function
is not modular, but only a ``piece'' of a modular invariant, a
fact which is central to our proof of Theorem \ref{MainTheorem}. It has the property that its ``trace", in a very natural sense, gives rise to
modular invariants (see Theorem \ref{partition_function_tau}) which resemble some modular partition
functions computed by Vafa and Witten (see \cite[Sec. 42.]{VafaWitten}).
The partition function (\ref{partition_function_Nf4_1st}) is naturally related to the mock theta function
\begin{equation*}
M(q^8):=q^{-1}\sum_{n=0}^{\infty}\frac{(-1)^{n+1}q^{8(n+1)^2}\prod_{k=1}^n
(1-q^{16k-8})}{\prod_{k=1}^{n+1}(1+q^{16k-8})^2}=-q^7+2q^{15}-3q^{23}+\cdots.
\end{equation*}
In Section \ref{qIdentities} we prove that
$$
 Q(q)-4 \, M(q),
$$
is a weight $1/2$ weakly holomorphic modular form. Therefore, we then find that
$$
\frac{1}{\eta(\tau)}\cdot \Big( Q(q) - 4 \, M(q)\Big)
$$
is the desired modular invariant, a weight $0$ modular form.
This paper is organized as follows. In Section~\ref{Section2} we recall
facts about $\mathrm{SO}(3)$-Donaldson invariants on $\mathbb{C}\mathrm{P}^2$, and we recall
certain generating functions of
G\"ottsche which play a central role in the proof of Theorem~\ref{MainTheorem}.
In Section~\ref{Section3}
we recall and classify the Seiberg-Witten curves, and in Section~\ref{Section4}
we use their properties to define the $u$-plane integral, which in turn
we employ to construct the vital generating series when $N_f\in \{0, 2, 3, 4\}$.
The proof of Theorem~\ref{MainTheorem} is then reduced to a criterion
(see Theorem~\ref{criterion})
involving these new generating functions, in the case
where $N_f=0$, and those of
G\"ottsche. To verify this criterion, we employ the theory of
harmonic Maass forms which is briefly recalled in
Section~\ref{HarmonicMaassForms}. Specifically, we relate the relevant
generating functions to weight 1/2 harmonic Maass forms
which are described using a non-holomorphic Jacobi form
constructed by Zwegers.
In Section~\ref{Section6} we recall the construction of this Jacobi form,
and in
Section~\ref{qIdentities} we use it to represent
the generating functions in terms of explicit harmonic Maass forms.
In Section~\ref{TheProof} we conclude with the proof of
Theorem~\ref{MainTheorem}, and we give some numerical examples.
\section*{Acknowledgments}
The authors thank Kathrin Bringmann, Amanda Folsom,
Lothar G\"ottsche, David Morrison, Marcos Mari\~no and Brian Rice  for their comments on
earlier versions of this paper. The authors thank Michael Griffin who simplified
the deduction of Theorem 1.1 from Theorem 4.26. The authors also thank the referee
for many helpful suggestions which improved the paper.
\section{$\mathrm{SO}(3)$-Donaldson invariants on $\mathbb{C}\mathrm{P}^2$}\label{Section2}
In this section we recall basic facts about the $\mathrm{SO}(3)$-Donaldson invariants on $\mathbb{C}\mathrm{P}^2$ (cf. \cite{ FreedUhlenbeck, Goettsche2, Kotschick}),
and we recall
a closed formula expression for these invariants which is due to G\"ottsche \cite{Goettsche}. We also discuss briefly
the $\mathrm{SO}(3)$-Donaldson theory with massless monopoles.
\subsection{The generating function.}\label{Sec_Donaldson}
The Fubini-Study metric $g$ on $\mathbb{C}\mathrm{P}^2$ is K\"ahler with the
K\"ahler form $K= \frac{i}{2} g_{a\bar{b}} \, dz^a \wedge dz^{\bar{b}}$.
We denote the first Chern class of the dual of the hyperplane bundle over
$\mathbb{C}\mathrm{P}^2$ by $\operatorname{H} = K/\pi$,  so that
$\int_{\mathbb{C}\mathrm{P}^2} \operatorname{H}^2 =1$, $c_1(\mathbb{C}\mathrm{P}^2)=3 \operatorname{H}$, and $p_1(\mathbb{C}\mathrm{P}^2)=3 \operatorname{H}^2$.
The Poincar\'e dual $\check{\operatorname{H}}$ of $\operatorname{H}$ is a generator of the rank-one homology group
$H_2(\mathbb{C}\mathrm{P}^2)$. $\mathrm{SO}(3)$-bundles on four-dimensional manifolds are classified by the second Stiefel-Whitney class
$w_2(P) \in H^2(\mathbb{C}\mathrm{P}^2;\mathbb{Z}_2)$ and the first Pontrjagin class $p_1(P) \in H^4(\mathbb{C}\mathrm{P}^2)$, such that
\begin{equation*}
 p_1(P)[\mathbb{C}\mathrm{P}^2] \equiv w_2^2(P)[\mathbb{C}\mathrm{P}^2]  \mod 4\;.
\end{equation*}
Since $\mathbb{C}\mathrm{P}^2$ is simply connected there is an integer class
$\alpha \in H^2(\mathbb{C}\mathrm{P}^2)$ whose mod-two reduction is $w_2(P)$.
Then, there is a smooth complex two-dimensional vector bundle $\xi \to \mathbb{C}\mathrm{P}^2$ with the Chern classes
$c_1(\xi)=\alpha$ and $c_2(\xi)$, such that $c^2_1(\xi) -4 \, c_2(\xi)=p_1(P)$.
The bundle $\xi$ can be reduced to an $\mathrm{SU}(2)$-bundle iff $c_1(\xi)=0$. Thus, a $\mathrm{SO}(3)$-bundle which does not
arise as the associated bundle for the adjoint representation of a $\mathrm{SU}(2)$-bundle has to satisfy $w_2(P) \not = 0$.
From now on, we will always assume $c_1(\xi)=-\operatorname{H}$ and $c_2(\xi)=k \operatorname{H}^2$. Once a Hermitian metric on $\xi$ is fixed, let $\mathfrak{A}$
be the set of compatible unitary connections $d_A = d+ A$ on $\xi$, and $\mathfrak{G}$ the group of gauge transformations.
The metric
on $\mathbb{C}\mathrm{P}^2$ gives rise to a Hodge star operator $*:\Lambda^p_{\mathbb{C}\mathrm{P}^2} \to
\Lambda^{4-p}_{\mathbb{C}\mathrm{P}^2}$ with $*^2 = (-1)^{p}$. On $p$-forms, the adjoint operator of $d_A$
is $d_A^* = - * d_A*$. Let $\operatorname{End}_0{\xi}$ be the
bundle of the traceless endomorphisms of $\xi$.
In terms of the curvature $F_A \in \Lambda^2_{\mathbb{C}\mathrm{P}^2}(\operatorname{End}_0{\xi})$ on $\xi$,
the Pontrjagin number is given by
\begin{eqnarray}
 \Big(c^2_1(\xi)- 4 \, c_2(\xi)\Big)[\mathbb{C}\mathrm{P}^2]  = - \frac{1}{4\pi^2} \int_{\mathbb{C}\mathrm{P}^2} \; \tr\left(F_A \wedge F_A\right) \;.
\end{eqnarray}
A connection is called {\it anti-selfdual}
if the curvature satisfies $F_A^+=0$ (i.e. $* F_A = - F_A$).
A connection is called {\it reducible} if the bundle $P \to \mathbb{C}\mathrm{P}^2$ reduces to the direct sum $\lambda \oplus \epsilon$ of a
line bundle $\lambda \to \mathbb{C}\mathrm{P}^2$ and a trivial oriented real line bundle $\epsilon\to \mathbb{C}\mathrm{P}^2$; otherwise it is called {\it irreducible}.
We denote by $\mathfrak{A}^*$ the
space of irreducible connections. Since the gauge group $\mathfrak{G}$ acts freely on $\mathfrak{A}^*$,
the space $\mathfrak{B}^* = \mathfrak{A}^* / \mathfrak{G}$ is a Banach manifold.
The action functional of the Yang-Mills theory is
\begin{equation}
\label{YM}
  \int_{\mathbb{C}\mathrm{P}^2} \; \tr\left(F_A \wedge * F_A\right),
\end{equation}
and has the Euler-Lagrange equation $d_A^* F_A = 0$. For $k>0$, the Yang-Mills action (\ref{YM}) is minimized by the anti-selfdual connections.
The moduli space of anti-selfdual irreducible connections with $c_2(\xi)[\mathbb{C}\mathrm{P}^2]=k$ modulo gauge transformations is denoted by $\mathfrak{M}(-1,k)$.
The Atiyah-Singer complex of $d_A$ in dimension four is the three term complex
\begin{eqnarray}
\label{AHS}
 \Lambda^0_{\mathbb{C}\mathrm{P}^2}(\xi) \xrightarrow{d_A} \Lambda^1_{\mathbb{C}\mathrm{P}^2} (\xi) \xrightarrow{\pi_{+}
 d_A} \Lambda^{2+}_{\mathbb{C}\mathrm{P}^2} (\xi),
\end{eqnarray}
where the map $\pi_{+}$ is the orthogonal projection on the self-dual two-forms. It is a complex if $F_A$ is
anti-selfdual. The associated elliptic operator of Laplace type is a Fredholm operator, and the index of the elliptic complex equals
$ \operatorname{ind} \, d_A  = 2(4k-1) - 3(1+b_2^+)$. A local model of a regular neighborhood of $[A] \in \mathfrak{M}(-1,k)$ is given by the
intersection of the slice $\ker d_A^* \in T_A\mathfrak{A}^*$, which is locally orthogonal to
pure gauge transformations, with $\ker (\pi_+ d_A )$, which is the linearization of $F_A^+=0$.
Since $c_1^2(\xi) - 4 c_2(\xi) \not \equiv 0 \mod{8}$ the moduli spaces $\mathfrak{M}(-1,k)$ of rank-two stable vector bundles
on $\mathbb{C}\mathrm{P}^2$ are smooth, projective varieties of dimension $2d_k= 8 k - 8$ \cite{Stromme}.
\begin{remark}
It is known that the moduli space $\mathfrak{M}(c_1,c_2)$ of rank-two stable vector bundles $\xi$ over $\mathbb{C}\mathrm{P}^2$ with Chern classes $c_1, \, c_2$ only depends on the
discriminant $c_1^2 - 4 c_2$, with the discriminant being negative for stable bundles. The isomorphism between $\mathfrak{M}(c_1,c_2)$ and $\mathfrak{M}(c_1-2r,c_2- r \, c_1 + r^2)$
is given by twisting $\xi \mapsto \xi \otimes \mathcal{O}^*_{\mathbb{C}\mathrm{P}^2}(r)$ \cite{Klyachko}.
\end{remark}
On $\mathbb{C}\mathrm{P}^2 \times \mathfrak{B}^*$, there exists a universal bundle with a
$\mathrm{SO}(3)$-connection \cite{AtiyahSinger}; it
is the vector bundle $\mathfrak{P}$ with the structure group $\mathrm{SO}(3)$ in the commutative diagram:
\begin{equation}
\label{universal_bundle}
\begin{array}{ccc}
 P \times \mathfrak{A}^* & \rightarrow & \mathfrak{P} = (P \times \mathfrak{A}^*)/\mathfrak{G} \\
 \downarrow && \downarrow\\
 \mathbb{C}\mathrm{P}^2 \times \mathfrak{A}^*  & \rightarrow & \mathbb{C}\mathrm{P}^2 \times \mathfrak{B}^*
\end{array}
\end{equation}
where the action of $\mathfrak{G}$ and $\mathrm{SO}(3)$ on $P \times \mathfrak{A}^*$ commutes.
Using the connection $A \in \mathfrak{A}^*$ on $P$ and the canonical connection on $\mathfrak{A}^* \to \mathfrak{B}^*$,
one defines a $\mathfrak{G}$-invariant connection  on $P \times \mathfrak{A}^* \to \mathbb{C}\mathrm{P}^2 \times \mathfrak{A}^*$.
The connection and its curvature $\mathfrak{F}$ descend to a connection and a curvature form on the quotient bundle.
The Pontrjagin number of the universal bundle can be decomposed into its components
\begin{equation}
\label{universal_curv}
p_1(\mathfrak{P})  =  - \frac{1}{4\pi^2} \; \tr \left( \mathfrak{F}\wedge\mathfrak{F} \right)  = -4 \; \sum_{r=0}^4
\mathfrak{W}^{r,4-r},
\end{equation}
with $\mathfrak{W}^{r,4-r}  \in  \Lambda^{r,4-r}_{\mathbb{C}\mathrm{P}^2 \times \mathfrak{B}^*}$.
We evaluate the class $- \frac{1}{4} \, p_1(\mathfrak{P})$ on the $2$-cycle $[\check{\operatorname{H}}] \in H_2(\mathbb{C}\mathrm{P}^2)$ and
the $0$-cycle $[\operatorname{pt}] \in H_0(\mathbb{C}\mathrm{P}^2)$ to obtain
\begin{equation}
\label{observable}
 \mu(\check{\operatorname{H}})  =  \int_{\check{\operatorname{H}}}  \;\mathfrak{W}^{2,2} \in H^2(\mathfrak{B}^*;\mathbb{Q}) \;, \qquad
 \mu(\operatorname{pt}) =  \mathfrak{W}^{0,4}(\operatorname{pt}) \in H^4(\mathfrak{B}^*;\mathbb{Q}) \;.
\end{equation}
By applying \cite[Proposition 5.1.15]{DonaldsonKronheimer} to $\mathbb{C}\mathrm{P}^2$ is follows:
\begin{proposition}
\label{ring}
The rational cohomology ring $H^*(\mathfrak{B}^*; \mathbb{Q})$ for $\mathbb{C}\mathrm{P}^2$ is a polynomial algebra
generated by the four-dimensional generator $\mu(\operatorname{pt})$ and the two-dimensional generator
$\mu(\check{\operatorname{H}})$.  In particular, there is no odd-dimensional cohomology of $\mathfrak{B}^*$.
\end{proposition}
Furthermore, it is known \cite[Sec.~2]{Stromme} that for $\mathfrak{M}(-1,k)$ there is a universal complex rank-two vector bundle $\Xi \to \mathbb{C}\mathrm{P}^2
\times \mathfrak{B}^*$ lifting the universal bundle $\mathfrak{P}$ which is defined uniquely up to a linear bundle lifted from $\mathfrak{B}^*$ \cite{Stromme}.
It follows:
\begin{lemma}
\label{chernclasses}
For $\mathbb{C}\mathrm{P}^2$, it follows that as an element of $H^*(\mathbb{C}\mathrm{P}^2 \times \mathfrak{B}^*; \mathbb{Q})$ we have
\begin{equation}
 [ p_1(\mathfrak{P}) ] = \, (1-4k) \, \operatorname{H}^2 - \, 4 \,  \operatorname{H}\wedge \, \mu(\check{\operatorname{H}})
 -4 \, \mu(\operatorname{pt})  \;. \\
\end{equation}
There is a unique choice for $\Xi \to \mathbb{C}\mathrm{P}^2
\times \mathfrak{B}^*$ such that
\begin{equation}
 [ c_2(\Xi) ] = \, k \, \operatorname{H}^2  + \operatorname{H}\wedge \, \mu(\check{\operatorname{H}})  +  \mu(\operatorname{pt}) \;,\qquad
 [ c_1(\Xi) ] = - \, \operatorname{H} \;.
\end{equation}
\end{lemma}
\begin{proof}
Since there is no odd-dimensional cohomology for $\mathfrak{B}^*$ and $\mathbb{C}\mathrm{P}^2$, it follows
from the K\"unneth isomorphism that
\begin{equation*}
\begin{split}
 H^4(\mathbb{C}\mathrm{P}^2 \times \mathfrak{B}^*; \mathbb{Q}) & \cong \bigoplus_{i+j=2} H^{2i}(\mathbb{C}\mathrm{P}^2; \mathbb{Q}) \otimes H^{2j}(\mathfrak{B}^*; \mathbb{Q}) \;.
\end{split}
\end{equation*}
Thus, as elements of $H^*(\mathfrak{B}; \mathbb{Q})$ we have $[\mathfrak{W}^{1,3}]=[ \mathfrak{W}^{3,1}]=0$
in Equation (\ref{universal_curv}). If we denote the exterior derivatives on $\mathbb{C}\mathrm{P}^2$ and
$\mathfrak{B}^*$ by $d$ and $\mathsf{d}$ respectively, it follows that $(d+\mathsf{d}) \, \mathfrak{W}^{\bullet,4-\bullet}=0$.
For dimensional reason, we have $\mathsf{d} \, \mathfrak{W}^{4,0}=0$. However,
$H^0(\mathfrak{B}^*;\mathbb{Q})\cong \mathbb{Q}$ is trivial. Hence, $\mathfrak{W}^{4,0}$ is constant on each connected component of $\mathfrak{B}^*$.
The restriction to $\mathbb{C}\mathrm{P}^2$ yields $[-4 \, \mathfrak{W}^{4,0}] = (1-4k)\, \operatorname{H}^2$. Similarly, since $d \, \mathfrak{W}^{0,4}=0$
it must be constant along $\mathbb{C}\mathrm{P}^2$. From the definition (\ref{observable}), it follows $[-4 \, \mathfrak{W}^{4,0}] = -4 \, \mu(\operatorname{pt})$.
Finally, $[ \mathfrak{W}^{2,2}] = \operatorname{H}\wedge \mu(\check{\operatorname{H}})$ since by Poincar{\'e} duality we have
\begin{equation*}
 \int_{\check{\operatorname{H}}}  \;\mathfrak{W}^{2,2} = \int_{\mathbb{C}\mathrm{P}^2} \operatorname{H}\wedge \, \mathfrak{W}^{2,2} =  \mu(\check{\operatorname{H}}) \;
\int_{\mathbb{C}\mathrm{P}^2} \operatorname{H}\wedge \operatorname{H} =  \mu(\check{\operatorname{H}})\;.
\end{equation*}
Since there exists a universal bundle  $\Xi \to \mathbb{C}\mathrm{P}^2 \times \mathfrak{B}^*$ its Chern classes have to satisfy
$p_1(\mathfrak{P})=c_1^2(\Xi)-4 \, c_2(\Xi)$. Thus, we have $c_1(\Xi)=-\operatorname{H}+ \, 2 \, \delta \, \mu(\check{\operatorname{H}})$
and $c_2(\Xi)= k \, \operatorname{H}^2
+ (1- \delta)  \, \operatorname{H} \wedge \,\mu(\check{\operatorname{H}})
+ \mu(\operatorname{pt})
+ \delta^2 \, \mu(\check{\operatorname{H}})^2$. We normalize the universal bundle by setting $\delta=0$.
\end{proof}
Classes of the de Rham cohomology of $\mathfrak{M}(-1,k)$ are obtained from $H^\bullet(\mathfrak{B}^*)$
by restriction which we will still denote by $\mu(\check{\operatorname{H}})$ and $\mu(\operatorname{pt})$.
To define the integration over the moduli space, a compactification of
the moduli space of anti-selfdual instantons is needed.
Such a compactification was introduced by Donaldson \cite{DonaldsonKronheimer}, but since it is based on Uhlenbeck's results
it is usually called the Uhlenbeck compactification. It was shown that the $\mu$-map in (\ref{observable}) extends over the
compactification $\overline{\mathfrak{M}}(-1,k)$. For $2d_k=4 m + 2 n$,
the {\it Donaldson invariants} are defined as
\begin{eqnarray}
\label{donaldson}
 \pmb{\Phi}_{k,m, n} = \int_{
\overline{\mathfrak{M}}(-1,k)}
 \mu(\operatorname{pt})^{m} \wedge \mu(\check{\operatorname{H}})^{n} \;.
\end{eqnarray}
The map $\pmb{\Phi}$ extends to a linear function from the graded algebra $\operatorname{Sym}_*( H_0(\mathbb{C}\mathrm{P}^2) \oplus H_2(\mathbb{C}\mathrm{P}^2) )$,
where the elements of $H_i(\mathbb{C}\mathrm{P}^2)$ have degree $(4-i)/2$ to $\mathbb{Q}$. It is known that the $\mathrm{SO}(3)$-Donaldson invariants are topological invariants
of $\mathbb{C}\mathrm{P}^2$. For $z \in \operatorname{Sym}_*( H_0(\mathbb{C}\mathrm{P}^2) \oplus H_2(\mathbb{C}\mathrm{P}^2) )$, there is a definition of a geometric
representative $\mathcal{V}(z) \subset \mathfrak{M}(-1,k)$  for $\mu(z)$ \cite{KronheimerMrowka}. If $\overline{\mathcal{V}}(z)$ is the closure of $\mathcal{V}(z)$
in $\overline{\mathfrak{M}}(-1,k)$, then the Donaldson invariants are
\begin{equation}
 \pmb{\Phi}_{k}(z) = \left\lbrace \begin{array}{lcl} \# \Big( \overline{\mathcal{V}}(z) \cap \overline{\mathfrak{M}}(-1,k) \Big) & & \text{if $4m+2n=8(k-1)$} \\ 0 &&
\text{otherwise}\end{array} \right. \;.
\end{equation}
 \begin{definition}
The formal power series
\begin{eqnarray}
\label{generatingfunction}
 \pmb{\mathrm{Z}}(p,S) = \sum_{k \ge1} \; \sum_{m ,n \ge 0} \pmb{\Phi}_{k,m, n} \;
 \; \frac{p^{m}}{m!}\frac{S^{n}}{n!}
\end{eqnarray}
is the generating function for the $\mathrm{SO}(3)$-Donaldson invariants of $\mathbb{C}\mathrm{P}^2$.
\end{definition}
Using the blowup formula for the Donaldson invariants, G\"ottsche \cite{Goettsche,GoettscheNakajimaYoshioka} derived a closed formula expression for $\pmb{\Phi}_{k,m, n}$. His work was based on earlier work with Ellingsrud \cite{EllingsrudGoettsche} and Zagier \cite{GoettscheZagier}, and it extended
the results previously
obtained by Kotschick and Lisca \cite{KotschickLisca} up to an overall sign convention. We state \cite[Thm.~3.5, (1)]{Goettsche} 
in terms of the Jacobi theta-functions $\vartheta_2, \vartheta_3, \vartheta_4$ which are used in theoretical physics, and give a
formulation that is in a convenient form with respect to
$\operatorname{pt}$ and $\check{\operatorname{H}}$.
In this way, we obtain a closed formula expression for $\pmb{\Phi}_{k,m, n}$, which we shall later show equals the
the Moore-Witten prediction based on the $u$-plane integral.
\begin{theorem}[G\"ottsche \cite{Goettsche}]
\label{Goettsche_thm}
Assuming the notation and hypotheses above, then
we have that
the only non-vanishing coefficients
in the generating function in (\ref{generatingfunction})
satisfy
\begin{equation}\label{Goettsche}
\begin{split}
 \pmb{\Phi}_{k,m, 2n} = \sum_{l=0}^n \sum_{j=0}^l & \; \frac{(-1)^{n+j}}{2^{l-3} \, 3^{l}} \; \frac{(2n)!}{(2n-2l)! \; j! \; (l-j)!} \\
 \times & \operatorname{Coeff}_{q^0} \left( \frac{\vartheta_4^8(\tau) \, \left[ \vartheta_2^4(\tau) + \vartheta_3^4(\tau)\right]^{m+j}}
{\left[ \vartheta_2(\tau) \, \vartheta_3(\tau)\right]^{2m+2n+3}}
 \; E^{l-j}_2(\tau) \; F_{2(n-l)}(\tau) \right) \;,
\end{split}
\end{equation}
where $m,n \in \mathbb{N}_0$, $2(k-1)=m+n$, $\operatorname{Coeff}_{q^0}$ is the constant term in a series expansion in $q=\exp{(2\pi i \tau)}$.
The series $F_t(\tau)$ are defined in (\ref{Ft_old}), and
the Jacobi $\vartheta$-functions are
\begin{eqnarray*}
 \vartheta_2(\tau) = 2 \, \Theta_2\left( \frac{\tau}{8}\right) \,, \quad
 \vartheta_3(\tau) = \Theta_3\left( \frac{\tau}{8}\right)  \,, \quad
 \vartheta_4(\tau) = \Theta_4\left( \frac{\tau}{8}\right) \;,
\end{eqnarray*}
where $\Theta_2, \Theta_3, \Theta_4$ are defined in (\ref{theta2})-(\ref{theta4}).
Also, we have that $E_2(\tau)$ is the normalized Eisenstein series
\begin{displaymath}
E_2(\tau):=1-24\sum_{n=1}^{\infty}\sum_{d\mid n}d \; q^n \;.
\end{displaymath}
\end{theorem}
\begin{remark}
Equation (\ref{Goettsche}) uses the original
sign convention of \cite{EllingsrudGoettsche,KotschickLisca} for the Donaldson invariants
in which the Donaldson invariants for $d_k=4$ are all negative.
\end{remark}
\begin{proof}
The following table summarizes the quantities used by G\"ottsche \cite[Thm.~3.5, (1)]{Goettsche} and in this article:
\medskip
\begin{center}
{\footnotesize
\begin{tabular}{|l|l|c|l|l|}
\cline{1-2}\cline{4-5}
 G\"ottsche  & Present Paper & \qquad \qquad & G\"ottsche  & Present Paper \\
\cline{1-2}\cline{4-5}
 $z$ & $S$ &&  $\theta(\tau) $                                        & $\vartheta_4(\tau)$ \\
 $x$ & $p$ &&  $f(\tau)      $                                        & $\frac{1}{2\sqrt{i}} \vartheta_2(\tau) \, \vartheta_3(\tau)$ \\
 $n$ & $2\beta +1, \; \beta \ge 0$ &&  $ \frac{\Delta^2(2\tau)}{\Delta(\tau)\,\Delta(4\tau)}$ & $ - 16 \frac{\vartheta^8_4(\tau)}{\left[ \vartheta_2(\tau) \, \vartheta_3(\tau) \right]^4}$ \\
 $a$ & $2\alpha, \; \alpha \ge \beta+1$ && $ G_2(2\tau) $                                         & $-\frac{1}{24} \, E_2(\tau)$\\
 $\tau   $ & $\frac{\tau-1}{2}$ &&  $ e_3(2\tau) $ & $\frac{1}{12} \left[ \vartheta_2^4(\tau) + \vartheta_3^4(\tau) \right] $\\
 $q  $ & $-q^\frac{1}{2}$ &&  $ \frac{-3i \, e_3(2\tau)}{f(\tau)^2}$  & $\frac{\vartheta_2^4(\tau) + \vartheta_3^4(\tau)}{\left[ \vartheta_2(\tau) \, \vartheta_3(\tau) \right]^2} $ \\
\cline{1-2}\cline{4-5}
\end{tabular}
}
\end{center}
\medskip
\noindent
We use
\begin{eqnarray*}
 && \left( \frac{n}{2} \, \frac{\sqrt{i}}{f(\tau)} \right)^{2(n-l)} \, \left( - \frac{i}{2 \,f(\tau)^2} (2 \, G_2(2\tau)+ e_3(2\tau))\right)^l \\
& = & \frac{(-1)^{n+l}}{2^l \, 3^l} \, (2\beta+1)^{2(n-l)} \, \frac{\left(-E_2(\tau) + \left[\vartheta_2^4(\tau) + \vartheta_3^4(\tau) \right]\right)^l}
{\left[ \vartheta_2(\tau) \, \vartheta_3(\tau) \right]^{2n}} \;.
\end{eqnarray*}
An expansion of the exponential in
\cite[Thm.~3.5, (1)]{Goettsche} then yields (\ref{Goettsche}).
\end{proof}
\subsection{Donaldson theory with monopoles}
Since $w_2(\mathbb{C}\mathrm{P}^2)= \operatorname{H}$, $\mathbb{C}\mathrm{P}^2$ is not a spin manifold.
For the $\operatorname{spin}_{\mathbb{C}}$-structure on $\mathbb{C}\mathrm{P}^2$ with the canonical class $c=3 \, \operatorname{H}$, we
denote by $S_\mathbb{C}^{\pm} \to \mathbb{C}\mathrm{P}^2$ the complex spinor bundles of positive and negative chirality.
In other words, for the spinor decomposition $T_{\mathbb{C}} \mathbb{C}\mathrm{P}^2 = S_\mathbb{C}^{- \; *} \otimes S_\mathbb{C}^{+}$
 we have $c=c_1 ( \det S_\mathbb{C}^{\pm}) = 3 \, \operatorname{H}$.
The vacua of the {\em $\mathrm{SO}(3)$-Donaldson theory with $2 \, N_f$ massless monopoles}
are given by a smooth connection $A \in \mathfrak{A}$ and $N_f$ \emph{complex} sections $\Phi_i
\in C^{\infty}(\mathbb{C}\mathrm{P}^2, S^+_{\mathbb{C}} \otimes \xi)$ such that
\begin{equation}\label{asd_monopole}
 F_A^+ = \sum_{i=0}^{N_f}  \left[\Phi_i \otimes \overline{\Phi}_i\right]_{00}\,, \qquad \slashed{D}_A \Phi_i = 0 \; \;  \text{for} \; 1\le i \le N_f \;,
\end{equation}
where $\slashed{D}_A: C^{\infty}(\mathbb{C}\mathrm{P}^2, S^+_{\mathbb{C}} \otimes  \xi) \to C^{\infty}(\mathbb{C}\mathrm{P}^2, S^-_{\mathbb{C}} \otimes  \xi)$
is the Dirac operator coupled to $\xi$.
Here $\left[\Phi \otimes \overline{\Phi}\right]_{00}$ is the double trace free component in $\operatorname{End}{\xi} \otimes
\operatorname{End}_{\mathbb{C}}{S^+_{\mathbb{C}}}$, and we have used the identification
$\operatorname{End}_{0,\mathbb{C}}{S^+_{\mathbb{C}}} \cong \Lambda^{2+}_{\mathbb{C}\mathrm{P}^2}$. Due to the positive scalar curvature of the Fubini-Study metric, the Weizenbock formula implies that $\ker \slashed{D}_A$ always vanishes whence $\Phi_i=0$.
However, a reducible connection may induce a nonsurjective $\slashed{D}_A$, thus instigating the study of the vector bundle formed by
the cokernel of $\slashed{D}_A$ or equivalently $\ker \slashed{D}^*_A$, called the obstruction bundle. Since the kernel of $\slashed{D}_{A}$ vanishes for the
Fubini-Study metric on $\mathbb{C}\mathrm{P}^2$, the obstruction bundle is the index bundle $\operatorname{Ind} \slashed{D}^* \to \mathfrak{B}^*$ in $K(\mathfrak{B}^*)$.
The fibers of the obstruction bundle are obstruction spaces to the existence of a section for the Dirac operator.
\subsubsection{The index bundle of the twisted Dirac operator.}
In light of the isomorphism between $\mathfrak{M}(c_1,c_2)$ and $\mathfrak{M}(c_1-2r,c_2- r \, c_1 + r^2)$, we will again restrict
ourselves to the rank-two stable vector bundles $\xi \to \mathbb{C}\mathrm{P}^2$ with $c_1(\xi)=-\operatorname{H}, \,
c_2(\xi)=k \, \operatorname{H}^2$ and investigate the index bundle of $\slashed{D}_{A}$ for the
$\operatorname{spin}_{\mathbb{C}}$-structure with the canonical class $c=(2r+1) \, \operatorname{H}$.
The numerical index of $\slashed{D}^*_A$ is computed by the Atiyah-Singer index theorem
\begin{eqnarray*}\label{index}
 \operatorname{ind} \slashed{D}^*_{A} & = & - \int_{\mathbb{C}\mathrm{P}^2} e^{\frac{c}{2}} \; \widehat{\operatorname{A}}(\mathbb{C}\mathrm{P}^2) \; \operatorname{ch}(\xi)
= k - r \;,
\end{eqnarray*}
where $\widehat{\operatorname{A}}(\mathbb{C}\mathrm{P}^2)= 1 - p_1(\mathbb{C}\mathrm{P}^2)/24$ and
$\operatorname{ch}(\xi)= \operatorname{rk}(\xi) + c_1(\xi) + c_1^2(\xi)/2 - c_2(\xi)$.  Furthermore,
the Chern character of the index bundle can be computed from the family index theorem \cite{AtiyahSinger}:
\begin{theorem}
The Chern character of the index bundle is
\begin{eqnarray}\label{family_index}
 \mathrm{ch} \left( \operatorname{Ind} \slashed{D}^* \right)
& = & - \int_{\mathbb{C}\mathrm{P}^2} e^{\frac{c}{2}} \; \widehat{\operatorname{A}}(\mathbb{C}\mathrm{P}^2) \; \operatorname{ch}(\Xi) \;,
\end{eqnarray}
where $\Xi$ is the universal bundle introduced in (\ref{universal_bundle}). $\Xi$ has the total Chern class
$c(\Xi)=1 + c_1(\Xi) +  c_2(\Xi)$.
\end{theorem}
\noindent
It was proved in \cite[Prop.~3.10]{Leness}:
\begin{lemma}
We have
\begin{equation}
 \operatorname{ch}(\Xi) = 2 \, e^{-\frac{\operatorname{H}}{2}} \; \sum_{n=0}^\infty \frac{(-1)^n}{(2n)!}
\left[\left(k - \frac{1}{4}\right) \, \operatorname{H}^2 + \, \operatorname{H}\wedge \, \mu(\check{\operatorname{H}}) + \, \mu(\operatorname{pt})  \right]^n  \;.
\end{equation}
For $c=(2r+1) \, \operatorname{H}$, we have
\begin{equation*}
 \begin{split}
  \operatorname{ch}\left( \operatorname{Ind} \slashed{D}^* \right)_{[2l]} & = - \frac{(-1)^l}{(2l)!} \left[ \left( r^2 - \frac{k + \frac{l}{2}}{2l+1} \right)
 \,  \mu(\operatorname{pt})^l - \frac{l}{2(2l+1)} \, \mu(\check{\operatorname{H}})^2 \wedge \mu(\operatorname{pt})^{l-1}\right] \;, \\
  \operatorname{ch}\left( \operatorname{Ind} \slashed{D}^* \right)_{[2l+1]} & = \frac{(-1)^l\, r}{(2l+1)!} \; \mu(\check{\operatorname{H}})\wedge\,\mu(\operatorname{pt})^l \;.
 \end{split}
\end{equation*}
\end{lemma}
\noindent
It was also proved in \cite[Theorem~1.1]{Leness}:
\begin{theorem}
It follows that
\begin{equation}
 c_{k}(\operatorname{Ind} \slashed{D}^* ) = \sum_{i+2j+2l=k} \; f_{i,2j,2l} \; (2r)^i \; \mu(\check{\operatorname{H}})^{i} \; \wedge \, \mu(\check{\operatorname{H}})^{2j} \wedge \; \mu(\operatorname{pt})^{l} \;,
\end{equation}
where the coefficients $f_{i,2j,2l}$ are given as
\begin{equation*}
 F(x,y,z) = \sum_{i,j,l} \, f_{i,2j,2l} \, x^i \, y^{2j} \, z^{2l} = \exp{\left( \frac{1}{2} \, x \, J_1(z) + \frac{1}{4} \, y^2 \, J_2(z) + J_3(z) \right)} \;,
\end{equation*}
and the functions $J_i(z)$ are given by
\begin{equation*}
 \begin{split}
 J_1(z) & = \frac{\tan^{-1}{(z)}}{z} \;,\\
 J_2(z) & = \frac{z - \tan^{-1}{(z)}}{z^3} \;, \\
 J_3(z) & = - \frac{r^2 - k}{2} \; \ln{(1+z^2)} + \left(k - \frac{1}{4}\right) \; \left( \frac{\tan^{-1}{(z)}}{z} - 1 \right) \;.
 \end{split}
\end{equation*}
\end{theorem}
\begin{lemma}
\label{chernclasses_indexbundle}
For $r=0$, we have $\forall l \in \mathbb{N}_0: \, c_{2l+1}(\operatorname{Ind} \slashed{D}^* )=0$.
\end{lemma}
\begin{proof}
We include an alternative proof as convenience for the reader. We check that $e^{\frac{c}{2}} \; \widehat{\operatorname{A}}(\mathbb{C}\mathrm{P}^2)= 1 + \frac{\operatorname{H}}{2}$.
If $x_1, x_2$ are the formal Chern roots of $\Xi$ it follows that $x_1+x_2 = c_1(\Xi) =-\operatorname{H}$.
We compute for $n>2$:
\begin{equation*}
\operatorname{ch}(\Xi)_{[n]} + \frac{\operatorname{H}}{2} \operatorname{ch}(\Xi)_{[n-1]}
 = \frac{1}{n!} \left( x_1^{n} + x_2^{n} \right) - \frac{1}{2\, (n-1)!} (x_1 + x_2) \, \left( x_1^{n-1} + x_2^{n-1} \right)\;.
\end{equation*}
We use that $x_2=-(x_1+\operatorname{H})$, then expand in terms of $\operatorname{H}$,
and truncate the series to the order $O(\operatorname{H}^3)$. We obtain
\begin{equation*}
\operatorname{ch}(\Xi)_{[n]} + \frac{\operatorname{H}}{2} \operatorname{ch}(\Xi)_{[n-1]}
 = \frac{\left[ 1 + (-1)^n \right]}{n!} \,  x_1^n + \frac{\left[ 1 + (-1)^{n-1} + 2 \, (-1)^n\right]}{2 \, (n-1)!}
\, \operatorname{H} \wedge \, x_1^{n-1} \;.
\end{equation*}
Thus, $\operatorname{ch}(\Xi)_{[n]} + \frac{\operatorname{H}}{2} \operatorname{ch}(\Xi)_{[n-1]}=0$ if $n$ is odd.
Hence, we have $\operatorname{ch}(\operatorname{Ind} \slashed{D}^*)_{[2l+1]}=0$ for all $l$. It follows that
$c_{2l+1}(\operatorname{Ind} \slashed{D}^* )=0$ for all $l$.
\end{proof}
\begin{remark}
Since $\ker \slashed{D}_A$ vanishes and the rank of the index bundle $\operatorname{Ind} \slashed{D}^*$ is $k-r$, the dimension of $\ker \slashed{D}^*_A$
is $k-r$ and the determinant line bundle $\operatorname{DET} \operatorname{Ind}\slashed{D}^* = \Lambda^{k-r} \ker \slashed{D}^* \to \mathfrak{B}^*$ is well defined. It follows
from \cite{AtiyahSinger} that
\begin{equation*}
 c_1\left( \operatorname{DET} \operatorname{Ind}\slashed{D}^* \right) = c_1\left( \operatorname{Ind}\slashed{D}^* \right) \;.
\end{equation*}
Computing $c_1\left( \operatorname{Ind}\slashed{D}^* \right)$ for $r=0$ and the universal bundle $\Xi$ with
\begin{equation*}
\begin{split}
 c_1(\Xi) & =-\operatorname{H}+ \, 2 \, \delta \, \mu(\check{\operatorname{H}}) \;,\\
 c_2(\Xi) & = k \, \operatorname{H}^2  + (1- \delta)  \, \operatorname{H} \wedge \,\mu(\check{\operatorname{H}})+ \mu(\operatorname{pt})
+ \delta^2 \, \mu(\check{\operatorname{H}})^2 \;,
\end{split}
\end{equation*}
it follows
\begin{equation*}
 c_1\left( \operatorname{DET} \operatorname{Ind}\slashed{D}^* \right) = c_1\left( \operatorname{Ind}\slashed{D}^* \right) = k\; \delta \, \mu(\check{\operatorname{H}})\;.
\end{equation*}
Thus, the choice $\delta=0$ in Lemma \ref{chernclasses} is equivalent to the vanishing of $c_1\left( \operatorname{DET} \operatorname{Ind}\slashed{D}^* \right)$, also
called the vanishing of the local anomaly in physics. On the other hand, the normalization used in \cite{Stromme} was different. For $c= (2r+1) \operatorname{H}$ and
$\delta=-1$ we observe that $\operatorname{ind} \slashed{D}^*=k-r^2$ and
\begin{equation*}
 c_1\left( \operatorname{Ind}\slashed{D}^* \right) =  (r^2+ r - k) \, \mu(\check{\operatorname{H}})\;.
\end{equation*}
If $\alpha, \beta, \gamma$ equal $c_1\left( \operatorname{Ind}\slashed{D}^* \right)$ for $r=1,0,-1$ then $\beta=\gamma=-k \,\mu(\check{\operatorname{H}})$ and $
\gamma-\alpha=2\,\mu(\check{\operatorname{H}})$.
\end{remark}
\subsubsection{The vanishing locus of the obstruction.}
For $r=1$, we set
\begin{equation*}
  c_{\operatorname{top}}\left(\operatorname{Ind} \slashed{D}^*\right)   = c_{k-1}(\operatorname{Ind} \slashed{D}^*) \;.
\end{equation*}
For $r=0$, we set
\begin{equation*}
  c_{\operatorname{top}}\left(\operatorname{Ind} \slashed{D}^*\right)   = \left\lbrace \begin{array}{rcl} c_{k-1}(\operatorname{Ind} \slashed{D}^*)
  & & \text{if $k$ is odd} \\ c_{k}(\operatorname{Ind} \slashed{D}^*)  && \text{if $k$ is even} \end{array}\right. \;.
\end{equation*}
We also set
\begin{equation}\label{euler_class}
    c_{\operatorname{top}}\left(\operatorname{Ind} \slashed{D}^{* \; \oplus \; N_f}\right)
= c^{N_f}_{\operatorname{top}}\Big(\operatorname{Ind} \slashed{D}^{*}\Big) \;.
\end{equation}
We denote the Poincar{\'e} dual of the top-dimensional Chern class $c_{\operatorname{top}}$ by
$\mathfrak{J}(N_f,k,r)\subset \mathfrak{M}(-1,k)$. Hence,
$\mathfrak{J}(N_f,k,r)$ is a smooth sub-manifold of $\mathfrak{M}(-1,k)$ of
codimension $N_f(k-r)$ if $r = 1$ and of codimension $2 \, N_f \, \lfloor \frac{k}{2} \rfloor$ if $r=0$.
We interpret $\mathfrak{J}(N_f,k,r)$ as the vanishing locus of the obstruction for the existence of $N_f$ sections
for the Dirac operator. On the other hand, $\mathcal{V}(c_{\operatorname{top}})$ is the geometric representative for the
$\mu$-class used in the definition of the Donaldson invariants.
\subsubsection{The $\mathcal{N}=2$ invariants for $N_f=2,3,4$.}
Using the Euler class, for $N_f= 2, 3$ we define the invariants
\begin{eqnarray}
\label{donaldson2}
 \pmb{\Phi}^{N_f,c,0}_{k,m, n} = \int_{
\overline{\mathfrak{M}}(-1,k) }
 \mu(\operatorname{pt})^{m} \wedge \mu(\check{\operatorname{H}})^{n} \wedge  c_{\operatorname{top}} \left(\operatorname{Ind} \slashed{D}^{* \, \oplus N_f}\right)\;.
\end{eqnarray}
Using the closure $\overline{\mathfrak{J}}(N_f,k)$ of the vanishing locus in the Uhlenbeck compactification $\overline{\mathfrak{M}}(-1,k)$, we define the invariants
\begin{eqnarray}
\label{donaldson2b}
 \pmb{\Phi}^{N_f,c,1}_{k,m, n} = \int_{
\overline{\mathfrak{J}}(N_f,k,r) }
 \mu(\operatorname{pt})^{m} \wedge \mu(\check{\operatorname{H}})^{n}\;.
\end{eqnarray}
For a $z \in \operatorname{Sym}_*( H_0(\mathbb{C}\mathrm{P}^2) \oplus H_2(\mathbb{C}\mathrm{P}^2) )$,
we can write the invariants as the following intersection numbers
\begin{equation}
\label{donaldson3}
\begin{array}{rcrcccl}
 \pmb{\Phi}^{N_f,c,0}_{k}(z) & = \quad & \# \Big(\; \overline{\mathcal{V}}( c_{\operatorname{top}}) & \cap  & \overline{\mathcal{V}}(z) & \cap & \overline{\mathfrak{M}}(-1,k) \Big)  \;,\\
 \pmb{\Phi}^{N_f,c,1}_{k}(z) & = \quad & \# \Big(\; \overline{\mathfrak{J}}(N_f,k,r)          & \cap  & \overline{\mathcal{V}}(z) & \cap & \overline{\mathfrak{M}}(-1,k) \Big)  \;.
\end{array}
\end{equation}
The intersections of $\mathcal{V}(c_{\operatorname{top}})$ and $\mathfrak{J}$ with any
compact cycle in $Z \subset \mathfrak{M}(-1,k)$ will be cobordant. Hence, the intersection numbers will be the same.
If $Z\subset \mathcal{V}(z) \cap \mathfrak{M}(-1,k)$ is a non-compact cycle,
both $\mathcal{V}(c_{\operatorname{top}})\cap Z$ and $\mathfrak{J}\cap Z$ can still be compact.
However, the cobordism need not to be compactly supported. In this case, their closures
in the enlarged Uhlenbeck compactification are no longer cobordant,
and the intersection numbers $\# (\mathcal{V}(c_{\operatorname{top}})\cap Z)$ and $\#(\mathfrak{J}\cap Z)$ differ.
An error term arises which we denote by
\begin{equation}
\label{errorterm1}
 \pmb{\mathcal{E}}^{N_f,c}_{k}(z) = \pmb{\Phi}^{N_f,c,1}_{k}(z) - \pmb{\Phi}^{N_f,c,0}_{k}(z) \;.
\end{equation}
If $\mathcal{V}(z)$
has dimension greater than or equal to four, then the closure of $\mathcal{V}(z)$ in the Uhlenbeck
compactification will intersect the lower strata, and there is a non-vanishing error term. The lower strata of
$\overline{\mathfrak{M}}(-1,k)$ have the form $\mathfrak{M}(-1,k-l) \times \Sigma$ where $\Sigma \subset
\mathrm{Sym}^l(\mathbb{C}\mathrm{P}^2)$ is a smooth stratum. Thus, the error term
is a polynomial expression in the $N_f=0$ Donaldson invariants with smaller instanton numbers $k'<k$ and
\begin{equation}
\label{errorterm2}
 \pmb{\mathcal{E}}^{N_f,c}_{k}(z)  = P\left( \left\lbrace \pmb{\Phi}_{k'}(z') \right\rbrace_{k',z'} \right)\;,
\end{equation}
where $\deg z' = k'$.
\begin{definition}
For $N_f = 2, 3$, the formal power series
\begin{eqnarray}
\label{generatingfunction2}
 \pmb{\mathrm{Z}}^{N_f,c}(p,S) = \sum_{k \ge 1} \sum_{m ,n \ge 0} \pmb{\Phi}^{N_f,c,1}_{k,m,n} \;
 \; \frac{p^{m}}{m!}\frac{S^{n}}{n!}
\end{eqnarray}
is the generating functions for the $\mathrm{SO}(3)$-Donaldson invariants with $2 \, N_f$ massless monopoles of $\mathbb{C}\mathrm{P}^2$
for the $\operatorname{spin}_{\mathbb{C}}$-structure with the canonical class $c=(2r+1) \, \operatorname{H}$.
\end{definition}
\begin{lemma}
\label{odd_coefficients}
For $N_f = 2, 3$, $r=0$, $a=0, 1$, and $n=2l+1$, we have $\pmb{\Phi}^{N_f,c,a}_{k,m, 2l+1} =0$ for all integers $k, l$.
\end{lemma}
\begin{proof}
For $a=0$ it follows from Lemma \ref{chernclasses_indexbundle} that $c_{\operatorname{top}}$ is a form of a degree divisible by four.
The claim then follows from the dimension of the moduli space and the fact that $\mu(\operatorname{pt})$ and
$\mu(\check{\operatorname{H}})$ have degree four and two.
Since the dimension of the moduli space $\mathfrak{M}(-1,k)$ equals $8(k-1)$ it follows $4m + 2n \equiv 0 \mod 4$
whence $n$ is even. The case $a=1$ follows similarly.
\end{proof}
\noindent
We have the following Lemma:
\begin{lemma}
\label{higher_invariants}
Assume $c=\operatorname{H}$, $r=0$. For $k$ even and $m+n+2=k$, we have
\begin{equation}
 \pmb{\Phi}^{2,c,0}_{k,m, 2n} = \sum_{j+l=k}  \;
\left( \sum_{\begin{subarray}{c}j_1 + j_2 = j \\ l_1+l_2=l\end{subarray}} f_{0,2j_1,2l_1} \cdot f_{0,2j_2,2l_2}\right)
\; \pmb{\Phi}_{k,m+l, 2(n+j)} \;.
\end{equation}
For $k$ even and $2m+2n+4=k$, we have
\begin{equation}
 \pmb{\Phi}^{3,c,0}_{k,m, 2n} = \sum_{j+l=3k/2}  \;
\left( \sum_{\begin{subarray}{c}j_1 + j_2 + j_3 = j \\ l_1+l_2 + l_3=l\end{subarray}} f_{0,2j_1,2l_1} \cdot f_{0,2j_2,2l_2}\cdot f_{0,2j_3,2l_3}\right)
\; \pmb{\Phi}_{k,m+l, 2(n+j)} \;.
\end{equation}
\end{lemma}
\begin{proof}
For $N_f=2$ and $k=2\kappa$, we start by computing
\begin{equation*}
 c_{2\kappa}(\operatorname{Ind} \slashed{D}^* ) = \sum_{j+l=\kappa} \; f_{0,2j,2l} \; \mu(\check{\operatorname{H}})^{2j} \wedge \; \mu(\operatorname{pt})^{l} \;.
\end{equation*}
It follows
\begin{equation*}
 c^2_{k}(\operatorname{Ind} \slashed{D}^* ) = \sum_{j+l=k} \;
\left( \sum_{\begin{subarray}{c}j_1 + j_2 = j \\ l_1+l_2=l\end{subarray}} f_{0,2j_1,2l_1} \cdot f_{0,2j_2,2l_2}\right) \;
\mu(\check{\operatorname{H}})^{2j} \wedge \; \mu(\operatorname{pt})^{l} \;.
\end{equation*}
The case $N_f=3$ follows similarly.
\end{proof}
\begin{definition}
For $N_f =4$, the formal power series
\begin{eqnarray}
\label{generatingfunction4}
 \pmb{\mathrm{Z}}^{4,c}(q) =  \sum_{k \ge 1} \, \pmb{\Phi}^{4,c,1}_{k,0,0} \;
 \; q^{\frac{k}{2}}
\end{eqnarray}
is the generating function for the $\mathrm{SO}(3)$-Donaldson invariants with $2N_f=8$ massless monopoles of $\mathbb{C}\mathrm{P}^2$
for the $\operatorname{spin}_{\mathbb{C}}$-structure with the canonical class $c=(2r+1) \, \operatorname{H}$.
\end{definition}
\begin{remark}
For $c=3\, \operatorname{H}$ and $N_f =4$, the Euler class in Equation (\ref{euler_class}) is a top-dimensional form, i.e., a form of
degree $8(k-1)$ on the moduli space. The first four coefficients of (\ref{generatingfunction4}) for $r=1$ were determined in \cite{Gorodenzev}.
\end{remark}
\begin{lemma}
For $c=\operatorname{H}$, $r=0$, we have $\pmb{\Phi}^{4,c,1}_{2l,0,0} =0$
for all integers $l$.
\end{lemma}
\begin{proof}
For $N_f =4$ and $c=\operatorname{H}$, the Euler class in Equation (\ref{euler_class}) is a form
on $\mathfrak{B}^*$ of degree $8k$ when $k$ is even and $8(k-1)$ when $k$ is odd. Hence, $c_{\operatorname{top}}=0$ when $k$ is even and
the top class is restricted to the moduli space of dimension $8(k-1)$. The statement follows.
\end{proof}
\noindent
\subsection{The $\mathcal{N}=4$ partition function}
We will later show that the $u$-plane integral determines a partition function (\ref{generatingfunction4}) for the $\mathcal{N}=2$ topological $SO(3)$-theory on
$\mathbb{C}\mathrm{P}^2$. The partition function satisfies a holomorphic anomaly equation analogous to the partition function for the $\mathcal{N}=4$
topological $SO(3)$-theory on $\mathbb{C}\mathrm{P}^2$ discussed in
\cite[Sec.~4.2]{VafaWitten}. Vafa and Witten showed that the holomorphic part of the partition function $\pmb{\mathcal{Z}}$
is (c.f. \cite[Equation (4.19)]{VafaWitten})
\begin{eqnarray}
\label{generatingfunction_N4}
 \pmb{\mathcal{Z}}^+(\tau) = \frac{1}{q^{\frac{1}{4}}} \sum_{k \ge 1} \, \chi\Big(\overline{\mathfrak{M}}(-1,k)\Big) \;q^{k} \;,
\end{eqnarray}
where $\chi( \overline{\mathfrak{M}}(-1,k))$ is the Euler characteristic of the moduli space $\overline{\mathfrak{M}}(-1,k)$.
The partition function $\pmb{\mathcal{Z}}(\tau)$ satisfies the holomorphic anomaly
equation
\begin{eqnarray}
  \frac{d}{d \bar{\tau}} \; \pmb{\mathcal{Z}}\left(\tau \right) = \frac{3}{16\pi i} \; \; \dfrac{\overline{\vartheta_2(2\tau)}}{\im \tau^{\frac{3}{2}}}\;.
\end{eqnarray}
They conclude that $\pmb{\mathcal{Z}}(\tau) = \pmb{\mathcal{Z}}^+(\tau) +  \pmb{\mathcal{Z}}^-(\tau)$ is a mock modular form of weight $\frac{3}{2}$:
\begin{equation*}
\begin{split}
  \pmb{\mathcal{Z}}^+\left(\tau \right) & = \frac{1}{q^\frac{1}{4}} \, \sum_{\alpha \ge 1} 3\, \mathcal{H}_{4\alpha-1} \, q^{\alpha} \;,\\
  \pmb{\mathcal{Z}}^-\left(\tau \right) & = \frac{3}{4\sqrt{\pi}} \sum_{l \ge 0} (2l+1) \; \Gamma\left( - \frac{1}{2} , \pi \, (2l+1)^2 \, \im\tau \right)
\, q^{-\frac{\left(2n + 1\right)^2}{4}} \;,
\end{split}
\end{equation*}
where $\Gamma(-1/2,t)$ is defined in (\ref{Gamma}) and $\mathcal{H}_{4\alpha-1}$ are the Hurwitz class numbers.
The comparison with Equation (\ref{generatingfunction_N4}) implies $\chi( \overline{\mathfrak{M}}(-1,k)) = 3\, \mathcal{H}_{4k-1}$.
This formula was proved by Klyachko \cite[Theorem 1.1]{Klyachko}.
It is well known that the generating function for the stable rank-two sheaves is equal to the product of the generating function (\ref{generatingfunction_N4})
for the stable rank-two vector bundles and the contribution from the boundary. Yoshioka \cite[Theorem 0.4]{Yoshioka} determined that the contribution from the boundary
is equal to $\eta^{-2\chi}\left(\tau\right)$.
He also determined directly a formula for $\pmb{\mathcal{Z}}^+(\tau)/\eta^{2\chi}\left(\tau \right)$ \cite[Theorem 0.1]{Yoshioka}.
The first terms in the $q$-series expansion are
\begin{equation}
\label{generatingfunction_N4_2}
 \frac{\pmb{\mathcal{Z}}^+(\tau)}{\eta^6\left(\tau \right)}
= \frac{1}{q^{\frac{1}{2}}} \, \Big( q + 9 \, q^2 + 48 \, q^3 + 203 \, q^4 + 729 \, q^5 + 2346 \, q^6 +
6918 \, q^7 + \dots \Big) \;.
\end{equation}
The expansion
\begin{equation*}
\begin{split}
 \frac{\pmb{\mathcal{Z}}^+ \left( \tau \right)}{\eta^6\left( \tau \right)} = \frac{1}{q^\frac{1}{2}} \Big[ & \; 3 \, \mathcal{H}_3 \, q
+ \left(18 \, \mathcal{H}_3 + 3 \mathcal{H}_7 \right) \, q^2
 + \left( 81 \, \mathcal{H}_3 + 18 \, \mathcal{H}_7 + 3 \, \mathcal{H}_{11}\right) \, q^3 \\
 + & \left( 294 \, \mathcal{H}_3 + 81 \, \mathcal{H}_7 + 18 \, \mathcal{H}_{11} + 3 \, \mathcal{H}_{15} \right) \, q^4 + \dots \Big]
\end{split}
\end{equation*}
coincides with Equation (\ref{generatingfunction_N4_2}).
\section{The Seiberg-Witten curves}\label{Section3}
In this section we define and classify the Seiberg-Witten curves needed to define the regularized
$u$-plane integral. A {\it Seiberg-Witten curve} is a rational elliptic surface
with a section and an analytical marking whose fibers are of Kodaira-type $I_n$ and $I_n^*$ with $0 \le n \le 4$.
Physics predicts that the Seiberg-Witten curve enodes the algebraic topology of the moduli spaces of $\mathrm{SO}(3)$-Donaldson theory.
\subsection{Weierstrass elliptic fibrations.}
An elliptic curve $E$ in the Weierstrass form can be written as
\begin{eqnarray}
 \label{jacobian} y^2 = 4 x^3 - g_2 \, x - g_3 \;,
\end{eqnarray}
where $g_2$ and $g_3$ are numbers such that the discriminant $\Delta
=g_2^3- 27 g_3$ does not vanish. In homogeneous coordinates $[X:Y:W]$,
(\ref{jacobian}) becomes
\begin{eqnarray*}
  W Y^2 = 4 X^3 -  g_2 \, X W^2 - g_3 \, W^3.
\end{eqnarray*}
One can check that the point $P$ with the coordinates $[0:1:0]$ is always a
smooth point of the curve. We consider $P$ the base point of the elliptic
curve and the origin of the group law on  $E$. The two types of singularities that can occur as
Weierstrass cubics are a rational curve with a node, which appears
when the discriminant vanishes and $g_2 ,g_3 \not = 0$, or a cusp when $g_2=g_3=0$.
Next, we look at a family of cubic curves over $\mathbb{C}\mathrm{P}^1$.
The family is parametrized by the base space $\mathbb{C}\mathrm{P}^1$ and a line bundle
$N \to \mathbb{C}\mathrm{P}^1$. The quantities
$g_2$ and $g_3$ are promoted to global sections of
$N^{\otimes 4}$ and $N^{\otimes 6}$ respectively; the
discriminant becomes a section of $N^{\otimes 12}$. If the
sections are generic enough so that they do not always lie in the
discriminant locus, we obtain a Weierstrass fibration $\pi: Z \to \mathbb{C}\mathrm{P}^1$ with section.
Each fiber comes equipped with the base point $P$,
which defines a section $\sigma$ of the elliptic fibration which does not pass through the nodes or cusps. The bundle
$N$ is the conormal bundle of the section $\sigma$. We will always assume $N=\mathcal{O}_{\mathbb{C}\mathrm{P}^1}(-1)$.
We will also assume that in the coordinate chart $[u:1] \in \mathbb{C}\mathrm{P}^1$, the discriminant $\Delta$ is a polynomial of degree
$N_f+2$ in $u$, and $g_2$ and $g_3$ are polynomials in $u$ of degree at most $2$ and $3$ respectively. The space of all such Weierstrass elliptic surfaces has $N_f+1$ moduli.
To see this, first consider the case where $N_f=4$. From the seven parameters defining $g_2$ and $g_3$, two can be eliminated by scaling and a shift in the $u$-plane. Furthermore,
we can arrange the coefficient of $g_2$ of degree two and the coefficient of $g_3$ of degree three to be the modular invariants of an elliptic curve
with periods $1$ and $\tau_0$. The remaining four coefficients can be expressed in terms of four complex parameters \cite{SeibergWitten2}.
In physics, they are usually denoted $m_1, \dots, m_{N_f}$, and called the
{\it masses of the hypermultiplets}.
We will explain their geometric meaning below (c.f. Fact \ref{residues}). For $0\le N_f \le 3$, we obtain $4-N_f$ additional constraints
from the requirement that the discriminant has degree $N_f+2$.
A non-trivial elliptic fibration has to develop singular fibers;
the classification of the singular fibers is part of Kodaira's classification theorem of all possible singular fibers
of an elliptic fibration (cf. \cite{Miranda2}). For generic values of the masses, the polynomial $\Delta$ has only $N_f+2$ simple zeros
for $|u|< \infty$, where the elliptic fibration develops a node
(i.e. a singular fiber of Kodaira type $I_1$).
For special values for $m_1, \dots, m_{N_f}$, several singular fibers of Kodaira type $I_1$ can coalesce and form singular fibers of Kodaira type $I_k$ with $k \ge 2$,
where the discriminant has a zero of order $k$.
The second chart over the base space is $[1:v] \in \mathbb{C}\mathrm{P}^1$. The
intersection of the two charts is given by $u=1/v$ with $v \not=0$.
The Weierstrass coordinates transform according to $x \mapsto v^2 x$ and $y \mapsto v^3 y$;
since $g_2$ and $g_3$ are sections of $N^4$ and $N^6$ respectively,
they transform according to $g_{2} \mapsto v^4 \, g_{2}$ and $g_{3} \mapsto v^6 \, g_{3}$.
The discriminant $\Delta \mapsto v^{12} \, \Delta$ becomes a polynomial in $v$ of degree $10-N_f$.
From Kodaira's classification theorem it follows that the singular fiber $E_\infty$ over $u=\infty \; (v=0)$ is a cusp, a singular fiber of Kodaira type $I_{4-N_f}^*$.
In physics, a Weierstrass fibration of the type described is called a
{\em Seiberg-Witten curve}.
$Z$ has surface singular points whenever all partial derivatives in $u, x, y$ simultaneously vanish. Singular fibers of Kodaira type $I_1$ do not give rise to surface
singularities, whereas all singular fibers of Kodaira type $I_n$,
with $n \ge 2$, and $I^*_{n}$, with $n\ge 0$, do.
It is known \cite[Sec.~4.6]{Miranda2} that a Weierstrass fibration is rational
(i.e. birational to $\mathbb{C}\mathrm{P}^2$),
if $g_2$ and $g_3$ are polynomials in $u$ of degree at most $4$ and $6$ respectively.
The minimal  resolution $\widehat{Z}$ of $Z$ is the blow-up of $\mathbb{C}\mathrm{P}^2$ in nine points, and therefore has Picard number $10$. The section
$\sigma$ uses up one dimension, and so the number of components of any singular fiber is at most nine since the components are
always independent in the Neron-Severi group. However, not all configurations of singular fibers exist.
Conversely, by contracting every component of the fiber which
does not meet $\sigma$, we obtain the normal surface $Z$.
A complete list of the possible configurations of singular fibers of the rational elliptic surfaces with sections that appear as Seiberg-Witten curves will be presented
in Lemma~\ref{sw-curves}. We also list the constraints on the moduli,
which when substituted into the Weierstrass presentation
in \cite{SeibergWitten2}, realize the configuration of singular fibers and hence prove their existence.
\subsection{The Mordell-Weil Lattice.}
Let us fix a generic smooth fiber $E$ of $\pi: Z \to \mathbb{C}\mathrm{P}^1$. Let $\operatorname{MW}(\pi)$
be the group of sections of $\pi$ which can be naturally identified with the rational points of $E$ with the origin given by $\sigma$.
The group $\operatorname{MW}(\pi)$ is equipped with a natural bilinear height pairing \cite{OguisoShioda}.
The Picard group of $\widehat{Z}$ is an unimodular lattice with signature $(1,9)$. The section $\sigma$ and the smooth fiber class
$E$ generate a two dimensional lattice that splits off the Picard lattice as $\mathbb{I} \oplus -\mathbb{I}$.
The orthogonal complement $\langle \sigma, E \rangle ^\bot$ is isomorphic to the even negative definite rank $8$
unimodular lattice $E_8$ \cite{Miranda1}. Each reducible fiber consists of a certain number of components that form a sublattice in the
Picard lattice. In particular, a reducible fiber of type $I_k$ and $I_k^*$ generate a root lattice of type $A_{k-1}$ and $D_{k+4}$ respectively.
We will denote the direct sum of the sublattices of the Picard group generated by the components of the fibers not meeting the section $\sigma$ by $\mathbf{T}$.
Thus, $\mathbf{T}$ must have an embedding into the $E_8$ lattice.
It turns out that $\mathbf{T}^\bot$ is isomorphic to the dual lattice of the free part of the Mordell-Weil group
$\left[ \operatorname{MW}(\pi) / \operatorname{MW}(\pi)_{\operatorname{tor}} \right]^*$.
Oguiso and Shioda \cite{OguisoShioda} proved a complete structure theorem of the Mordell-Weil group including the torsion subgroup.
In the table of Lemma~\ref{sw-curves}, we can therefore present $\mathbf{T}, \operatorname{MW}(\pi)$ as well as their ranks, and the classification
number assigned by Oguiso and Shioda.
\begin{lemma}\label{sw-curves}
The following is a complete list of the singular fibers of the rational elliptic surfaces with sections that appear as the
Seiberg-Witten curves (up to permutation of the masses):

\medskip
{\footnotesize
\begin{tabular}{|c|c|ccc|ll|l|l|l|}
 \hline
   $N_f$ & No & $r$ & $\operatorname{rk} \mathbf{T}$ & $\# \operatorname{tor}$
 & $\mathbf{T}$ & $\operatorname{MW}(\pi)$ & $E_\infty$ & $E_{u,\operatorname{sing}}$ & mass constraints\\
 \hline
  4 &  9 &  4  & 4 & 1 &  $D_4                       $& $D_4^*$ & $I^*_0$
& $6 I_1 $     & generic     \\
  4 & 18 &  3  & 5 & 1 &  $D_4 \oplus A_1            $& $A_1^{* \oplus 3}                               $ & $I^*_0$
& $I_2, 4 I_1 $& $m_1 = m_2$ \\
  4 & 32 &  2  & 6 & 1 &  $D_4 \oplus A_2            $& $A^*_1 \oplus \langle \frac{1}{6} \rangle       $ & $I^*_0$
& $I_3, 3 I_1 $& $m_1 = m_2 = m_3$\\
  4 & 34 &  2  & 6 & 2 &  $D_4 \oplus A_1^{\oplus 2} $& $A_1^{*\oplus 2} \oplus \mathbb{Z}_2            $ & $I^*_0$
& $2 I_2,2 I_1 $& $m_1=m_2=0$ \\
  4 & 54 &  1  & 7 & 2 &  $D_4 \oplus A_3            $& $\langle \frac{1}{4} \rangle \oplus \mathbb{Z}_2$ & $I^*_0$
& $I_4, 2 I_1 $& $m_1=m_2=m_3=0$ \\
  4 & 57 &  1  & 7 & 4 &  $D_4 \oplus A_1^{\oplus 3} $& $A_1^{*} \oplus (\mathbb{Z}_2)^2                $ & $I^*_0$
& $3 I_2 $ & $m_1=m_2 , m_3=m_4=0$\\
 \hline
\end{tabular}

\medskip
\begin{tabular}{|c|c|ccc|ll|l|l|l|}
 \hline
  $N_f$ & No & $r$ & $\operatorname{rk} \mathbf{T}$ & $\# \operatorname{tor}$
& $\mathbf{T}$ & $\operatorname{MW}(\pi)$ & $E_\infty$ & $E_{u,\operatorname{sing}}$ & mass constraint \\
 \hline
  3 & 16 & 3 & 5 & 1 & $D_5                       $& $A_3^*                                          $ & $I^*_1$ & $ 5 I_1 $ & generic \\
  3 & 30 & 2 & 6 & 1 & $D_5 \oplus A_1            $& $A_1^* \oplus \langle \frac{1}{4} \rangle       $ & $I^*_1$ & $ I_2, 3 I_1 $ &$m_1=m_2$\\
  3 & 50 & 1 & 7 & 1 & $D_5 \oplus A_2            $& $\langle \frac{1}{12} \rangle                   $ & $I^*_1$ & $ I_3, 2 I_1 $ & $m_1=m_2=m_3$\\
  3 & 52 & 1 & 7 & 2 & $D_5 \oplus A_1^{\oplus 2} $& $\langle \frac{1}{4} \rangle \oplus \mathbb{Z}_2$ & $I^*_1$ & $ 2 I_2, I_1 $ & $m_1=m_2=0$\\
  3 & 72 & 0 & 8 & 4 & $D_5 \oplus A_3            $& $\mathbb{Z}_4                                   $ & $I^*_1$ & $ I_4, I_1 $   & $ m_1=m_2=m_3=0$\\
 \hline
\end{tabular}

\medskip
\begin{tabular}{|c|c|ccc|ll|l|l|l|}
 \hline
  $N_f$ & No & $r$ & $\operatorname{rk} \mathbf{T}$ & $\# \operatorname{tor}$
& $\mathbf{T}$ & $\operatorname{MW}(\pi)$ & $E_\infty$ & $E_{u,\operatorname{sing}}$ & mass constraint \\
 \hline
  2 & 26 & 2 & 6 & 1 & $D_6                       $& $A_1^{*\oplus 2}          $ & $I^*_2$ & $ 4 I_1 $ & generic \\
  2 & 48 & 1 & 7 & 2 & $D_6 \oplus A_1            $& $A_1^* \oplus \mathbb{Z}_2$ & $I^*_2$ & $ I_2, 2 I_1$ & $m_1=m_2$\\
  2 & 71 & 0 & 8 & 4 & $D_6 \oplus A_1^{\oplus 2} $& $(\mathbb{Z}_2)^2         $ & $I^*_2$ & $ 2I_2$ & $m_1=m_2=0$\\
 \hline
\end{tabular}

\medskip
\begin{tabular}{|c|c|ccc|ll|l|l|l|}
 \hline
  $N_f$ & No & $r$ & $\operatorname{rk} \mathbf{T}$ & $\# \operatorname{tor}$
& $\mathbf{T}$ & $\operatorname{MW}(\pi)$ & $E_\infty$ & $E_{u,\operatorname{sing}}$ & mass constraint\\
 \hline
  1 & 46 & 1 & 7 & 1 & $D_7$ & $\langle \frac{1}{4} \rangle$ & $I^*_3$ & $ 3 I_1 $ & - \\
 \hline
\end{tabular}

\medskip
\begin{tabular}{|c|c|ccc|ll|l|l|l|}
 \hline
  $N_f$ & No & $r$ & $\operatorname{rk} \mathbf{T}$ & $\# \operatorname{tor}$
& $\mathbf{T}$ & $\operatorname{MW}(\pi)$ & $E_\infty$ & $E_{u,\operatorname{sing}}$ \\
 \hline
  0 & 64 & 0 & 8 & 2 & $D_8$ & $\mathbb{Z}_2$ & $I^*_4$ & $2 I_1$ \\
 \hline
\end{tabular}
\medskip
}

\noindent
In the table above, we have put $r = \operatorname{rk}\operatorname{MW}(\pi)$, $\# \operatorname{tor}=\left| \operatorname{MW}(\pi)_{\text{tor}} \right|$
and $\langle m \rangle$ for a rank $1$ lattice $\mathbb{Z}x$ with $\langle x,x \rangle=m$.
The embedding of $\mathbf{T}$ into $-E_8$ is unique up to the action of the Weyl group.
The following configurations of singular fibers do not exist on a rational elliptic surface:
\begin{center}
\medskip
{\footnotesize
\begin{tabular}{|c|c|l|c|c|c|l|}
 \cline{1-3} \cline{5-7}
 $N_f$ &  $E_\infty$ & $E_{u,\operatorname{sing}}$ &&
 $N_f$ &  $E_\infty$ & $E_{u,\operatorname{sing}}$ \\
 \cline{1-3} \cline{5-7}
 4  & $I_0^*$ & $I_6$           &\ \ \ \ \ \  \ \ & 2  & $I_2^*$ & $I_4$      \\
 4  & $I_0^*$ & $I_5, I_1$      && 2  & $I_2^*$ & $I_3, I_1$ \\
 \cline{5-7}
 4  & $I_0^*$ & $I_3, I_3$      && 1  & $I_3^*$ & $I_3 $     \\
 4  & $I_0^*$ & $I_4, I_2$      && 1  & $I_3^*$ & $I_2, I_1$ \\
 \cline{5-7}
 4  & $I_0^*$ & $I_3, I_2, I_1$ && 0  & $I_4^*$ & $I_2$ \\
 \cline{1-3} \cline{5-7}
 3  & $I_1^*$ & $I_5$           \\
 3  & $I_1^*$ & $I_3, I_2$      \\
 \cline{1-3}
 \end{tabular}}
\end{center}
\end{lemma}
\begin{proof}
Using the parametrization of Seiberg and Witten \cite{SeibergWitten2}, one determines
the zeros of $g_2, g_3, \Delta$ and their degrees which will determine the singularities of
Kodaira type $I_k$ and $I_k^*$. The obtained configuration of singular fibers
appear in
\cite{OguisoShioda}, and they give $\mathbf{T}$ and $\operatorname{MW}(\pi)$.
The uniqueness of the embedding follows from
the fact that the Seiberg-Witten curves by construction always contain a singular fiber of type $I_k^*$
and \cite[Thm. 3.3]{OguisoShioda}. The impossible configurations are found in a list
produced by Persson \cite{Persson}, and are further explained by Miranda \cite{Miranda2}.
\end{proof}
Let $\mathbf{\Lambda}$ be the set of simple roots of the Mordell-Weil lattice
(i.e. a basis with the property that every vector in the root lattice
is a linear combination of elements of $\mathbf{\Lambda}$ with all coefficients non-negative). For each $\alpha \in \mathbf{\Lambda}$,
we denote the corresponding
section by $s_{\alpha} \in \operatorname{MW}(\pi)$. The sections $s_{\alpha}$ generate $\operatorname{MW}(\pi)$, and we denote the union of all sections
by $\mathfrak{S}=\sum_{\alpha \in \mathbf{\Lambda}} s_{\alpha}$.
\subsection{Analytical marking and Seiberg-Witten one-form.}
We will denote  by $\mathrm{UP}$ the base curve $\mathbb{C}\mathrm{P}^1$ with small open discs around the points
with a singular fiber removed. In physics, $\mathrm{UP}$ is known as
the {\em $u$-plane}. The restriction of $\pi: Z \to \mathbb{C}\mathrm{P}^1$ to $\mathrm{UP}$,
which we will still denote by $Z$, is a smooth four-dimensional manifold with a three-dimensional boundary with multiple components, and
so $\pi: Z \to \mathrm{UP}$ is a proper, surjective, holomorphic map of smooth complex manifolds. We denote by $\Omega_Z$ the canonical
bundle of $Z$. The spaces $H^0(\Omega_Z)$ is the space of global holomorphic two-forms on $Z$.
We fix a smooth fiber $E_u$ of the fibration, and thus have $\dim H^1(E_u)=2$. Since a base point is given in each fiber by the section $\sigma$, we can
choose a symplectic basis $\lbrace A_u, B_u \rbrace$ of the homology $H_1(E_u)$ with respect to the intersection form, called a {\it homological marking}.
We cannot define $A_u, B_u$ globally over $\mathrm{UP}$. The cycles are transformed by monodromies around the points
with singular fibers. However, we can define globally an analytical marking.
An {\it analytical marking} is a choice
of a non-zero holomorphic one-form on the smooth fiber $E_u$. We choose the canonical holomorphic differential
$dx/y \in H^0(\Omega_{E_u})$ in the coordinates of (\ref{jacobian}) where $\Omega_{E_u}$ is the holomorphic cotangent
bundle of $E_u$.  Similarly,  let $\Omega_{\mathrm{UP}}$ be the restriction
of the canonical bundle on $\mathbb{C}\mathrm{P}^1$ to $\mathrm{UP}$. The bundle $\Omega_{Z/\mathrm{UP}} = \Omega_Z \otimes (\pi^* \Omega_{\mathrm{UP}})^{-1}$
restricts to the canonical bundle $\Omega_{E_u}$ on each smooth fiber $E_u$. Thus, given the rational elliptic surface $Z \to \mathrm{UP}$
and the analytic marking, we can associate to it the holomorphic symplectic two-form $\Omega_{SW} = du \wedge dx/y \in H^0(\Omega_Z)$.
The Picard-Fuchs equations on $Z \to \mathrm{UP}$ ask whether there is a meromorphic one-form $\lambda_{SW} \in H^0(\Omega_{Z-\mathfrak{S}/\mathrm{UP}})$,
called the {\em Seiberg-Witten differential}, such that $\Omega_{SW}$ is the derivative of $\lambda_{SW}$:
\begin{definition}\label{periodintegrals}
The Picard-Fuchs equations on $Z \to \mathrm{UP}$ for $\lambda_{SW} \in H^0(\Omega_{Z-\mathfrak{S}/\mathrm{UP}})$ are
\begin{eqnarray}\label{PicardFuchs}
  d \; \int_{A_u} \lambda_{SW}  = \int_{A_u} \Omega_{SW} \;, \qquad  d \int_{B_u} \lambda_{SW}  = \int_{B_u} \Omega_{SW} \;,
\end{eqnarray}
where $d$ is the exterior derivative on $\mathrm{UP}$.
\end{definition}
In general, $\lambda_{SW}$ has poles located at the intersection of sections and fibers, and it is only well-defined as
a section of $\Omega_{Z-\mathfrak{S}/\mathrm{UP}}$. We denote the period integrals of the elliptic fiber by
$2 \, \omega  \, du =  \int_{A_u} \Omega_{SW} $ and $2 \, \omega' \, du = \int_{B_u} \Omega_{SW}$,
and the integrals of $\lambda_{SW}$ by $2 \, \mathbf{a}  =  \int_{A_u} \lambda_{SW}$ and  $2\, \mathbf{a}_D = \int_{B_u} \lambda_{SW}$.
The modular parameter of the elliptic fiber $E_u$ is $\tau = \omega'/\omega \in \mathbb{H}$.
We have that $\omega, \omega'$ are
sections of a holomorphic rank-two vector bundle over $\mathrm{UP}$, called the {\it period bundle}. The vector bundle is equipped with a
flat connection $\nabla: H^0\left(\Omega_{Z-\mathfrak{S}/\mathrm{UP}}\right) \to H^0\left(\Omega_{\mathrm{UP}}\right)
\otimes H^0\left(\Omega_{Z-\mathfrak{S}/\mathrm{UP}}\right)$ which satisfies the equation $\nabla \lambda_{SW} = \Omega_{SW}$, and
is known as the {\it Gauss-Manin connection}.
The holonomy of the connection around the singular fibers determines a local system on $\mathrm{UP}$ and a representation
of the fundamental group $\pi_1(\mathrm{UP})\to \mathrm{SL}_2(\mathbb{Z})$ \cite{Atiyah, Stiller}. The connection has regular singularities at the base points of
the singular fibers \cite{Stiller2, M}.
\begin{remark}
The variable $a$ used in \cite{SeibergWitten1, SeibergWitten2} is related to $\mathbf{a}$ by $\mathbf{a} =2\sqrt{2} \pi \, a$.
\end{remark}
For the Seiberg-Witten curves
listed in Lemma \ref{sw-curves}, the following two facts were proved in \cite{SeibergWitten2} by an explicit computation.
For their interpretation in terms of the Gauss-Manin connection on the period bundle, we refer to Shimizu \cite[Lemma 3.1.5, Rem.~3.2.6]{Shimizu}.
\begin{fact}
The Picard-Fuchs equations (\ref{PicardFuchs}) have a solution $\lambda_{SW} \in H^0(\Omega_{Z-\mathfrak{S}/\mathrm{UP}})$ which is unique up to
a holomorphic one-form on $Z$.
\end{fact}
\begin{fact}\label{residues}
For each $s_\alpha \in \operatorname{MW}(\pi)$ with $\alpha \in \mathbf{\Lambda}$, the residue $m_\alpha=\operatorname{Res}_{s_\alpha}(\lambda_{SW})$
is flat. Up to the permutation of the roots, the residues are the parameters $m_1, \dots, m_{N_f}$ in the description of the rational
elliptic surfaces in Lemma \ref{sw-curves}.
\end{fact}
\section{The $u$-plane integral}\label{Section4}
In this section we define the regularized $u$-plane integral on $\mathbb{C}\mathrm{P}^2$. The $u$-plane integral will depend on
the choice of a Seiberg-Witten curve (i.e. a rational elliptic surface
with a section and an analytical marking).
To compute the $u$-plane integral, the integral has to be renormalized and then extended across the points in the $u$-plane
with singular fibers. For $\mathbb{C}\mathrm{P}^2$, we will show that
the regularized $u$-plane integral receives contributions only from the point $u=\infty$, where the fiber is cuspidal.
The regularized $u$-plane integral computes the generating function of a $\mathcal{N}=2$ supersymmetric, topological $U(1)$-gauge theories on $\mathbb{C}\mathrm{P}^2$.
In physics, the theory is called the \emph{low energy effective field theory}.
We compute the $u$-plane integrals on $\mathbb{C}\mathrm{P}^2$ for the Weierstrass presentations of the cases
$\# \, 54, \, 57, \, 64, \, 71, \, 72$ in Lemma \ref{sw-curves}.
\subsection{Massless Seiberg-Witten curves for $N_f=0,2,3$.}
\label{section_sw_curves}
Here we consider the Seiberg-Witten curves when $N_f=0,2,3$.
\begin{lemma}\label{presentation}
In Lemma \ref{sw-curves}, the cases $\# \, 64, \, 71, \, 72$ are the only elliptic modular surfaces. The Weierstrass presentations
are as follows:
{\footnotesize
\begin{equation*}
\begin{array}{|c|rrrr||lcl|lcl|}
\hline
 \multicolumn{5}{|l||}{N_f =  0,\# 64}
   &
 \multicolumn{6}{|l|}{u^{(0)} \in \mathrm{UP}_{\epsilon}^{(0)}=\mathbb{C}\mathrm{P}^1 - B_{\epsilon}(- 1)   - B_{\epsilon}(1)  - B_{\epsilon}(\infty)}  \\
\hline
 E_{\operatorname{sing}} & I_1 & I_1 & I_4^*  &
   &  g^{(0)}_2   & = & \frac{\left(u^{(0)}\right)^2}{12} - \frac{1}{16} &   \tau^{(0)} & \in & \mathbb{H}/\Gamma_0(4) \\
\cline{1-5}
 u^{(0)}_{\operatorname{sing}}                         & -1  & 1   & \infty &
   &  g^{(0)}_3      & = &  \frac{\left(u^{(0)}\right)^3}{216} -\frac{u^{(0)}}{192}   & u^{(0)}   &  =  & \frac{1}{2} \frac{\vartheta_2^4 + \vartheta_3^4}
{\left[\vartheta_2 \, \vartheta_3\right]^2} \\
\tau^{(0)}_{\operatorname{sing}}                       & 0   & 2   & \infty &
   &  \Delta^{(0)}   & = &  \frac{1}{4096}\left[\left(u^{(0)}\right)^2-1\right] &  \omega^{(0)} & = & \sqrt{2} \, \pi \, \vartheta_2 \; \vartheta_3 \\
\hline
\hline
\multicolumn{5}{|l||}{N_f =  2,\# 71}
 &
 \multicolumn{6}{|l|}{u^{(2)}  \in \mathrm{UP}^{(1)}_\epsilon=\mathbb{C}\mathrm{P}^1 -  B_{\epsilon}(- 1)   - B_{\epsilon}(1)  - B_{\epsilon}(\infty)}  \\
\hline
 E_{\operatorname{sing}} & I_2 & I_2 & I_2^*  &
   &  g_2^{(2)}   & = &  \frac{\left(u^{(2)}\right)^2}{12} + \frac{1}{4} & \tau^{(2)} & = & \frac{\tau^{(0)}}{2} \in \mathbb{H}/\Gamma(2) \\
\cline{1-5}
 u^{(2)}_{\operatorname{sing}}                         & -1  & 1   & \infty &
   &  g_3^{(2)}      & = &  \frac{\left(u^{(2)}\right)^3}{216} - \frac{u^{(2)}}{24} &  u^{(2)}   &  =  & u^{(0)} \\
\tau^{(2)}_{\operatorname{sing}}                 & 0   & 1   & \infty &
   &  \Delta^{(2)}   & = &  \frac{1}{64}\left[\left(u^{(2)}\right)^2-1\right]^2 & \omega^{(2)} & = & \omega^{(0)}\\
\hline
\hline
 \multicolumn{5}{|l||}{N_f =  3,\# 72}
  &
 \multicolumn{6}{|l|}{u^{(3)} \in \mathrm{UP}_{\epsilon}^{(3)}=\mathbb{C}\mathrm{P}^1 - B_{\epsilon}(0)   - B_{\epsilon}(1)  - B_{\epsilon}(\infty) }  \\
\hline
 E_{\operatorname{sing}} & I_4 & I_1 & I_1^*  &
&  g_2^{(3)}   & = &  \frac{\left(u^{(3)}\right)^2}{12} - \frac{5 u^{(3)}}{4} + \frac{11}{16}& \tau^{(3)} & = & -\frac{1}{\tau^{(0)}}\\
 \cline{1-5}
 u^{(3)}_{\operatorname{sing}}                      & -\frac{1}{2}  & \frac{1}{2}   & \infty &
   &  g_3^{(3)}      & = &  \frac{\left(u^{(3)}\right)^3}{216} + \frac{7 \left(u^{(3)}\right)^2}{48} & u^{(3)}   &  =  & - \frac{2}{u^{(0)}-1} - \frac{1}{2}\\
\tau^{(3)}_{\operatorname{sing}}                 & 0   & - \frac{1}{2}   & \infty &
 & &-& \frac{29 u^{(3)}}{96} + \frac{7}{64}  & \omega^{(3)} & = & \frac{\pi}{2} \; \left( \vartheta_3^2 - \vartheta_4^2 \right)  \\
\cline{9-10}
   &&&&&  \Delta^{(3)}   & = &  \multicolumn{3}{l|}{- \frac{1}{512} \left(2 u^{(3)}- 1\right)  \left( 2 u^{(3)}+1\right)^4}&  \\
\hline
\end{array}
\end{equation*}
}
\noindent
The relation between the discriminant and the Dedekind eta-function is given by $\Delta=\pi^{12} \, \eta^{24}(\tau)/\omega^{12}$.
$B_{\epsilon}(u_{\operatorname{sing}})$ is the $\epsilon$-disc
around the point $u_{\operatorname{sing}}$ with singular fiber.
\end{lemma}
\begin{proof}
 Shioda \cite{Shioda} proved that elliptic modular surfaces have maximal Picard number.
 The Seiberg-Witten curves $\# \, 64, 71, 72$ are the only extremal ones. The parametrization in the
 cases $\# \, 64, 71$ was given in \cite{SeibergWitten1, SeibergWitten2}. The cases $\# \, 64$ and $\# \, 71$
 are related by  a two-isogeny. We have obtained the parametrization for $\# 72$ by twisting the fibration
 in $\# \, 64$. The Weierstrass presentations of any two surfaces that are related by a twist
 are related by `starring' one and `unstarring' another singular fiber of the first
 fibration to obtain the second. The Euler number of the singular fibers increase or decrease by six respectively.
 One multiplies $g_2$ and $g_3$ by $(u+1)^2$ and $(u+1)^3$ respectively to obtain the new $g_2$ and $g_3$ up to
 a change of the chart on $\mathbb{C}\mathrm{P}^1$.
\end{proof}
\begin{remark} The rational elliptic surfaces in the cases $\# \, 54, \, 57$ in Lemma~\ref{sw-curves}
are said to be in {\it non-canonical form} since they can be obtained by twisting the elliptic fibrations $\# \, 64, \, 71$
in Lemma~\ref{presentation}. We will discuss the three cases in Lemma~\ref{sw_potential}.
\end{remark}
\begin{remark}
We label the points in the $u$-plane with a singular fiber by an index $s \in \lbrace D, M, \infty \rbrace$. The two points in the $u$-plane with $|u_{\operatorname{sing}}|<\infty$
where the fiber is of Kodaira-type $I_{K_s}$ are called the points where the \emph{dyon} ($u_{\operatorname{sing}}= 1/2 \; \text{or} \;1$) and
\emph{monopole} ($u_{\operatorname{sing}}=-1/2 \; \text{or} \; -1$)
becomes massless. There is a duality transformation $\tau = A_{s}\cdot \tau_s$ with
$A_s \in \mathrm{SL}_2(\mathbb{Z})$ such that approaching the singular point corresponds to
$\im \tau_s \to \infty$. We set $u_{s}(\tau_s) = u( \tau)$ and $A_s \cdot \omega_s(\tau_s) = \omega(\tau)$.
Over $u=\infty$, the fiber is of Kodaira-type $I^*_{4-N_f}$, and approaching the singular point corresponds to
$\im \tau \to \infty$. We set $\tau_{\infty}=\tau$ and $u_\infty = u$.
 If $N_f$ is fixed, we will often suppress the superindex $(N_F)$.
\end{remark}
\begin{remark}
The $\epsilon$-disc $B_{\epsilon}(u_{\operatorname{sing}})$ around the point $u_{\operatorname{sing}}$ with a singular fiber of the Kodaira-type $I_{K_s}$ is constructed as follows:
the boundary of the $\epsilon$-disc corresponds to the line in the $\tau_s$-coordinate with
$\re \, \tau_s \in [0 ; K_s]$ and $\im{\tau_s} = - \ln{\epsilon}$.
The $\epsilon$-disc around around the point $u_{\operatorname{sing}}=\infty$ with a singular fiber of the Kodaira-type $I^*_{4-N_f}$ is constructed as follows:
the boundary of the $\epsilon$-disc corresponds to the line in the $\tau$-coordinate with $\re \, \tau \in [0 ; 4-N_f]$ and
$\im{\tau} = -\ln{\epsilon}$.
\end{remark}
\begin{lemma}\label{choice_of_tau}
For $|q_s|<1$, there are normally convergent power series expansions,
where $\alpha_{s,m}, \beta_{s,m} \in \mathbb{Q}$, such that
\begin{equation}
\begin{split}
 u_{s}(\tau_s) & = \sum_{m \ge 0} \alpha_{s,m} \; q_s^{\frac{m}{K_s}} \;, \\
 \omega^2_s(\tau_s) & =  \sum_{m \ge 0} \beta_{s,m} \; q_s^{\frac{m}{K_s}} \;.
\end{split}
\end{equation}
For $|q|<1$, there are normally convergent power series expansions,
where $\alpha_{m}, \beta_{m} \in \mathbb{Q}$, such that
\begin{equation}
\begin{split}
 \left[ q^{\frac{1}{4-N_f}} \; u(\tau) \right] & =\sum_{m \ge 0} \alpha_{m} \; q^{\frac{m}{4-N_f}} \;, \\
 \left[ q^{-\frac{1}{4-N_f}} \;\omega^2(\tau)\right]  & =  \sum_{m \ge 0} \beta_{m} \; q^{\frac{m}{4-N_f}} \;.
\end{split}
\end{equation}
\end{lemma}
\begin{proof}
The proof is a direct consequence of the explicit Weierstrass presentations
in Lemma \ref{presentation}.
\end{proof}
\begin{lemma}\label{sw_potential}
In the cases $\# \, 64, \, 71, \, 72$, the Seiberg-Witten meromorphic one-form is
\begin{equation}\label{sw-form}
\lambda_{SW} = \frac{(N_f+2) \, u}{3}  \, \frac{dx}{y} - \frac{\delta_{3,N_f}}{2}\, \frac{dx}{y}
 - 4 (4-N_f) \frac{x \, dx}{y} \;,
\end{equation}
and it has vanishing residues. The integrals over the A-cycle and B-cycle are
\begin{subequations}
\label{a}
\begin{align}
\mathbf{a} & = \frac{(N_f+2) \, u}{3}  \, \omega - \frac{\delta_{3,N_f}}{2}\, \omega
 + (4-N_f) \frac{\pi^2 \, E_2(\tau)}{3 \, \omega} \;,\\
\mathbf{a}_D & = \frac{(N_f+2) \, u}{3}  \, \omega' - \frac{\delta_{3,N_f}}{2}\, \omega'
 + (4-N_f) \frac{\pi^2 \, \tau \, E_2(\tau)}{3 \, \omega} - 2 \pi i \, \frac{4-N_f}{\omega} \;,
\end{align}
\end{subequations}
such that $d\mathbf{a}/du = \omega$, $d\mathbf{a}_D/du = \omega'$, and
\begin{equation}
 \omega' \; \mathbf{a} - \omega \; \mathbf{a}_D = 2 \pi i \, (4-N_f) \;.
\end{equation}
\end{lemma}
\begin{proof}
Equation (\ref{sw-form}) was proved in \cite{AlvarezGaumeMarinoZamora}. The case $N_f=2$ follows from the case $N_f=0$
using the equation $E_2(\tau/2)=2 \, E_2(\tau) - u \, \omega^2/\pi^2$. Equations (\ref{a}) follows from the fact that the integrals of the
meromorphic form $x \,dx/y$ over the A-cycle and B-cycle equal $\int_{A_u} x \, dx/y =-2 \, \underline{\eta}$ and
$\int_{B_u} x \, dx/y =-2 [ \tau \, \underline{\eta} - \pi i/(2\omega)]$ with $\underline{\eta}=\pi^2 \, E_2(\tau)/(12\omega)$.
\end{proof}
\begin{remark}
In the cases $\# \, 64,\, 71, \, 72$, the Picard-Fuchs equations (\ref{PicardFuchs}) are
second order hypergeometric ordinary differential equations with three regular singularities \cite{Stiller}.
\end{remark}
\subsection{The Renormalization Group}
To describe the $u$-plane integral, we have to introduce the \emph{two-observable} $\widehat{T}$ of the low energy effective field theory for
the families $\# \, 64, \, 71, \, 72$ described in Lemma \ref{presentation}. To do so, we rescale (\ref{jacobian}) describing the elliptic fiber over $u$ by $\Lambda^6$ and
introduce $x_{\Lambda} = \Lambda^2 x$, $y_{\Lambda} = \Lambda^3 y$, $g_{\Lambda,2}= \Lambda^4 g_2$, $g_{\Lambda,3}= \Lambda^6 g_3$.
Then, we observe that (\ref{jacobian}) is still satisfied in the rescaled quantities.
We also rescale the coordinate of the base space $u_{\Lambda} = \Lambda^2 \, u$. In physics, rescaling in this manner is called the action of
\emph{the renormalization group}. One of the important properties of the renormalization group is that the holomorphic symplectic two-form is not invariant but rescaled
$\Omega_{\Lambda,SW}=\Lambda \, \Omega_{SW}$. We set $\mathbf{a}_{\Lambda} = \Lambda \, \mathbf{a}$ and consider
$u_{\Lambda}$ a function of $\mathbf{a}_{\Lambda}$.
\begin{definition}
The contact term of the renormalization group flow is
\begin{equation}
 T = \left. \frac{1}{4-N_f} \left( \Lambda \; \frac{\partial u_{\Lambda}}{\partial \Lambda}
 \right)_{\mathbf{a}_{\Lambda}=\operatorname{const}} \right|_{\Lambda=1} \;.
\end{equation}
\end{definition}
\begin{lemma} We have that
$T =  - \frac{\pi^2 \, E_2(\tau)}{3 \, \omega^2} + \, \frac{u}{3} \, + \frac{\delta_{3, N_f}}{2}$\;.
\end{lemma}
\begin{proof}
From the equation
$ d\mathbf{a}_{\Lambda} =  \frac{\mathbf{a}_{\Lambda}}{\Lambda} \, d\Lambda + \Lambda \, \frac{d\mathbf{a}}{du} \, du
 =   \frac{\mathbf{a}_{\Lambda}}{\Lambda} \, d\Lambda + \Lambda^2 \, \omega_{\Lambda} \, du$ it follows
\begin{eqnarray*}
 \Lambda \, du_{\Lambda} & =& 2 \, \Lambda^2 \, u \, d\Lambda + \Lambda^3 \, du =
  \left( 2 \, u_{\Lambda} - \frac{\mathbf{a}_{\Lambda}}{\omega_{\Lambda}} \right) \, d\Lambda + \frac{\Lambda}{\omega_{\Lambda}}
 \, d\mathbf{a}_{\Lambda} \;,
\end{eqnarray*}
whence
\begin{eqnarray}
\label{T}
 T & =  \frac{1}{4-N_f} \left( 2 \, u - \frac{\mathbf{a}}{\omega} \right) - \delta_{3,N_f}  =  - \frac{\pi^2 \, E_2(\tau)}{3 \, \omega^2} + \, \frac{u}{3} \, + \frac{\delta_{3, N_f}}{2} \;.
\end{eqnarray}
\end{proof}
\begin{lemma}\label{T_asymptotics}
We have that
$T$ satisfies $T = O(u^{-1})$ as $u \to \infty$.
\end{lemma}
\begin{proof}
On the elliptic curve $E_u$ in (\ref{jacobian}), we have the following relation between $g_2, g_3$ and the Eisenstein series $E_4(\tau), E_6(\tau)$ of weight
four and six
\begin{eqnarray*}
  \frac{g_3}{g_2} = \dfrac{\frac{\pi^6}{216 \, \omega^6} \, E_6(\tau)}{\frac{\pi^4}{12 \, \omega^4} E_4(\tau)} = \frac{\pi^2}{18 \, \omega^2} \, \frac{E_6(\tau)}{E_4(\tau)} \;,
\end{eqnarray*}
whence
\begin{eqnarray}\label{T_exp}
 T  & = - 6 \; \frac{g_3}{g_2} \; \frac{E_2(\tau) \, E_4(\tau)}{E_6(\tau)} + \, \frac{u}{3} \, + \frac{\delta_{3, N_f}}{2} \;.
\end{eqnarray}
It is easy to check from the presentation in Lemma~\ref{choice_of_tau} that near the singular fiber of type $I_{4-N_f}^*$, the behavior of $q=\exp(2\pi i \tau)$
is as follows:
\begin{eqnarray}\label{u-expansion}
  q = \frac{1}{u^{4-N_f}} \; \Big( c_0 + O(u^{-1}) \Big)
  \qquad (u \to \infty)\;.
\end{eqnarray}
It follows from the presentation in Lemma~\ref{presentation} that near the singular fiber of type $I_{4-N_f}^*$, we have the following behavior
\begin{eqnarray*}
  \frac{g_3}{g_2} = \frac{u}{18} + \frac{31}{12} \, \delta_{3,N_f} + O(u^{-1}) \qquad (u \to \infty)\;.
\end{eqnarray*}
For $N_f=0,1,2$ the claim follows. Using the equation
$\Delta = \pi^{12} \, \eta^{24}(\tau)/\omega^{12}$, we compute
$c_0=-1/16$ in (\ref{u-expansion}) for $N_f=3$ and case  $\# 72$.
It follows that
\begin{eqnarray*}
  \frac{E_2(\tau) \, E_4(\tau)}{E_6(\tau)} & = &\frac{\left(1-\frac{24 \, c_0}{u}\right) \left(1+\frac{240 \, c_0}{u}\right)}{ \left(1-\frac{504 \, c_0}{u}\right)} + O(u^{-2}) \\
 & = &  1 - \frac{45}{u} + O(u^{-2}) \qquad (u \to \infty)\;,
\end{eqnarray*}
and the claim follows from the expansion of (\ref{T_exp}) as $u \to \infty$.
\end{proof}
\begin{remark}
The asymptotic behavior proved in Lemma \ref{T_asymptotics},
together with the holomorphicity in $u$, was originally used by
Moore and Witten \cite{MooreWitten} to define $T$. It was later
shown
 \cite{LNS, MarinoMoore} that $T$ can be described as the contact term of the renormalization group flow which allows a definition in terms of the geometry of the elliptic surface alone.
However, for $N_f=3$ these two methods of determining $T$ generally only agree up to a shift in $u$. The variable $u^{(3)}$ defined in Lemma~\ref{presentation} is chosen in a way
that makes the two methods agree. The variable $u^{(3)}$ relates to the variable $u^{(3)}_{SW}$ used by Seiberg and Witten \cite{SeibergWitten2} by $u^{(3)} = u^{(3)}_{SW} -1/2$.
\end{remark}
\noindent
Finally, to obtain an expression which is a modular invariant in $\tau$,
we make the following definition.
\begin{definition}
For the elliptic surfaces in the cases $\# \, 64, \, 71, \, 72$ in Lemma \ref{presentation},
the two-observable of the low energy effective field theory is
\begin{equation}\label{That}
 \widehat{T} = - \frac{\pi^2 \, \widehat{E}_2(\tau)}{3 \, \omega^2} + \, \frac{u}{3} \, + \frac{\delta_{3, N_f}}{2} \;,
\end{equation}
where $\widehat{E}_2(\tau)$ is the non-holomorphic weight 2 Eisenstein series defined by
\begin{displaymath}
\widehat{E}_2(\tau)= E_2(\tau) - \frac{3}{\pi \, \im\tau} \;.
\end{displaymath}
\end{definition}
\subsection{The $u$-plane integral.}
We briefly describe some features of the low energy effective field theory on $\mathbb{C}\mathrm{P}^2$.
In particular, we outline the physics argument on how the generating function is
computed using the language of determinant line bundles so that the
reader can see how the idea arose in physics that elliptic integrals
compute the Donaldson invariants.
The bosonic fields of the classical field theory are a connection $\mathfrak{a}$ on a line bundle $L \to \mathbb{C}\mathrm{P}^2$ with the curvature $F_{\mathfrak{a}}$ and a complex
valued function $\varphi \in C^{\infty}(\mathbb{C}\mathrm{P}^2;\mathbb{C})$.
A family of bosonic gauge-invariant actions is parametrized by a complex coupling coupling $\tau \in \mathbb{H}$  which is the modular parameter
$\tau$ of the elliptic fiber $E_u$ of a Seiberg-Witten curve $Z \to \mathrm{UP}_{\epsilon}$ in Lemma \ref{presentation}, and it is given by
\begin{equation}
\label{LEFT}
\begin{split}
\mathbf{S} & =
 \frac{i \bar{\tau}}{16\pi} \int_{\mathbb{C}\mathrm{P}^2} F_\mathfrak{a}^+\wedge F_\mathfrak{a}^+  + \frac{i \tau}{16\pi} \int_{\mathbb{C}\mathrm{P}^2} F_\mathfrak{a}^- \wedge F_\mathfrak{a}^- \\
& +  \frac{i}{4}  \int_{\mathbb{C}\mathrm{P}^2} F_\mathfrak{a} \wedge \alpha + \frac{\im{\tau}}{8\pi} \int_{\mathbb{C}\mathrm{P}^2} d\varphi \wedge * d\bar{\varphi} \;,
\end{split}
\end{equation}
where $\alpha \in H^2(\mathbb{C}\mathrm{P}^2)$ is an integer class whose mod-two reduction is $w_2(\mathbb{C}\mathrm{P}^2)$.
The critical points of the action are self-dual connections
(i.e. $F_\mathfrak{a}^-=0$), and a constant $\varphi$.
The field theory is connected to the $\mathrm{SO}(3)$-Donaldson theory in
Section~\ref{Sec_Donaldson} through two global constraints: the topological condition $c_1(L)[\check{\operatorname{H}}]=
(2 \,c_2(\xi)+c^2_1(\xi))[\mathbb{C}\mathrm{P}^2]$, and the constraint $\varphi=\mathbf{a}$, where $\mathbf{a}$ is the integral of $\lambda_{SW}$ over the A-cycle in Lemma~\ref{sw_potential}.
The first condition implies $c_1(L) = [ -\frac{1}{2\pi} \, F_{\mathfrak{a}}] = (2k+1) \, \operatorname{H}$, the latter will be further explained below.
If we view $\tau$ as a function of $\mathbf{a}$, each critical point of the action in (\ref{LEFT}) is labeled by the discrete topological data $k$ and the
continuous modulus $\mathbf{a}$. The action $\mathbf{S}$ at a critical point $(k,\mathbf{a})$ is
\begin{equation*}
 \mathbf{S}^{(0)} = i \pi \bar{\tau} \; \left(k + \frac{1}{2}\right)^2 - i \pi \left(k+\frac{1}{2}\right)  \;.
\end{equation*}
The path integral for the supersymmetric extension of the action in (\ref{LEFT}) can be defined with mathematical rigor by the stationary phase approximation \cite{Witten1}.
The quadratic approximation of the action around a critical point is the Hessian of the supersymmetric extension of $\mathbf{S}$. It determines
a free field theory in the collected variations of the bosonic fields $\widetilde{\Phi}$ and fermionic fields
$\widetilde{\Psi}$ of the form
\begin{eqnarray*}
  \mathbf{S}^{(2)} & = & \int_{\mathbb{C}\mathrm{P}^2} \operatorname{vol} \quad \Big( \left< \widetilde{\Phi} , \Delta_{(k,\mathbf{a})} \widetilde{\Phi} \right>
+ \left< \widetilde{\Psi} , \, \slashed{D}_{(k,\mathbf{a})} \widetilde{\Psi}\right> \Big) \;,
\end{eqnarray*}
where $\Delta_{(k,\mathbf{a})}$ is a family of second-order elliptic operators,
and $\slashed{D}_{(k,\mathbf{a})}$ is a family of real skew-symmetric first-order operators depending
on the moduli. We point out that the operators $\Delta_{(k,\mathbf{a})}$ and
$\slashed{D}_{(k,\mathbf{a})}$ describe the Hessian of the supersymmetric action that
describes the gauge theory, but also the coupling to gravity.
The functional integration over the fluctuations (i.e. the
coordinates of the normal bundle at the critical points),
is an infinite-dimensional Gaussian integral. We define the semi-classical path
integral\footnote{The Gaussian integral for the free field theory is actually well-defined and agrees with the ratio of the Pfaffians in
the definition we give.}
\begin{eqnarray*}
 \int \left[ \mathcal{D} \widetilde{\Phi} \; \mathcal{D} \widetilde{\Psi}\right] \; e^{-\mathbf{S}^{(2)}}
\;\; \text{to be} \; \;  \frac{\operatorname{pfaff}' \slashed{D}_{(k,\mathbf{a})}}{\sqrt{\sideset{}{'}\det \Delta_{(k,\mathbf{a})}}} \;.
\end{eqnarray*}
The above expression is a section of the
product of the Pfaffian line bundles over the moduli space. To proceed further, one has to be able to integrate
this section over the moduli space. To do so, the line bundle needs to be flat. In physics, this is called the {\it vanishing of the local anomaly}.
Since the topology of the moduli space is not trivial, the section can still have holonomy around non-trivial loops in the moduli space.
To make sure that the bundle is globally trivial as well, the holonomy around all curves needs to be the identity. In physics, this is called
the {\it  vanishing of the global anomaly}.
Finally, it is not enough to have a trivial line bundle.
We also need a canonical trivialization (i.e. a trivializing
section). All these points are in fact satisfied \cite{M2}, and we regard the ratio of determinants as a function on the moduli space.
We then carry out the integral over the continuous moduli, and
we carry out the sum over the discrete moduli to obtain the semi-classical approximation
of the \emph{partition function}
\begin{equation}
\label{partition_function}
\begin{split}
&  \sum_{k \in \mathbb{Z}} \quad \int  d\mathbf{a}\wedge d\bar{\mathbf{a}} \quad
  e^{-\mathbf{S}^{(0)}} \quad   \int \left[ \mathcal{D} \widetilde{\Phi} \; \mathcal{D} \widetilde{\Psi}\right] \; e^{-\mathbf{S}^{(2)}} \\
= \; & \sum_{k \in \mathbb{Z}} \quad \int d\mathbf{a}\wedge d\bar{\mathbf{a}} \quad
  e^{-\mathbf{S}^{(0)}} \;  \frac{\operatorname{pfaff}' \slashed{D}_{k,\mathbf{a}}}{\sqrt{\sideset{}{'}\det \Delta_{k,\mathbf{a}}}} \;.
\end{split}
\end{equation}
In the case which is relevant for the definition of the $u$-plane integral, physical considerations guarantee that the semi-classical approximation
is in fact exact \cite{MooreWitten, Witten1}, and (\ref{partition_function}) is the full quantum partition function. Moore and Witten \cite{MooreWitten} also showed that
\begin{equation}\label{determinant}
 \frac{\operatorname{pfaff}' \slashed{D}_{k,\mathbf{a}}}{\sqrt{\sideset{}{'}\det \Delta_{k,\mathbf{a}}}}
 = - \frac{4i}{(2\pi)^{\frac{3}{2}}} \quad \frac{1}{\sqrt{\im\tau}} \frac{d\bar{\tau}}{d\bar{\mathbf{a}}} \, \left(\int_{\check{\operatorname{H}}} F_{\mathfrak{a}} \right) \quad \frac{\Delta^{\frac{1}{8}}}{\omega^{\frac{3}{2}}} \;,
\end{equation}
where the last term on the right side in (\ref{determinant}) describes the coupling to gravity on $\mathbb{C}\mathrm{P}^2$, which had already been determined by Witten \cite{Witten2}.
We not only want to determine the partition function, but we also wish to
construct the generating function with the insertion of the observables
that correspond to the observables in (\ref{observable}) for the $\mathrm{SO}(3)$-Donaldson theory on $\mathbb{C}\mathrm{P}^2$.
Moore and Witten \cite{MooreWitten} determined that $\mu(\operatorname{pt}) \mapsto 2u$ as a consequence of the choice $\varphi=\mathbf{a}$, and
$\mu(\check{\operatorname{H}}) \mapsto \widehat{T}$, where $\widehat{T}$ was defined in (\ref{That}).
Thus, we obtain  the \emph{generating function}
\begin{equation}\label{generating_function0}
\begin{split}
 & \sum_{k \in \mathbb{Z}} \quad \int  d\mathbf{a}\wedge d\bar{\mathbf{a}} \quad
  e^{-\mathbf{S}^{(0)}} \quad   \int \left[ \mathcal{D} \widetilde{\Phi} \; \mathcal{D} \widetilde{\Psi}\right] \; e^{-\mathbf{S}^{(2)}} \;\; e^{2 \, u \, p + S^2 \, \widehat{T}}\\
= \; & \sum_{k \in \mathbb{Z}} \quad \int d\mathbf{a}\wedge d\bar{\mathbf{a}} \quad
  e^{-\mathbf{S}^{(0)}} \;  \frac{\operatorname{pfaff}' \slashed{D}_{k,\mathbf{a}}}{\sqrt{\sideset{}{'}\det \Delta_{k,\mathbf{a}}}} \;\; e^{2 \, u \, p + S^2 \, \widehat{T}}\;.
\end{split}
\end{equation}
We use the formula
\begin{equation}\label{eta3}
 \sum_{k\in \mathbb{Z}} e^{i \pi \left( k + \frac{1}{2} \right)} \; \left(k + \frac{1}{2} \right) \; e^{- i \pi \bar{\tau} \, \left(k+ \frac{1}{2} \right)^2} = i \; \overline{\eta^3(\tau)} \;,
\end{equation}
the substitution rule, (\ref{determinant}) and (\ref{partition_function}) to make the following definition:
\begin{definition}\label{integral_wo_sings}
The generating function of the massless $N_f$ low energy effective field theory on $\mathbb{C}\mathrm{P}^2$ is
\begin{equation}\label{uplane1}
 \widetilde{\pmb{\mathrm{Z}}}^{N_f}_{\mathrm{UP}_{\epsilon}}(p,S)
 =  - \frac{8}{\sqrt{2\pi}} \; \int_{\mathrm{UP}_{\epsilon}}  \frac{du\wedge d\bar{u}}{\sqrt{\im\tau}} \; \frac{d\bar{\tau}}{d\bar{u}} \; \frac{\Delta^{\frac{1}{8}}}{\omega^{\frac{1}{2}}} \;
  e^{2 \, u \, p + S^2 \, \widehat{T}} \; \; \overline{\eta^3(\tau)} \;.
\end{equation}
By $\mathrm{UP}_\epsilon$, we denote the base curve $\mathbb{C}\mathrm{P}^1$ with small open $\epsilon$-discs around the points
with the singular fiber of Kodaira-type $I_{K_s}$ and the point $u=\infty$ with the singular fiber of Kodaira-type $I^*_{4-N_f}$,
removed. Also, we respectively let
$\tau, \Delta, \omega, \widehat{T}$ denote the modular parameter, the discriminant, the
period integral over the A-cycle, and the two-observable of the regular fiber defined for the families $\# \, \,64, \, 71, \, 72$  by
(\ref{That}).
\end{definition}
\begin{lemma}
The generating function in (\ref{uplane1}) is well-defined.
\end{lemma}
\begin{proof}
For the integral to be well-defined we have to check two facts:
the integral in (\ref{uplane1}) is locally well-defined since the poles of the integrand only appear at the singular locus
which have been removed. Secondly, the integrand is written in terms which have monodromy in $\tau$. The $\tau$-dependent part of the integral is
\begin{equation}\label{tau_part}
 \frac{\pi^{\frac{3}{2}}}{\sqrt{\im\tau}} \; \frac{d\bar{\tau}}{d\bar{u}} \; \frac{|\eta(\tau)|^6}{\omega^2} \; e^{S^2 \, \widehat{T}} \;.
\end{equation}
It is easy to show that the term is invariant under any
modular transformation $\Gamma \subset \mathrm{SL}_2(\mathbb{Z})$ which
leaves $u$ invariant.
\end{proof}
In summary, formula (\ref{LEFT}) is the bosonic part of the low energy effective field theory. A physics argument says that computing
the non-abelian Donaldson invariants can be reduced to a saddle-point computation in an abelian theory given by the supersymmetric extension of
(\ref{LEFT}) where the coupling constant is taken to be the $\tau$-parameter of a Seiberg-Witten family of curves.
\subsection{The Mock modular form.}
To compute the generating function,
we integrate by parts using the
``holomorphic part'' of a weight 1/2 harmonic Maass form.
The weight 1/2 harmonic Maass form is given in Theorem~\ref{QasMu}, and the
holomorphic part of a weight 1/2 harmonic Maass form is defined in Equation~(\ref{holomorphic_part}).
We set $Q(\tau) = Q^+(\tau) + Q^-(\tau)$, where
$Q^+(\tau)$ with $q=\exp{(2\pi i \tau)}$ is the holomorphic part of this Maass form.
The holomorphic part has a series expansion of the form
\begin{equation}\label{Qplus}
 Q^+\left(\tau\right)  = \frac{1}{q^{\frac{1}{8}}} \; \sum_{\alpha \ge 0} H_\alpha \; q^{\frac{\alpha}{2}}
 = \frac{1}{q^{\frac{1}{8}}}\left( 1 + 28 \, q^{\frac{1}{2}} + 39 \, q + 196 \, q^{\frac{3}{2}} + 161 \, q^2 + \dots \right) \;.
\end{equation}
The ``non-holomorphic part'' $Q^{-}$ is
\begin{equation}\label{Qminus}
\begin{split}
 Q^{-}\left(\tau\right) & = \frac{1}{q^{\frac{1}{8}}} \sum_{\alpha \ge 0} H_{-\alpha} \, q^{-\alpha}  =
- \frac{2i}{\sqrt{\pi}} \; \sum_{l \ge 0} (-1)^l \; \Gamma\left( \frac{1}{2} , \pi \, \frac{(2l+1)^2}{2} \, \im\tau \right) \; q^{-\frac{(2l+1)^2}{8}}
\end{split}
\end{equation}
where $\Gamma(1/2,t)$ is defined in (\ref{Gamma}).
The non-holomorphic part
$Q^{-}$ has an exponential decay since
\begin{equation*}
  \Gamma\left( \frac{1}{2} , t\right) = \frac{e^{-t}}{\sqrt{t}} \left( 1 + O(t^{-1}) \right)
\qquad (t \to \infty) \;.
\end{equation*}
Setting $\zeta_{2k}=\exp(\pi i/k)$ it follows that $\zeta_8^2 \, Q^{+}(\tau+2) = Q^+(\tau)$ and $\zeta_8 \, Q^{-}(\tau+1)=Q^-(\tau)$.
Thanks to Theorem~\ref{QasMu}, we know that
\begin{equation}\label{DE}
 \sqrt{2} i \; \frac{d}{d\bar{\tau}} \; Q\left(\tau\right) = \frac{1}{\sqrt{\im \tau}} \, \overline{\eta^3(\tau)} \;.
\end{equation}
We have the following analog of \cite[(9.18)]{MooreWitten} for $Q$:
\begin{lemma}\label{EkQ}
The function
\begin{equation*}
 \mathcal{E}^k_{\frac{1}{2}} \left[ Q \right] = \sum_{j=0}^k (-1)^j \; \binom{k}{j} \; \frac{\Gamma\left(\frac{1}{2}\right)}{\Gamma\left(\frac{1}{2}+j\right)} \; 2^{2j} \; 3^j
 \; E_2^{k-j}(\tau) \; \left(q \, \frac{d}{dq} \right)^j Q\left(\tau\right)
\end{equation*}
is modular of weight $2k+1/2$ for $\Gamma(2) \cap \Gamma_0(4)$ and satisfies
\begin{equation}
 \sqrt{2} i \; \frac{d}{d\bar{\tau}}\; \mathcal{E}^k_{\frac{1}{2}} \left[ Q \right]  = \frac{1}{\sqrt{\im \tau}} \, \widehat{E}_2^k(\tau) \; \overline{\eta^3(\tau)} \;.
\end{equation}
\end{lemma}
\begin{proof}
The proof of the modularity is very similar to the proof of the Lemma in \cite[Sec.~9.1]{MooreWitten}.
Let $f$ be modular of weight $(s, 0)$. Cohen's operator acting on $f$ is given by
\begin{equation*}
 \mathcal{F}_{r}[f] = \sum_{j=0}^r \, \binom{r}{j}  \frac{\Gamma(s+r)}{\Gamma(s+j)} \, \left( \frac{1}{2iy} \right)^{r-j} \, \left( \frac{d}{d\tau} \right)^j f \;.
\end{equation*}
If $f$ is modular of weight $(s, 0)$ for $\Gamma(2) \cap \Gamma_0(4)$ then $\mathcal{F}_r[f]$ is modular of weight $(2r + s, 0)$ for $\Gamma(2) \cap \Gamma_0(4)$.
One can easily check that
\begin{equation*}
 \mathcal{E}^k_{\frac{1}{2}} \left[ Q \right] = \sum_{m=0}^k  \; \binom{k}{m} \; \frac{\Gamma\left(\frac{1}{2}\right)}{\Gamma\left(\frac{1}{2}+k-m\right)} \; \left( \frac{6i}{\pi} \right)^{k-m}
 \; \widehat{E}_2^{m}(\tau) \; \mathcal{F}_{k-m}[Q] \;.
\end{equation*}
Moreover, we apply $d/d\bar{\tau}$ to $Q(\tau)$ in each summand. We use the fact that
$\overline{\eta^3(\tau)}$ does not depend on $\tau$. Using  $\left( \frac{d}{d\tau} \right)^j
\frac{1}{\sqrt{\im\tau}} = (-1)^j  \frac{\Gamma\left(\frac{1}{2}+j\right)}{\Gamma\left(\frac{1}{2}\right)}\; \frac{\sqrt{2i}}{(\tau-\bar{\tau})^{j+\frac{1}{2}}}
$ the claim follows.
\end{proof}
\noindent
For the different coordinate patches around the cusps of Kodaira-type $I^*_{4-N_f}$ and $I_{K_s}$
in Lemma \ref{choice_of_tau} we need to define the transforms of $Q$ which we will use to integrate by parts in equation (\ref{generating_fct1}).
If $F(\tau)$ is a solution to (\ref{DE}), then taking  $\tau \mapsto -1/\tau$ and $\tau \mapsto \tau+n$, implies
that $\widetilde{F}(\tau)=(-1/\tau)/\sqrt{-i\tau}$ and $\zeta_8^n F(\tau+n)$ are solutions as well. Their non-holomorphic parts equal
$Q^-(\tau)$ in equation (\ref{Qminus}) since the transforms satisfy the same differential equation in $\partial_{\bar{\tau}}$.
\begin{lemma}\label{ModularForm}
The function $Q(\tau)$ satisfies
\begin{eqnarray*}
\frac{1}{\sqrt{-i \, \tau}} \; Q\left(- \frac{1}{\tau}\right) & = & \frac{1}{\sqrt{-i \, \tau}} \; \zeta_8^2 \; Q\left(2- \frac{1}{\tau}\right) \;.
\end{eqnarray*}
We have the following $q$-expansion for the holomorphic part
\begin{equation*}
\left[ \frac{1}{\sqrt{-i \, \tau}} \; Q\left(- \frac{1}{\tau}\right) \right]^+
  =   \frac{1}{q^{\frac{1}{8}}}\left( \frac{5}{2} + \frac{111}{2} \, q + \frac{413}{2} \, q^2 + 819 \, q^3 + \frac{4407}{2} \, q^4 + \dots \right) \;.
\end{equation*}
\end{lemma}
\begin{proof}
It is proved in Theorem~\ref{QasMu} that
\begin{eqnarray*}
 Q(\tau) & = & - \frac{7}{2} \, A_{3,8}(\tau) + \frac{3}{2} \, A_{7,8}(\tau) - \frac{1}{2} \, B(\tau) + 4 \, M(\tau) \\
 & = &  \frac{1}{q^{\frac{1}{8}}} \left( 1 + 28 \, q^{\frac{1}{2}} + 39 \, q + 196 \, q^{\frac{3}{2}} + 161 \, q^2 + \dots \right) \;,
\end{eqnarray*}
where the functions $A_{3,8}, A_{7,8}, B, M$ are defined in (\ref{A-functions}), (\ref{B-functions}), and (\ref{M-function}).
We have the following transformation properties
\begin{eqnarray*}
 \zeta_8 \, A_{3,8}(\tau+1) & = & - A_{3,8}(\tau) \;,\\
 \zeta_8 \, B(\tau+1) & = & B(\tau) \;, \\
 \zeta_8 \, A_{7,8}(\tau+1) & = & A_{7,8}(\tau) \;,\\
 \zeta_8 \, M(\tau+1) & = & M(\tau) \;.
\end{eqnarray*}
The expansions of the modular forms in $-1/\tau$ can be computed easily. For example, we have
\begin{eqnarray*}
\frac{1}{\sqrt{-i \,\tau}} \; A_{3,8}\left( - \frac{1}{\tau} \right) =   - \frac{1}{2}   \frac{\eta^8\left(\frac{\tau}{2}\right)}{\eta^7(\tau)}
  =  \frac{1}{q^{\frac{1}{8}}} \left( - \frac{1}{2} +4 \, q^{\frac{1}{2}} - \frac{27}{2}  \, q + 28 \, q^{\frac{3}{2}} \dots \right)\;,
\end{eqnarray*}
and similar equations for $B(\tau), A_{7,8}(\tau)$. For the expansion of $M$ we use Theorem \ref{ZwegersThm} to show that
\begin{eqnarray*}
 \frac{1}{\sqrt{-i \,\tau}} \; M\left( - \frac{1}{\tau} \right)
 & = & - \frac{\sqrt{i}}{4 \, q^\frac{1}{32}} \; \widehat{\mu}\left( \frac{1}{2}, \frac{1}{4}-\frac{\tau}{8}; \frac{\tau}{4}\right)
 - \frac{1}{4 \, \sqrt{i} \, q^\frac{1}{32}} \; \widehat{\mu}\left( \frac{1}{2},  \frac{3}{4}-\frac{\tau}{8}; \frac{\tau}{4}\right)\\
&= & \frac{1}{q^{\frac{1}{8}}} \left(  \frac{1}{2} \, q^{\frac{1}{4}} - q^{\frac{1}{2}} + 2 \, q^{\frac{3}{4}} - 3 \, q + \dots \right) \;.
\end{eqnarray*}
This proves the series expansion. The second series expansion
then follows from $\zeta_8^2 \, Q(\tau+2) = Q(\tau)$\;.
\end{proof}
\begin{definition}
We set $Q^{(0)}_{\infty} (\tau) =  Q^{(2)}_{\infty} (\tau) = Q(\tau)$, and
\begin{eqnarray*}
  Q^{(0)}_{D} (\tau) = Q^{(0)}_{M} (\tau) & = & \frac{1}{\sqrt{-i \, \tau}} \; Q\left(- \frac{1}{\tau}\right) \;.
\end{eqnarray*}
Similarly, we define
\begin{eqnarray*}
 Q^{(2)}_{M} (\tau) =  \frac{1}{\sqrt{-i \, \tau}} \; Q\left(- \frac{1}{\tau}\right) \;, \quad
 Q^{(2)}_{D} (\tau) =  \frac{1}{\sqrt{-i \, \tau}} \; \zeta_8 \; Q\left(1- \frac{1}{\tau}\right) \;.
\end{eqnarray*}
For $N_f=3$, we set
\begin{equation*}\label{Qplus_3}
  Q^{(3)}_{\infty}\left(\tau\right) = Q^{(0)}_M\left(\tau\right) \,, \quad
  Q^{(3)}_M\left(\tau\right) =  Q^{(0)}_{\infty}\left(\tau\right) \,, \quad
  Q^{(3)}_D\left(\tau\right) =  Q^{(0)}_{D}\left(\tau\right) \,.
\end{equation*}
\end{definition}
\noindent
\begin{lemma}
\label{transfo_Z_rho}
For $\pmb{\mathrm{Z}}(\tau)=\eta^3(\tau) \, Q(\tau)$ and
$\rho(\tau)=\sqrt{2} \; \eta^2(2\tau)/\eta^2(\tau)$, we have that
\begin{equation*}
 \pmb{\mathrm{Z}}(\tau) = \pmb{\mathrm{Z}}(\tau+1) + 14 \;\eta^4(\tau) \, \rho^4(\tau)  \;.
\end{equation*}
We have the following transformation properties under modular transformations
\begin{equation*}
\begin{array}{rcrcrcr}
 \pmb{\mathrm{Z}}(\tau+2) & = & \pmb{\mathrm{Z}}(\tau) \;,&\qquad&
 \pmb{\mathrm{Z}}\left(-\frac{1}{\tau}+ 1\right) &=& \tau^2 \;\;\; \pmb{\mathrm{Z}}\left(\tau+1\right) \;.
\end{array}
\end{equation*}
The expression
\begin{equation}\label{invariant_combination}
  \frac{1}{\eta^4(\tau)} \, \sum_{k=0}^3 (-1)^k \; \pmb{\mathrm{Z}}(\tau+k) = 28 \, \rho^4(\tau)=   \frac{4\, q^{\frac{3}{8}}}{\eta(\tau)}  \; \sum_{m=0}^{\infty} \, H_{2m+1} \, q^{m}
\end{equation}
is a weight $0$ modular form.
\end{lemma}
\begin{proof}
The equation
\begin{equation*}
 \pmb{\mathrm{Z}}(\tau) = \pmb{\mathrm{Z}}(\tau+1) + 14 \;\eta^4(\tau) \, \rho^4(\tau)
\end{equation*}
follows from the proof of Lemma \ref{ModularForm}. Using the remark following Theorem \ref{QasMu} it follows that
\begin{eqnarray*}
\sum_{k=0}^3 \, (-1)^k \; \pmb{\mathrm{Z}}(\tau+k)   & = & -14 \, A_{3,8}(\tau) \, \eta^3(\tau) = 112 \, \frac{\eta(2\tau)^8}{\eta(\tau)^4}\;.
\end{eqnarray*}
\noindent
To prove that the indicated $q$-series is modular, we apply
the following principle (see the Section \ref{HarmonicMaassForms}): A harmonic
weak Maass form whose non-holomorphic part is zero, is a weakly
holomorphic modular form. To this end we now compute the contribution
to the Fourier expansion of the ``period integrals" to see
that indeed they vanish thereby giving us a weight 1/2 weakly holomorphic
modular form, which after division by $\eta(\tau)$ then becomes a weakly
holomorphic modular form of weight $0$.
To this end, notice
that (\ref{Qminus}) and (\ref{Qplus})
imply that $Q^{\pm}(\tau+2)=\zeta_8^{-2} \,Q^{\pm}(\tau)$.
Since $\eta^3(\tau+2)=\zeta_8^2 \eta^3(\tau)$, we have $\pmb{\mathrm{Z}}(\tau+2)=\pmb{\mathrm{Z}}(\tau)$. To see that
the non-holomorphic parts
vanish, we check that for $Q^-(\tau)$ in Equation (\ref{Qminus}) we have $Q^{-}(\tau+1)=\zeta_8^{-1} \,Q^{-}(\tau)$
whence $\pmb{\mathrm{Z}}^-(\tau+1)=\pmb{\mathrm{Z}}^-(\tau)$.
Thus, we obtain
\begin{eqnarray*}
2 \, \left[  \pmb{\mathrm{Z}}(\tau) - \pmb{\mathrm{Z}}(\tau+1) \right]
& = & 2 \, \left[  \pmb{\mathrm{Z}}^+(\tau) - \pmb{\mathrm{Z}}^+(\tau+1) \right]\\
 & = &  2 \,  \eta^3(\tau) \left[ Q^+(\tau)  - \zeta_8 \, Q^+(\tau+1) \right] \\
 & = &  \frac{4 \, \eta^3(\tau)}{q^{\frac{1}{8}}} \, \sum_{m=0}^{\infty} H_{2m+1} \; q^{m+\frac{1}{2}} \;.
\end{eqnarray*}
The S-duality equation for $\pmb{\mathrm{Z}}\left(\tau\right)$
follows from the identity
\begin{equation}
\begin{split}
\frac{1}{\sqrt{-i \, \tau}} \; \zeta_8 \; Q\left(1- \frac{1}{\tau}\right)
 & =  - \zeta_8 \; Q(\tau+1) \;,
\end{split}
\end{equation}
which can be proved using the results in Section 6 and Theorem 7.2.
\end{proof}
\subsection{The regularization procedure.}
We now extend the integral by taking the limit $\epsilon \to 0$ to include the contributions from the singularities.
To do so we will have to \emph{regularize} the $u$-plane integral. We will denote the result by
\begin{equation}\label{uplane2}
\begin{split}
  \pmb{\mathrm{Z}}^{N_f}_{\mathrm{UP}}(p,S)
 & =  - \frac{8}{\sqrt{2\pi}} \; \int_{\mathbb{C}\mathrm{P}^1}^{\operatorname{reg}}
\frac{du\wedge d\bar{u}}{\sqrt{\im\tau}} \; \frac{d\bar{\tau}}{d\bar{u}} \; \frac{\Delta^{\frac{1}{8}}}{\omega^{\frac{1}{2}}} \;
  e^{2 \, u \, p + S^2 \, \widehat{T}} \; \overline{\eta^3(\tau)} \;.
\end{split}
\end{equation}
\begin{remark}
Equation (\ref{uplane2}) is the analog of \cite[(9.1)]{MooreWitten} for the $\mathrm{SO}(3)$-Donaldson theory and follows the
same normalization.
\end{remark}
\noindent
We expand the exponential in the generating function in (\ref{uplane1}) as follows:
\begin{equation}\label{generating_fct1}
\begin{split}
 \widetilde{\pmb{\mathrm{Z}}}^{N_f}_{\mathrm{UP}_{\epsilon}}(p,S) = & \sum_{m, n \ge 0} \frac{p^m}{m!}  \frac{S^{2n}}{(2n)!} \;\; e^{S^2 \, \frac{\delta_{3,N_f}}{2}}
\quad
\sum_{k=0}^n (-1)^{k+1} \; \frac{2^{m+\frac{5}{2}} \, \pi^{2k+1}}{3^n} \;
 \frac{(2n)!}{k! \; (n-k)!} \\
 \times & \int_{\mathrm{UP}_{\epsilon}} du\wedge d\bar{u} \;  \frac{d\bar{\tau}}{d\bar{u}} \; \frac{u^{m+n-k} \, \eta^3(\tau)}{\omega^{2k+2}} \; \frac{\widehat{E}_2^k(\tau) \,
 \overline{\eta^3(\tau)}}{\sqrt{\im \tau}} \;.
\end{split}
\end{equation}
The integrals in (\ref{generating_fct1}) are all of the form
\begin{equation}\label{typical_term}
 \int_{\mathrm{UP}_{\epsilon}} d\tau\wedge d\bar{\tau} \; \dfrac{u^{m+n-k} \, \Delta \, \eta^3(\tau)}{\omega^{2k}} \;
\dfrac{\widehat{E}_2^k(\tau) \, \overline{\eta^3(\tau)}}{\sqrt{\im \tau}} \;.
\end{equation}
We use the substitution rule and Lemma \ref{EkQ} to write the integrand as $d\tau\wedge d\bar{\tau} \, \partial_{\bar{\tau}} (\dots)$.
With $\tau=x + i y$ and for an integrand $f$, we have
\begin{equation*}
 d\tau\wedge d\bar{\tau} \, \partial_{\bar{\tau}} f = - i \, dx \wedge dy \; \left( \partial_{x} + i \, \partial_{y} \right) \, f
= - \, d\Big(  f \, dx + i \, f \, dy \Big) \;.
\end{equation*}
Thus, the integral (\ref{typical_term}) reduces to the following integral over the boundary of $\mathrm{UP}_{\epsilon}$ which is the union
of the boundary of the $\epsilon$-discs around the points with singular fibers described in section \ref{section_sw_curves}:
\begin{equation}\label{typical_term2}
 X_{\epsilon}:= \sqrt{2} i \, \int_{\partial \mathrm{UP}_{\epsilon}}  \dfrac{u^{m+n-k} \, \Delta \, \eta^3(\tau)}{\omega^{2k}} \;
\mathcal{E}^k_{\frac{1}{2}}[Q] \; dz\;.
\end{equation}
The regularized $u$-plane integral is defined by
\begin{equation}\label{reg_term}
 \int^{\operatorname{reg}}_{\mathrm{UP}_{\epsilon}} d\tau\wedge d\bar{\tau} \; \dfrac{u^{m+n-k} \, \Delta \, \eta^3(\tau)}{\omega^{2k}} \;
\dfrac{\widehat{E}_2^k(\tau) \, \overline{\eta^3(\tau)}}{\sqrt{\im \tau}} := \lim_{\epsilon \to 0} X_{\epsilon}\;.
\end{equation}
Integrating each summand in (\ref{generating_fct1}) and taking the limit $\epsilon \to 0$, we obtain a well-defined procedure to
compute (\ref{uplane2}).
\subsection{The generating functions for $N_f=0,2,3$.}
\begin{lemma}\label{integration_by_parts}
In the notation of Lemmas \ref{choice_of_tau} and \ref{integral_wo_sings},
it follows that
\begin{equation*}
\begin{split}
 \widetilde{\pmb{\mathrm{Z}}}^{N_f}_{\mathrm{UP}_{\epsilon}}(p,S)  =  & \sum_{m, n \ge 0} \frac{p^m}{m!}  \frac{S^{2n}}{(2n)!} \; e^{S^2 \, \frac{\delta_{3,N_f}}{2}} \quad
\sum_{i=0}^n (-1)^{i} \; \frac{2^{m+4} \, \pi^{2i+2}}{3^n} \;
 \frac{(2n)!}{i! \; (n-i)!} \\
  \times  & \; \sum_{\operatorname{sing}} \operatorname{Coeff}^{N_f}_{m,n,i}\Big(u_{\operatorname{sing}}, \epsilon\Big),
\end{split}
\end{equation*}
with
\begin{equation}\label{coeff_wo_limit}
\begin{split}
\operatorname{Coeff}^{N_f}_{m,n,i}\Big(u_{\operatorname{sing}}, \epsilon \Big)
 =  \int_0^{K_s} dx_s \quad q_s \, \frac{du_s}{dq_s} \frac{\eta^3(\tau_s) \, u_s^{m+n-i}}{\omega_s^{2i+2}} \;
\mathcal{E}^i_{\frac{1}{2}} \left[ Q_s \right] \; \Big|_{\tau_s=x_s-i \, \ln{\epsilon}} \;.
\end{split}
\end{equation}
\end{lemma}
\begin{proof}
We use the substitution rule
and Lemma \ref{EkQ} to write each summand in the integrand as $d\tau\wedge d\bar{\tau} \, \partial_{\bar{\tau}} (\dots)$. Thus, the integral reduces to an integral over the boundary components of $\mathrm{UP}_{\epsilon}$. To carry out the integration along the boundary, we switch to the variable $\tau_s$ as described in
Lemma \ref{choice_of_tau}. In the $\tau_s$-coordinate the boundary is a line parallel to the $x_s$-axis passed from the right to the left. We use $2\pi i q_s \, \partial_{q_s} =
\partial_{\tau_s}$ and the claim follows.
\end{proof}
\begin{lemma}\label{evaluation}
The limit $ \operatorname{Coeff}^{N_f}_{m,n,i}(u_{\operatorname{sing}}) = \lim_{\epsilon \to 0} \operatorname{Coeff}^{N_f}_{m,n,i}(u_{\operatorname{sing}}, \epsilon)$ is
well-defined and equals
\begin{equation}\label{eval}
\begin{split}
 & \operatorname{Coeff}^{N_f}_{m,n,i}(u_{\operatorname{sing}})
  = K_s \;  \sum_{j=0}^i \, (-1)^j \; \binom{i}{j} \; \frac{\Gamma\left(\frac{1}{2}\right)}{\Gamma\left(\frac{1}{2}+j\right)} \; 2^{2j} \; 3^j \\
 \times & \operatorname{Coeff}_{q_s^0}\left[ \, q_s \frac{du_s}{dq_s} \frac{\eta^3(\tau_s) \, u_s^{m+n-i}}{\omega_s^{2i+2}} \;E^{i-j}_2(\tau_s) \;
 \left( q_s \frac{d}{dq_s}\right)^j  Q^+_s\left(\tau_s\right) \right] \;.
\end{split}
\end{equation}
In particular, it follows that for all
$u_{\operatorname{sing}} \not = \infty$ that
\begin{equation}
\operatorname{Coeff}^{N_f}_{m,n,i}(u_{\operatorname{sing}}) =0 \;.
\end{equation}
\end{lemma}
\begin{proof}
We first consider the contribution to integral from a singular point $u_{\operatorname{sing}} =\infty$.
We write $\mathcal{E}^k_{1/2}[Q]= q^{-1/8} ( \sum_{\alpha\ge 0} H_{k,-\alpha}(y) \, q^{-\alpha} + \sum_{\alpha\ge 0} H_{k,\alpha} \, q^{\alpha/(4-N_f)})$
where $H_{k,\alpha} \in \mathbb{Q}$. For $\alpha \in \mathbb{N}$, the coefficients $H_{k,-\alpha}(y)$ have the same exponential decay as the coefficient
$H_{-\alpha}(y)$ in equation (\ref{Qminus}). Using Lemma \ref{choice_of_tau} and (\ref{eta3}), it follows, for $|q|<1$, that the integrand in (\ref{generating_fct1})
has a normally convergent power series expansion
\begin{equation}\label{expansion2}
\begin{split}
& \; q^{\frac{m+n+2}{4-N_f}+1} \, \frac{du}{dq} \frac{\eta^3(\tau) \, u^{m+n-i}}{\omega^{2i+2}} \;
\mathcal{E}^i_{\frac{1}{2}}\left[ Q_{\infty} \right] \\
= \; & \sum_{\beta \ge 0} C_{\beta} \; q^{\frac{\beta}{4-N_f}} \;
\left( \sum_{\alpha\ge 0} H_{i,-\alpha}(y) \, q^{-\alpha} + \sum_{\alpha\ge 0} H_{i,\alpha} \, q^{\frac{\alpha}{4-N_f}} \right),
\end{split}
\end{equation}
where the $C_{\beta} \in \mathbb{Q}$.  For $\alpha \in \mathbb{Z}$, we have
that $\int_0^{4-N_f} dx \; q^{\alpha/(4-N_f)} = (4-N_f) \; \delta_{\alpha,0}$. Thus, the term-by-term integration in equation (\ref{coeff_wo_limit}) with respect to $dx$
extracts the term proportional to $q^0$. The term contains the linear combination of products of the coefficients $H_{k,-\alpha}(y)$ and $C_{\alpha}$, and it has exponential decay in $y$.
It follows that $\lim_{\epsilon \to 0} \operatorname{Coeff}^{N_f}_{m,n,i}(u_{\operatorname{sing}}, \epsilon )$ is
well-defined, equation (\ref{eval}) follows.
\noindent
Secondly, we consider the contribution to the integral from a singular point $u_{\operatorname{sing}} \not= \infty$.
We write $\mathcal{E}^k_{1/2}[Q_s]= q_s^{-1/8} ( \sum_{\alpha\ge 0} H_{s,k,-\alpha}(y_s) \, q_s^{-\alpha}
+ \sum_{\alpha\ge 0} H_{s,k,\alpha} \, q_s^{\alpha/K_s})$ where $H_{s,k,\alpha} \in \mathbb{Q}$. For $\alpha \in \mathbb{N}$,
the coefficients $H_{k,-\alpha}(y)$ have the same exponential decay as the
coefficients $H_{-\alpha}(y)$ in equation (\ref{Qminus}).  Using
Lemma~\ref{choice_of_tau} and (\ref{eta3}), it follows that for $|q_s| <1$ the integrand in (\ref{generating_fct1}) has a normally
convergent power series expansion of the form
\begin{equation}\label{expansion}
\begin{split}
 & q_s \, \frac{du_{s}}{dq_s} \;  \frac{\eta^3(\tau_s) \, u_{s}^{m+n-i} }{\omega_{s}^{2i+2}} \, \mathcal{E}^i_{\frac{1}{2}} \left[ Q_s \right] \\
= \; & \sum_{\beta \ge 1} C_{s,\beta} \; q_s^{\frac{\beta}{K_s}} \;
\left( \sum_{\alpha\ge 0} H_{s,i,-\alpha}(y_s) \, q_s^{-\alpha} + \sum_{\alpha\ge 0} H_{s,i,\alpha} \, q_s^{\frac{\alpha}{K_s}} \right),
\end{split}
\end{equation}
with $C_{s,\beta} \in \mathbb{Q}$. For $\alpha \in \mathbb{Z}$, we have that
$\int_0^{K_s} dx_s \; q_s^{\alpha/K_s} = K_s \, \delta_{\alpha,0}$. Thus, the term-by-term integration in equation (\ref{coeff_wo_limit}) with respect to $dx_s$ receives only one
contribution since the summation for $\beta$ starts at $1$. The term proportional to $q_s^0$ contains the linear combinations of the coefficients $H_{s,k,-\alpha}(y_s)$ and $C_{s,\alpha}$,
but not $H_{s,k,\alpha}$. Thus, this contribution has an exponential decay in $y_s$. It follows that
$\lim_{\epsilon \to 0} \operatorname{Coeff}^{0}_{m,n,i}(u_{\operatorname{sing}}, \epsilon )=0$.
\end{proof}
\begin{definition}
The regularized generating function of the massless $N_f$ low energy effective field theory on $\mathbb{C}\mathrm{P}^2$ is
\begin{equation*}
\begin{split}
 \pmb{\mathrm{Z}}^{N_f}_{\mathrm{UP}}(p,S)  =  & \sum_{m, n \ge 0} \frac{p^m}{m!}  \frac{S^{2n}}{(2n)!}\;\; e^{S^2 \; \frac{\delta_{3,N_f}}{2}}
 \quad \sum_{i=0}^n (-1)^{i} \; \frac{2^{m+4} \, \pi^{2i+2}}{3^n} \; \frac{(2n)!}{i! \; (n-i)!} \\
  \times  &  \; \operatorname{Coeff}^{N_f}_{m,n,i}\Big(u_{\operatorname{sing}}=\infty \Big) \;.
\end{split}
\end{equation*}
\end{definition}
\begin{theorem} \label{koeffizienten}
We have that
\begin{equation}\label{generating_function}
 \pmb{\mathrm{Z}}^{N_f}_{\mathrm{UP}}(p,S)  =  \sum_{m, n \ge 0} \pmb{\mathrm{D}}^{N_f}_{m,2n} \; \frac{p^m}{m!}  \frac{S^{2n}}{(2n)!} \;  e^{\left( - p + \frac{S^2}{3}\right) \; \delta_{3,N_f}}
\end{equation}
where the coefficients $\pmb{\mathrm{D}}^{N_f}_{m,n}$ are as follows:
\begin{equation}
\label{uplane0}
\begin{split}
 \pmb{\mathrm{D}}^{0}_{m,2n}
= & \sum_{i=0}^n \sum_{j=0}^i \frac{(-1)^{i+j+1}}{2^{n-2j-1} \; 3^{n-j}} \;
 \frac{(2n)!}{(n-i)! \; j! \; (i-j)!} \frac{\Gamma\left(\frac{1}{2}\right)}{\Gamma\left(j+\frac{1}{2}\right)}\\
\times & \operatorname{Coeff}_{q^0}\left[ \, \frac{\vartheta^9_4(\tau) \,
\left[\vartheta^4_2(\tau) + \vartheta^4_3(\tau)\right]^{m+n-i}}{\left[\vartheta_2(\tau) \, \vartheta_3(\tau)\right]^{2m+2n+3}} \;E^{i-j}_2(\tau) \;
 \left( q \frac{d}{dq}\right)^j  Q_{\infty}^{(0) \, +}\left(\tau\right) \right]\;,
\end{split}
\end{equation}
\begin{equation*}
\begin{split}
 \pmb{\mathrm{D}}^{2}_{m,2n}
= & \sum_{i=0}^n \sum_{j=0}^i \frac{(-1)^{i+j+1} \, 2^{-n+3j+2}}{3^{n-j}} \;
 \frac{(2n)!}{(n-i)! \; j! \; (i-j)!} \frac{\Gamma\left(\frac{1}{2}\right)}{\Gamma\left(j+\frac{1}{2}\right)}\\
\times & \operatorname{Coeff}_{q^0}\left[ \, \frac{\vartheta^{10}_4(\tau) \,
\left[\vartheta^4_2(\tau) + \vartheta^4_3(\tau)\right]^{m+n-i}}{\vartheta_2\left(\frac{\tau}{2}\right) \; \left[\vartheta_2(\tau) \, \vartheta_3(\tau)\right]^{2m+2n+3}}
\;E^{i-j}_2\left(\frac{\tau}{2}\right) \;
 \left( q \frac{d}{dq}\right)^j  Q^{(2) \, +}_{\infty}\left(\frac{\tau}{2}\right) \right]\;,
\end{split}
\end{equation*}
and
\begin{equation*}
\begin{split}
 \pmb{\mathrm{D}}^{3}_{m,2n}
= & \sum_{i=0}^n \sum_{j=0}^i (-1)^{i+j} \, (-1)^{m+n-j} \, \frac{2^{3m+2n+2j+5}}{3^{n-j}} \;
 \frac{(2n)!}{j! \; (i-j)! \; (n-i)!} \frac{\Gamma\left(\frac{1}{2}\right)}{\Gamma\left(j+\frac{1}{2}\right)}\\
\times & \operatorname{Coeff}_{q^0}\left[ \, \frac{\vartheta^{9}_2(\tau) \,
\left[\vartheta_3(\tau) \, \vartheta_4(\tau)\right]^{2m+2n-2i+3}}{\left[\vartheta_3^2(\tau) - \vartheta_4^2(\tau)\right]^{2m+2n+6}}
\;E^{i-j}_2\left(\tau\right) \;
 \left( q \frac{d}{dq}\right)^j  Q_{\infty}^{(3) \, +}\left(\tau\right) \right]\;.
\end{split}
\end{equation*}
\end{theorem}
\begin{proof}
In Lemmas~\ref{integration_by_parts} and \ref{evaluation},
we have already obtained an expansion of the coefficients. We use the explicit
Weierstrass presentations from Lemma~\ref{presentation} to write each coefficient as a ratio of Jacobi $\vartheta$-functions. We use
$2 \; \eta^3(\tau)=\vartheta_2(\tau) \, \vartheta_3(\tau) \, \vartheta_4(\tau)$ and $q \, du/dq = - \vartheta_4^8(\tau)/[8 \, \vartheta^2_2(\tau) \, \vartheta^2_3(\tau)]$
for $N_f=0$. For $N_f=2$, we use $\tau^{(2)}=\tau^{(0)}/2$. For $N_f=3$, we use $q \, du/dq = \vartheta_2^8(\tau) \, \vartheta_3^2(\tau) \, \vartheta_4^2(\tau)
/[\vartheta_3^2(\tau) - \vartheta_4^2(\tau)]^4$.
\end{proof}
\subsubsection{The invariants for $N_f=0$.}
\begin{lemma}
\label{invariants_Nf_0}
For all $m,n \in \mathbb{N}$ with $m+n \not \equiv 0\mod 2$ it follows $\pmb{\mathrm{D}}^{0}_{m,2n}=0$.
\end{lemma}
\begin{proof}
The claim follows from the fact that a factor $\vartheta_2^{2m+2n+3}(\tau) \, q^{\frac{1}{8}} = q^{\frac{m+n+2}{4}} (1 + O(\sqrt{q}))$
is present in the denominator of the expression defining $\pmb{\mathrm{D}}^0_{m,2n}$.
\end{proof}
\noindent
In the table below we list the first non-vanishing coefficients of the generating function in (\ref{generating_function}) for $N_f=0$.
They agree with the coefficients $-\Phi^{\mathbb{P}_2,H}_{H}$ computed in \cite[Thm.~4.4.]{EllingsrudGoettsche}.
{\footnotesize
\begin{center}
\begin{tabular}{|l||l|l|}
 \hline
 \raisebox{-0.5ex}{$p^m \, S^{2n}$} & \raisebox{-0.5ex}{$\pmb{\mathrm{D}}^{0}_{m,2n}$} & \raisebox{-0.5ex}{$\pmb{\mathrm{D}}^{0}_{m,2n}$}\\ [1ex]
 \hline
 \raisebox{-0.5ex}{$1$}   &  \raisebox{-0.5ex}{$-1$}                  &  \raisebox{-0.5ex}{$-\frac{1}{4} H_1 + 6 H_0$} \\ [1ex]
 \hline
  \raisebox{-0.5ex}{$S^4$} &   \raisebox{-0.5ex}{$- \frac{3}{16}   $}   &  \raisebox{-0.5ex}{$- \frac{49}{64} H_2 +  \frac{9}{4}  H_1 -  \frac{2133}{64}  H_0$} \\ [1ex]
 $p \, S^2$ &  $- \frac{5}{16}   $    & $ - \frac{7}{64} H_2 +  \frac{1}{4}  H_1 -   \frac{195}{64}  H_0$ \\ [0.5ex]
 $p^2$ &  $- \frac{19}{16}  $    & $  -\frac{1}{64} H_2 -  \frac{1}{4}  H_1 +   \frac{411}{64}  H_0$ \\ [0.5ex]
 \hline
  \raisebox{-0.5ex}{$S^8$}&  \raisebox{-0.5ex}{$ - \frac{232}{256}   $ } &
\raisebox{-0.5ex}{$- \frac{14641}{1024} H_3 + \frac{2401}{128}  H_2 +  \frac{44631}{1024}  H_1 + \frac{108741}{128} H_0$}\\ [1ex]
 $p \, S^6$&$ - \frac{152}{256}   $  & $- \frac{1331}{1024}  H_3 -    \frac{49}{128} H_2 +  \frac{10341}{1024}  H_1 - \frac{1749}{128}   H_0$\\ [0.5ex]
 $p^2 \, S^4$&$ - \frac{136}{256}   $  & $- \frac{121}{1024}   H_3 -    \frac{91}{128} H_2 +  \frac{2895}{1024}   H_1 - \frac{3687}{128}   H_0$\\ [0.5ex]
 $p^3 \, S^2$&$ - \frac{184}{256}   $  & $- \frac{11}{1024}    H_3 - \frac{29}{128}    H_2 +  \frac{589}{1024}    H_1 - \frac{753}{128}    H_0$\\ [0.5ex]
 $p^4$&$ - \frac{680}{256}   $  & $- \frac{1}{1024}     H_3 - \frac{7}{128}     H_2 -  \frac{505}{1024}    H_1 + \frac{1725}{128}   H_0$\\ [0.5ex]
 \hline
  \raisebox{-0.5ex}{$S^{12}$}&  \raisebox{-0.5ex}{$ - \frac{69525}{4096} $} &
\raisebox{-0.5ex}{$- \frac{11390625}{16384} H_4  + \frac{44838675}{16384} H_2 +  \frac{6075}{4}  H_1 - \frac{76478175}{2048} H_0$}\\ [1ex]
 $p \, S^{10}$&$ - \frac{26907}{4096} $ & $- \frac{759375}{16384}   H_4  - \frac{43923}{512}  H_3 + \frac{4833213}{16384}  H_2 +  \frac{185733}{512}  H_1 + \frac{5340591}{2048} H_0$ \\  [0.5ex]
 $p^2 \, S^8$&$ - \frac{12853}{4096} $ & $- \frac{50625}{16384}    H_4  - \frac{9317}{512}   H_3 + \frac{462707}{16384}   H_2 +  \frac{43587}{512}   H_1 + \frac{1179489}{2048} H_0$\\ [0.5ex]
 $p^3 \, S^6 $&$ - \frac{7803}{4096}  $ & $- \frac{3375}{16384}     H_4  - \frac{363}{128}    H_3 + \frac{861}{16384}      H_2 + \frac{2829}{128}     H_1 - \frac{69201}{2048}   H_0$\\ [0.5ex]
 $p^4 \, S^4 $&$ - \frac{6357}{4096}  $ & $- \frac{225}{16384}      H_4  - \frac{99}{256}     H_3 - \frac{21549}{16384}    H_2 + \frac{1653}{256}     H_1 - \frac{108639}{2048}  H_0$\\ [0.5ex]
 $p^5 \, S^2 $&$ - \frac{8155}{4096}  $ & $- \frac{15}{16384}       H_4  - \frac{25}{512}     H_3 - \frac{9475}{16384}     H_2 + \frac{815}{512}      H_1 - \frac{29265}{2048}   H_0$\\ [0.5ex]
 $p^6$&$ - \frac{29557}{4096} $ & $- \frac{1}{16384}        H_4  - \frac{3}{512}      H_3 - \frac{3021}{16384}     H_2 - \frac{619}{512}      H_1 + \frac{71649}{2048}   H_0$\\ [0.5ex]
\hline
\end{tabular}
\end{center}
}
\subsubsection{The invariants for $N_f=2$.}
In the table below we list the first coefficients of the generating function for $N_f=2$.
{\footnotesize
\begin{center}
\begin{tabular}{|l||l|l|}
 \hline
 \raisebox{-0.5ex}{$p^m \, S^{2n}$} & \raisebox{-0.5ex}{$\pmb{\mathrm{D}}^{2}_{m,2n}$} & \raisebox{-0.5ex}{$\pmb{\mathrm{D}}^{2}_{m,2n}$}\\ [1ex]
 \hline
 \raisebox{-0.5ex}{$1$}   &  \raisebox{-0.5ex}{$-3$}  &  \raisebox{-0.5ex}{$ -\frac{1}{4} H_2 + \frac{27}{4} H_0$} \\ [1ex]
 \hline
 \raisebox{-0.5ex}{$S^2$} & \raisebox{-0.5ex}{$0$} &  \raisebox{-0.5ex}{$- \frac{11}{16} H_3 + \frac{77}{16} H_1$} \\[1ex]
 $p$          & $0$                     & $- \frac{1}{16} H_3 + \frac{7}{16} H_1$ \\[0.5ex]
 \hline
  \raisebox{-0.5ex}{$S^4$}  &  \raisebox{-0.5ex}{$- \frac{21}{16}$}      & \raisebox{-0.5ex}{$- \frac{225}{64}  H_4 +  \frac{1043}{64}   H_2 -  \frac{567}{8}  H_0$}\\[1ex]
  $p \, S^2$&  $- \frac{27}{16}$     & $- \frac{15}{64}   H_4 +  \frac{61}{64}     H_2 -  \frac{9}{8}    H_0$\\[0.5ex]
  $p^2$     &  $- \frac{53}{16}$     & $- \frac{1}{64}   H_4 -  \frac{13}{64}    H_2 +  \frac{57}{8}   H_0$\\[0.5ex]
 \hline
  \raisebox{-0.5ex}{$S^6$} &  \raisebox{-0.5ex}{$0$}     &  \raisebox{-0.5ex}{$- \frac{6859}{256} H_5 + \frac{22869}{256} H_3 + \frac{12555}{128} H_1$} \\[1ex]
  $p \, S^4$       & $0$     & $- \frac{361}{256} H_5 + \frac{759}{256}  H_3 + \frac{2217}{128}  H_1$ \\[0.5ex]
  $p^2 \, S^2$     & $0$     & $- \frac{19}{256}  H_5 - \frac{115}{256}  H_3 + \frac{659}{128}  H_1$ \\[0.5ex]
  $p^3$     & $0$     & $- \frac{1}{256}   H_5 - \frac{33}{256}   H_3 + \frac{129}{128}  H_1$ \\[0.5ex]
 \hline
  \raisebox{-0.5ex}{$S^8$} &  \raisebox{-0.5ex}{$-\frac{3955}{256}$}     &  \raisebox{-0.5ex}{
  $- \frac{279841}{1024} H_6 + \frac{664875}{1024} H_4 + \frac{366667}{256} H_2 +\frac{4203535}{1024} H_0$} \\[1ex]
  $p \, S^6  $     &
  $-\frac{1925}{256}$     &
  $- \frac{12167}{1024} H_6 + \frac{10125}{1024}  H_4 + \frac{37709}{256}  H_2 - \frac{195895}{1024}$ \\[0.5ex]
  $p^2 \, S^4$     & $-\frac{1219}{256}$     &
  $- \frac{529}{1024}  H_6 - \frac{2565}{1024}  H_4 + \frac{5051}{256}  H_2 - \frac{61409}{1024} H_0$ \\[0.5ex]
  $p^3 \, S^2$     & $-\frac{949}{256} $     &
  $- \frac{23}{1024}   H_6 - \frac{451}{256}   H_4 + \frac{541}{256}  H_2 - \frac{1735}{1024} H_0$ \\[0.5ex]
  $p^4       $     & $-\frac{1811}{256}$     &
  $- \frac{1}{1024} H_6 - \frac{53}{1024} H_4 - \frac{85}{256} H_2 + \frac{15151}{1024} H_0$ \\[0.5ex]
 \hline
\end{tabular}
\end{center}
}
\begin{lemma}
\label{invariants_Nf_2}
For all $m,n \in \mathbb{N}$ with $m+n \not\equiv 0\mod 2$ it follows $\pmb{\mathrm{D}}^{2}_{m,2n}=0$.
\end{lemma}
\noindent
The lemma follows from the explicit formulas for these $q$-series.
These series are given in terms of sums of products of derivatives
of modular forms and harmonic Maass forms which are explicitly given
in terms of theta functions and the $\mu$-function of Zwegers (see Section
6). We leave the proof of this lemma to the reader.
\subsubsection{The invariants for $N_f=3$.}
\begin{lemma}\label{Sduality}
For all $m,n \in \mathbb{N}$ it follows:
\begin{equation}
\begin{split}
  \pmb{\mathrm{D}}^{3}_{m,2n} \left[ \, Q^{(3)}_{\infty}(\tau) \right] & =
  \pmb{\mathrm{D}}^{3}_{m,2n} \left[ \, - \,Q\left(\tau\right)\right] \;.
\end{split}
\end{equation}
\end{lemma}
\begin{proof}
We have the $S$-duality equation
\begin{equation}
\label{S-duality}
\begin{split}
\frac{1}{\sqrt{-i \, \tau}} \; \zeta_8 \; Q\left(1- \frac{1}{\tau}\right)
 & =  - \zeta_8 \; Q(\tau+1) \;.
\end{split}
\end{equation}
In Lemma \ref{transfo_Z_rho}, we will show $Q(\tau) = \zeta_8 \, Q(\tau+1) + 14 \, \eta(\tau) \, \rho^4(\tau)$
where we have set $\rho(\tau)=\sqrt{2}\, \eta^2(2\tau)/\eta^2(\tau)$. It then follows
\begin{equation}
\begin{split}
\pmb{\mathrm{D}}^{3}_{m,2n} \left[ \, - \,Q\left(\tau\right)\right] & = \pmb{\mathrm{D}}^{3}_{m,2n} \left[ \,- \,\zeta_8 \, Q\left(\tau+1\right)\right] -
14 \, \pmb{\mathrm{D}}^{3}_{m,2n} \left[ \, \eta(\tau) \, \rho^4(\tau) \right] \;.
\end{split}
\end{equation}
On the other hand, we have
\begin{equation*}
\begin{split}
\pmb{\mathrm{D}}^{3}_{m,2n} \left[ \,- \,\zeta_8 \, Q\left(\tau+1\right)\right]
& = \pmb{\mathrm{D}}^{3}_{m,2n} \left[ \frac{1}{\sqrt{-i \, \tau}} \; \left\lbrace Q\left(- \frac{1}{\tau}\right) - 14 \, \eta\left(- \frac{1}{\tau}\right) \,
    \rho^4\left(- \frac{1}{\tau}\right) \right\rbrace \right]\\
& = \pmb{\mathrm{D}}^{3}_{m,2n} \left[ \frac{1}{\sqrt{-i \, \tau}} Q\left(- \frac{1}{\tau}\right) - 14 \, \eta\left(\tau\right) \,
    \dfrac{1}{\rho^4\left(\frac{\tau}{2}\right)} \right]\;.
\end{split}
\end{equation*}
Hence, we obtain
\begin{equation*}
\begin{split}
\pmb{\mathrm{D}}^{3}_{m,2n} \left[ \, - \,Q\left(\tau\right)\right] & =  \pmb{\mathrm{D}}^{3}_{m,2n} \left[ \, Q^{(3)}_{\infty}(\tau) \right]
- 14 \, \pmb{\mathrm{D}}^{3}_{m,2n} \left[ \, \eta(\tau) \, \left( \rho^4(\tau) +  \dfrac{1}{\rho^4\left(\frac{\tau}{2}\right)} \right) \right] \;.
\end{split}
\end{equation*}
One then checks that the last term on the right hand side vanishes.
\end{proof}
\smallskip
\noindent
In the table below we list the first coefficients of the generating function for $N_f=3$.
{\footnotesize
\begin{center}
\begin{tabular}{|l||l|l|}
 \hline
 \raisebox{-0.5ex}{$p^m \, S^{2n}$} & \raisebox{-0.5ex}{$\pmb{\mathrm{D}}^{3}_{m,2n}$} & \raisebox{-0.5ex}{$\pmb{\mathrm{D}}^{3}_{m,2n}$}\\ [1ex]
 \hline
 \raisebox{-0.5ex}{$1$}   &  \raisebox{-0.5ex}{$-\frac{5}{4}$} &  \raisebox{-0.5ex}{$-\frac{1}{16} H_4 + \frac{3}{16} H_2 + \frac{3}{2} H_0$} \\ [1ex]
 \hline
  \raisebox{-0.5ex}{$S^2$} &   \raisebox{-0.5ex}{$-\frac{95}{96}   $}   &  \raisebox{-0.5ex}{$ - \frac{23}{128} H_6 +  \frac{119}{384}  H_4 +  \frac{45}{32}  H_2 + \frac{313}{128} H_0$} \\ [1ex]
 $p$ &  $  \phantom{-} \frac{45}{32}   $    & $ - \frac{1}{128} H_6 +  \frac{11}{128}  H_4 -   \frac{5}{32}  H_2 - \frac{209}{128} H_0$ \\ [0.5ex]
 \hline
  \raisebox{-0.5ex}{$S^4$} &   \raisebox{-0.5ex}{$ - \frac{1787}{768} $}   &  \raisebox{-0.5ex}{$  - \frac{961}{1024} H_8 +  \frac{851}{1024}  H_6 +  \frac{133}{24}  H_4
+ \frac{4587}{1024} H_2 - \frac{171}{128} H_0$} \\ [1ex]
 $p \, S^2$ &  $ \phantom{-}\frac{201}{256}   $    & $ - \frac{31}{1024} H_8 +  \frac{743}{3076}  H_6 -   \frac{5}{24}  H_4 - \frac{4577}{3072} H_2 - \frac{991}{384} H_0$ \\ [0.5ex]
 $p^2$ &  $ - \frac{489}{256}  $    & $  - \frac{1}{1024} H_8 +  \frac{19}{1024}  H_6 - \frac{1}{8} H_4 +   \frac{171}{1024}  H_2 + \frac{277}{128} H_0$ \\ [0.5ex]
 \hline
  \raisebox{-0.5ex}{$S^6$}&  \raisebox{-0.5ex}{$ - \frac{189187}{18432} $ } &
\raisebox{-0.5ex}{$ - \frac{59319}{8192} H_{10} + \frac{12493}{8192}  H_8 + \frac{70403}{2048}  H_6 + \frac{3091945}{73728} H_4 + \frac{600451}{12288} H_2 - \frac{970759}{12288} H_0$}\\ [1ex]
 $p \, S^4$&$   \phantom{-}\frac{2211}{2048} $  & $- \frac{1521}{8192}  H_{10} + \frac{9579}{8192} H_8 -  \frac{2065}{6144}  H_6 - \frac{128731}{24576}   H_4 - \frac{29563}{12288} H_2 +
\frac{35039}{12288} H_0$\\ [0.5ex]
 $p^2 \, S^2$&$ - \frac{1627}{2048}   $  & $ - \frac{39}{8192}   H_{10} + \frac{1751}{24576} H_8 - \frac{2087}{6144}  H_6 + \frac{4051}{24576}  H_4 + \frac{7953}{4096} H_2
+ \frac{40585}{12288} H_0$\\ [0.5ex]
 $p^3 $&$ \phantom{-} \frac{5843}{2048}  $  & $- \frac{1}{8192}  H_{10} + \frac{27}{8192} H_8 - \frac{75}{2048} H_6 + \frac{1575}{8192} H_4 - \frac{825}{4096} H_2
- \frac{12987}{4096} H_0$\\ [0.5ex]
 \hline
\end{tabular}
\end{center}
}
\medskip
\begin{remark}
In the case $N_f=3$, using the variable $u^{(3)}_{SW} = u^{(3)} + 1/2$ of \cite{SeibergWitten2} to define the $u$-plane integral
leads to the generating function in Equation (\ref{generating_function}) with the term $\exp(-p)$ removed.
\end{remark}
\subsection{The partition function for $N_f=4$.}
Here we consider the case where $N_f=4$. We set $\alpha=\vartheta_2^4(\tau_0)$, $\beta=\vartheta_3^4(\tau_0)$.
\begin{lemma}\label{presentation2}
The Weierstrass presentation in the case $\# \,57$
is related to the Weierstrass presentation in the case $\# \,71$ in Lemma~\ref{presentation} by:
{\footnotesize
\begin{equation*}
\begin{array}{|c|ccccc||lcl|}
\hline
 \multicolumn{6}{|l||}{N_f =  4,\; \# \, 57}
   &
 \multicolumn{3}{|l|}{\mathrm{UP}^{(4)}_\epsilon =\mathbb{C}\mathrm{P}^1 - B_{\epsilon}(0)- B_{\epsilon}( u_\pm )
- \widetilde{B}_{\epsilon}(\infty) \phantom{\widehat{B^{(4)}}}} \\
\hline
 E_{\operatorname{sing}} & I_2 & I_2 & I_2  & I_0^* &  &   &  & \\
\cline{1-6}
\tau^{(4)}_{\operatorname{sing}} & \infty   & 0   & 1 & \tau_0 && u^{(4)}
&=&  \frac{8 \, \alpha \, \beta \,m^2}{(\alpha-\beta) \, u^{(2)} + \alpha + \beta} \\
u^{(2)}_{\operatorname{sing}}                       & \infty  & -1   & 1 & - \frac{\alpha+\beta}{\alpha-\beta} &&
 u_{\pm}            & = & 4\, \beta \, m^2 \;, \;\; 4\, \alpha \, m^2 \\
 u^{(4)}_{\operatorname{sing}}  & 0 &  u_-  &  u_+  & \infty && \omega^{(4)}(\tau) & = &  \frac{2 \, \sqrt{2} \, \omega^{(2)}(\tau)}{\sqrt{\alpha - \beta} \, \sqrt{u^{(4)}}}\\
\hline
\end{array}
\end{equation*}
\begin{equation*}
\begin{split}
  \Delta^{(4)}   & =  \frac{1}{256} \, \alpha^2 \, \beta^2 \, (\alpha-\beta)^2  \, \frac{1}{2^{12}}
 \, \left(u^{(4)}\right)^2
 \, \left[u^{(4)} - 4 \, \beta  \, m^2\right]^2
 \, \left[u^{(4)} - 4 \, \alpha \, m^2\right]^2 \\
  & = \frac{\eta^{24}(\tau_0)}{2^{12}} \, u^2
 \, \left[u^{(4)} - 4 \,\alpha \, m^2\right]^2
 \, \left[u^{(4)} - 4 \,\beta  \, m^2\right]^2
\end{split}
\end{equation*}
}
\end{lemma}
\begin{proof}
 The parametrization in the case $ \# \, 57$
 are obtained from $\# 71$ in Lemma~\ref{presentation}
 by transferring a star. This determines their Weierstrass presentations.
 \end{proof}
\begin{remark}
The case $\# \, 57$ agrees with the case described in \cite[Sect.~17.4]{SeibergWitten2} if we set $m_1 = m_2 = \frac{m}{2}$ and $m_3=m_4=0$.
\end{remark}
\begin{remark}
The $\epsilon$-disc $\widetilde{B}_{\epsilon}(\infty)$ around the point $\infty$ with a singular fiber of the Kodaira-type $I^*_{0}$ is constructed as follows:
the boundary of the $\epsilon$-disc corresponds to the circle in the $u$-coordinate with $u=R e^{i\theta}$, $R=1/\epsilon$, and $\theta \in [0 ; 2\pi]$.
It is easy to show that Lemmas~\ref{integral_wo_sings}, \ref{integration_by_parts}, \ref{evaluation} still hold for the elliptic families in Lemma~\ref{presentation2}
for all singular points  with
$|u_{\operatorname{sing}}| < \infty$. Thus, the partition function receives contributions only from the singular point $u=\infty$.
\end{remark}
\begin{definition}\label{integral_wo_sings-pSb}
The partition function of the $N_f=4$ low energy effective field theory on $\mathbb{C}\mathrm{P}^2$ is
\begin{equation}\label{uplane3b}
\begin{split}
 \widetilde{\pmb{\mathrm{Z}}}^{4}_{\mathrm{UP}_\epsilon}(m,\tau_0)
 & = - \frac{8}{\sqrt{2\pi}} \; f(m^2,\tau_0) \; \int_{\mathrm{UP}_\epsilon}  \frac{du\wedge d\bar{u}}{\sqrt{\im\tau}} \;
\frac{d\bar{\tau}}{d\bar{u}} \; \frac{\Delta^{\frac{1}{8}}}{\omega^{\frac{1}{2}}}  \; \overline{\eta^3(\tau)} \;,
\end{split}
\end{equation}
where $f(m^2,\tau_0) = \frac{1}{4^{4} \, m^4 \, \eta^{12}(\tau_0)}$ is a universal normalization factor.
The regularized partition function of the massive $N_f=4$ low energy effective field theory on $\mathbb{C}\mathrm{P}^2$
is
\begin{equation*}
\pmb{\mathrm{Z}}^{4}_{\mathrm{UP}}(m,\tau_0) = \lim_{\epsilon \to 0}
\widetilde{\pmb{\mathrm{Z}}}^{4}_{\mathrm{UP}_\epsilon}(m,\tau_0) \;.
\end{equation*}
\end{definition}
\begin{remark}
The function $f(m^2,\tau_0)$ is a universal normalization factor appearing in the $u$-plane integral
which is familiar from the $\mathcal{N}=4$ case \cite{LabastidaLozano}.
\end{remark}
\begin{lemma}
It follows that
\begin{equation}
\begin{split}
\widetilde{\pmb{\mathrm{Z}}}^{4}_{\mathrm{UP}_{\epsilon}}  (m,\tau_0)
& =  16 \, f(m^2,\tau_0) \,  \int_0^{1} d\hat{\theta}  \; \; \left.\left(  u \, \Delta^{\frac{1}{6}} \; \frac{\pmb{\mathrm{Z}}(\tau)}{\eta^4(\tau)}  \right)
\right|_{u =\frac{1}{\epsilon} e^{2\pi i\hat{\theta}}} \;.
\end{split}
\end{equation}
\end{lemma}
\begin{proof} The integral reduces to an integral over
the boundary component of $\mathrm{UP}_{\epsilon}$ at $u=\infty$. To carry out the integration along the boundary we set $u=R \, e^{i\theta}$
(with $\theta$ running clockwise from $0$ to $-2\pi$) and use
\begin{equation}
 du \wedge d\bar{u} \; \partial_{\bar{u}} f =  - d\left( f \, \frac{u}{R} \, dR + i \, f \, u \, d\theta \right) \;.
\end{equation}
\end{proof}
\begin{lemma}
The regularized partition function of the massive $N_f=4$ low energy effective field
theory on $\mathbb{C}\mathrm{P}^2$ is
\begin{equation}
\begin{split}
  \pmb{\mathrm{Z}}^{4}_{\mathrm{UP}}(m,\tau_0)
  =  \left. \left[ \frac{1}{2}\, \frac{q}{\eta^4(\tau)} \, \frac{d}{dq} \left( \frac{q}{\eta^4(\tau)}  \, \frac{d}{dq} \right)
+ g(\tau)  \right] \, \frac{\pmb{\mathrm{Z}}(\tau)}{\eta^4(\tau)}\right|_{\tau=\tau_0}
\;,
\end{split}
\end{equation}
where
\begin{equation*}
 g(\tau) = - \dfrac{1}{2^2 \, 3^2 \, \eta^8(\tau)} \, \left[ \alpha^2 - \alpha \beta + \beta^2 \right] \;.
\end{equation*}
\end{lemma}
\begin{proof}
In the case $\# \, 57$, we have
\begin{equation*}
\begin{split}
\frac{4 \, u^{(4)}}{\eta^4(\tau_0)} \; \Delta^{\frac{1}{6}}  & =  \left(u^{(4)}\right)^2 - \left( \alpha + \beta \right) \frac{4 \, m^2 \, u^{(4)}}{3} \\
&  - \left( \alpha^2 - \alpha \, \beta + \beta^2\right) \frac{16 \,m^4}{9} + m^4 O\left( \frac{m^2}{u^{(4)}} \right)\;.
\end{split}
\end{equation*}
Using
\begin{equation*}
u^{(2)} + \frac{\alpha+\beta}{\alpha-\beta} = \frac{8 \, \alpha \, \beta\, m^2 }{(\alpha-\beta) \; u^{(4)}} \;,
\end{equation*}
we obtain the series expansion
\begin{equation*}
\begin{split}
 &  \frac{\pmb{\mathrm{Z}}(\tau)}{\eta^4(\tau)}   = \frac{\pmb{\mathrm{Z}}(\tau_0)}{\eta^4(\tau_0)}  \\
+ & \; \left[\eta^4(\tau) \, \frac{d \tau\;\,}{du^{(2)}} \right]_{\tau_0} \;
\left[\frac{1}{\eta^4(\tau)} \, \frac{d}{d\tau} \left(\frac{\pmb{\mathrm{Z}}(\tau)}{\eta^4(\tau)} \right)\right]_{\tau_0} \;
  \frac{8 \, \alpha \, \beta \, m^2}{(\alpha-\beta) \, u^{(4)}}
\\
+ & \;  \frac{1}{2} \; \left[\frac{d}{du^{(2)}}\left(\eta^4(\tau) \, \frac{d \tau\;\,}{du^{(2)}} \, \right)\right]_{\tau_0} \;
 \left[\frac{1}{\eta^4(\tau)} \, \frac{d}{d\tau} \left(\frac{\pmb{\mathrm{Z}}(\tau)}{\eta^4(\tau)} \right)\right]_{\tau_0} \;
\left(\frac{8 \, \alpha \, \beta \, m^2}{(\alpha-\beta) \, u^{(4)}}\right)^2
\\
+ & \; \frac{1}{2} \; \left[\eta^4(\tau) \, \frac{d \tau\;\,}{du^{(2)}}\right]^2_{\tau_0} \;
 \left[\frac{1}{\eta^4(\tau)} \, \frac{d}{d\tau} \left\lbrace \frac{1}{\eta^4(\tau)} \, \frac{d}{d\tau} \left(\frac{\pmb{\mathrm{Z}}(\tau)}{\eta^4(\tau)} \right)
\right\rbrace \right]_{\tau_0} \;
\left(\frac{8 \, \alpha \, \beta \, m^2}{(\alpha-\beta) \, u^{(4)}}\right)^2
\\
+ & \; O\left( \left[\frac{m^2}{u^{(4)}} \right]^3\right) \;,
\end{split}
\end{equation*}
where
\begin{equation*}
 2\pi i\; \eta^4(\tau) \;\frac{d \tau\;\,}{du^{(2)}} = - \frac{\pi^2}{4} \, \frac{\eta^4(\tau)}{\omega^{(2)}(\tau)^2 \, \Delta^{(2)}(\tau)^\frac{1}{2}}
= - \frac{1}{4 \; \left[\Delta^{(2)}\right]^{\frac{1}{3}}}\;.
\end{equation*}
A tedious calculation using $2^4 \, \eta^{12}(\tau_0)= \alpha \, \beta \, (\alpha-\beta)$, and
\begin{equation*}
  \left[ \Delta^{(2)} \right]_{\tau_0} = \frac{1}{4} \; \frac{(\alpha \,\beta)^2}{(\alpha - \beta)^4}\;, \quad \qquad
\left[ \partial_{u^{(2)}} \Delta^{(2)} \right]_{\tau_0} = - \frac{1}{4} \; \frac{\alpha \,\beta \,(\alpha+\beta)}{(\alpha - \beta)^3}\;,
\end{equation*}
then gives
\begin{equation*}
\begin{split}
 \left[\frac{u^{(4)}}{m^4 \; \eta^4(\tau_0)} \;   \Delta^{\frac{1}{6}} \; \frac{\pmb{\mathrm{Z}}(\tau)}{\eta^4(\tau)} \right]_{[u^{(4)}]^{0}}
= & - \frac{4}{9} \, \left( \alpha^2 -  \alpha \, \beta + \beta^2 \right) \; \frac{\pmb{\mathrm{Z}}(\tau_0)}{\eta^4(\tau_0)} \\
& + 8 \,\eta^8(\tau_0) \;
 \left[\frac{q}{\eta^4(\tau)} \, \frac{d}{dq} \left\lbrace \frac{q}{\eta^4(\tau)} \, \frac{d}{dq} \left(\frac{\pmb{\mathrm{Z}}(\tau)}{\eta^4(\tau)} \right)
\right\rbrace\right]_{\tau_0} \;.
\end{split}
\end{equation*}
The statement the follows from using the regularized integral
\begin{equation*}
 \lim_{\epsilon \to 0} \int_0^{1} d\hat{\theta}  \; u^k  \Big|_{u =e^{2\pi i\hat{\theta}}/\epsilon} = \delta_{k,0} \;.
\end{equation*}
\end{proof}
\noindent
The family of elliptic curves for the conformally invariant case corresponds to the
\emph{massless} $N_f=4$ low energy effective field theory on $\mathbb{C}\mathrm{P}^2$. It follows:
\begin{theorem} \label{partition_function_tau_a}
The regularized generating function of the massless $N_f=4$ low energy effective field
theory on $\mathbb{C}\mathrm{P}^2$ is given by
\begin{equation}\label{partition_function_tau}
\begin{split}
  \pmb{\mathrm{Z}}^{4}_{\mathrm{UP}}(0,\tau)  =  \left[ \frac{1}{2}\, \frac{q}{\eta^4(\tau)} \, \frac{d}{dq} \left( \frac{q}{\eta^4(\tau)}  \, \frac{d}{dq} \right)
+ g(\tau)  \right] \, \frac{\pmb{\mathrm{Z}}(\tau)}{\eta^4(\tau)}
\;,
\end{split}
\end{equation}
where
\begin{equation*}
  g(\tau) =
- \dfrac{1}{2^2 \, 3^2 } \, \left[ \left(\frac{\vartheta_2(\tau)}{\eta(\tau)}\right)^8 - \left(\frac{\vartheta_2(\tau)}{\eta(\tau)}\right)^4
\left(\frac{\vartheta_3(\tau)}{\eta(\tau)}\right)^4 + \left(\frac{\vartheta_3(\tau)}{\eta(\tau)}\right)^8 \right]\;.
\end{equation*}
We have the following transformation properties under modular transformations:
\begin{equation*}
\begin{split}
  \pmb{\mathrm{Z}}^{4}_{\mathrm{UP}}(0,\tau+2) = \pmb{\mathrm{Z}}^{4}_{\mathrm{UP}}(0,\tau) \;,\qquad
  \pmb{\mathrm{Z}}^{4}_{\mathrm{UP}}\left(0,-\frac{1}{\tau}+ 1\right)   =   - \, \pmb{\mathrm{Z}}^{4}_{\mathrm{UP}}(0,\tau + 1) \;.
\end{split}
\end{equation*}
\end{theorem}
\begin{proof}
The family of elliptic curves for the conformally invariant case is obtained from $\# \, 57$ in the limit $m \to 0$.
We start from the generating function in the case $\# \, 57$ to obtain a duality group $\Gamma(2)$.
In the limit $m \to 0$, the fibration becomes a rational elliptic surface with
a constant $j$-invariant and two singularities of Kodaira-type $I_0^*$, one at $u=0$ and one at $u=\infty$.
We have $g_3 = \pi^6 \, E_6(\tau)/[216\, \omega^6]$, $g_2 = \pi^4 \, E_4(\tau)/[12\, \omega^4]$,
and $\omega=2\pi/\sqrt{u}$  and $\tau$ arbitrary.
Higher terms proportional to $p^k S^{2l}$ in the generating function $\pmb{\mathrm{Z}}^{4}_{\mathrm{UP}}(p,S)$ will be proportional to
$m^{2n}$ with $n=k+l\ge 1$. The reason is that the integrand contains higher powers of $u$ or $\widehat{T}$.
The series expansion then introduces higher derivatives $u^{(0) \; n} \, \partial_{u^{(0)}}^n$ which are proportional to $m^{2n}$. Thus,
in the limit $m \to 0$ only the constant term in the generating function $\pmb{\mathrm{Z}}^{4}_{\mathrm{UP}}(p,S)$ survives.
The transformation under the shift follows from Lemma \ref{transfo_Z_rho} and the fact that $g(\tau)/\eta^4(\tau)$
is invariant under $\tau \mapsto \tau+2$. The remaining transformation follows from Lemma \ref{transfo_Z_rho} and the fact
that $g(1+\tau)=g(1-\frac{1}{\tau})$.
\end{proof}
\subsection{Criterion for Theorem~\ref{MainTheorem}}
From a physics point of view, at a high energy scale the $\mathrm{SO}(3)$-Donaldson theory
is described by the low energy effective field theory. We show that the physics speculation by
Moore and Witten on how to compute Donaldson invariants from the Seiberg-Witten curves can be turned into
a mathematical statement which we then prove.
The cuspidal contributions to the generating function of the low energy effective field theory
should be equal to the generating function of the $\mathrm{SO}(3)$-Donaldson theories with $2 \, N_f$ massless monopoles.
We conjecture the following relationship between the Donaldson invariants
defined in Equations (\ref{donaldson}) and (\ref{donaldson2b}) and the $u$-plane invariants given
in Equation (\ref{generating_function}):
\begin{conjecture}
Assume $c=\operatorname{H}$. For $N_f=0, 2, 3$ and $2m+2n+4=(4-N_f)k$, we have
\begin{subequations}\label{conjecture}
\begin{align}
\label{first}
  \pmb{\mathrm{D}}^{0}_{m,2n}   & =    \pmb{\Phi}_{k,m,2n} \;,\\
 \label{second}
  \pmb{\mathrm{D}}^{N_f}_{m,2n} & =  2^{\frac{(N_f+2)\,k}{2}}  \pmb{\Phi}^{N_f,c,1}_{k,m,2n} \;.
\end{align}
\end{subequations}
\end{conjecture}
\noindent
For $N_f=0$, this conjecture is the conjecture of Moore and Witten.
Lemma \ref{invariants_Nf_0} and Lemma \ref{invariants_Nf_2} show that for $N_f=0, 2$ the invariants vanish on both sides
for $m+n \not \equiv 0 \mod 2$. For $N_f=0$, the conjecture is equivalent to
the assertion that the generating functions in (\ref{generatingfunction}) and
in (\ref{generating_function}) are equal. We will prove Equation (\ref{first}) by proving:
\begin{theorem}\label{criterion}
Theorem~\ref{MainTheorem} is equivalent to the vanishing of constant
terms, for every pair of non-negative integers $m$ and $n$, of
the series
\begin{equation}
\label{difference}
\begin{split}
 \sum_{k=0}^n \sum_{j=0}^k & \; (-1)^{j}  \; \frac{(2n)!}{(n-k)! \; j! \; (k-j)!}
 \frac{\vartheta_4^8(\tau) \, \left[ \vartheta_2^4(\tau) + \vartheta_3^4(\tau)\right]^{m}}
{\left[ \vartheta_2(\tau) \, \vartheta_3(\tau)\right]^{2m+2n+3}}
 \; E^{k-j}_2(\tau)  \\
\times & \left[ \frac{(-1)^{n+1}}{2^{k-3} \, 3^{k}} \; \frac{(n-k)!}{(2n-2k)!}
 \; \left[ \vartheta_2^4(\tau) + \vartheta_3^4(\tau)\right]^{j} \; F_{2(n-k)}(\tau) \right. \\
+ & \left. \frac{(-1)^{k+1}}{2^{n-2j-1} \; 3^{n-j}} \;
 \frac{\Gamma\left(\frac{1}{2}\right)}{\Gamma\left(j+\frac{1}{2}\right)}
 \vartheta_4(\tau) \, \left[\vartheta^4_2(\tau) + \vartheta^4_3(\tau)\right]^{n-k} \; \left( q \frac{d}{dq}\right)^j  Q^+\left(\tau\right) \right] \;,
\end{split}
\end{equation}
where series $F_t(\tau)$ are defined in (\ref{Ft_old}).
\end{theorem}
\begin{proof}
In Theorem~\ref{koeffizienten} we established explicit formula (\ref{uplane0}) for the contribution from the cusp at $\tau=\infty$ in the regularized $u$-plane
integral. Theorem~\ref{Goettsche_thm} contained explicit formula (\ref{Goettsche}) for the generating function of the Donaldson invariants of $\mathbb{C}\mathrm{P}^2$.
The difference of Equations (\ref{Goettsche}) and (\ref{uplane0}) is the constant coefficient in Equation (\ref{difference}) where we have used that 
$Q_{\infty}^{(0) \, +}\left(\tau\right)= Q^+\left(\tau\right)$. Hence the vanishing of
the constant coefficient term in Equation (\ref{difference}) is equivalent to Theorem~\ref{MainTheorem}.
\end{proof}
\noindent
\begin{remark}
For $N_f= 2, 3$, it is difficult to compute the geometrically defined invariants $\pmb{\Phi}^{N_f,c,1}_{k,m,2n}$ directly. In contrast, the invariants $\pmb{\Phi}^{N_f,c,0}_{k,m,2n}$ defined in Equation
(\ref{donaldson2}) are easy to compute by using Lemma \ref{higher_invariants}. However, there is no
known formula to compute the error term in Equation (\ref{errorterm1}) to relate the latter to the former. However, we successfully
performed the following check of the Equation (\ref{second}) for $N_f=2,3$.
\smallskip
The Donaldson invariants $\pmb{\Phi}^{N_f,c,1}_{k,m,2n}$ and
$\pmb{\Phi}^{N_f,c,0}_{k,m,2n}$ are rational polynomials in the coefficients $H_k$ of the holomorphic part of the Maass form $Q(\tau)$.
By definition, the same statement is true for the invariants $\pmb{\mathrm{D}}^{N_f}_{m,2n}$ as well.
Geometrically, the invariants $\pmb{\Phi}^{N_f,c,1}_{k,m,2n}$ and $\pmb{\Phi}^{N_f,c,0}_{k,m,2n}$ differ by the contributions from the lower strata in the Uhlenbeck compactification.
But the contribution from the highest stratum in the Uhlenbeck compactification agrees whence the coefficients of the highest $H_k$ in $\pmb{\Phi}^{N_f,c,1}_{k,m,2n}$
and $\pmb{\Phi}^{N_f,c,0}_{k,m,2n}$ must agree. We proved the following weak version of a Moore-Witten type conjecture~(\ref{second}) for the Donaldson invariants for $N_f=2,3$:
In terms of the coefficients $H_k$ of the holomorphic part of the weight $1/2$ harmonic Maass form $Q(\tau)$,
we have for $c=\operatorname{H}$, $k$ even, and $2m+2n+4=(4-N_f)\,k$ that
\begin{equation}
 \operatorname{Coeff}_{H_k}\left[ \pmb{\mathrm{D}}^{N_f}_{m,2n} - 2^{\frac{(N_f+2)\,k}{2}} \;  \pmb{\Phi}^{N_f,c,0}_{k,m,2n} \right] = 0 \;.
\end{equation}
The proof will appear somewhere else.
\end{remark}
\section{Harmonic Maass forms}\label{HarmonicMaassForms}
We shall use Theorem~\ref{criterion} to prove Theorem~\ref{MainTheorem}.
To this end, we make use of the theory of harmonic Maass forms, which
we briefly recall here.
\subsection{Definitions and facts}
Following Bruinier and Funke \cite{BF},
we define the notion of a harmonic weak Maass form of weight $k\in
\frac{1}{2}\Z$ as follows. We let $\tau=u+iv\in \H$ with $u, v \in \R$,
and throughout we let $q:=e^{2\pi i \tau}$.
We define the
weight $k$ hyperbolic Laplacian $\Delta_k$ by
\begin{equation}
\Delta_k := -v^2\left( \frac{\partial^2}{\partial u^2} +
\frac{\partial^2}{\partial v^2}\right) + ikv\left(
\frac{\partial}{\partial u}+i \frac{\partial}{\partial v}\right).
\end{equation}
For odd integers $d$, define $\epsilon_d$ by
\begin{equation}\label{epsilond}
\epsilon_d:=\begin{cases} 1 \ \ \ \ &{\text {\rm if}}\ d\equiv
1\pmod 4,\\
i \ \ \ \ &{\text {\rm if}}\ d\equiv 3\pmod 4. \end{cases}
\end{equation}
\begin{definition}\label{defhwmf}
If $N$ is a positive integer (with $4\mid N$ if $k\in
\frac{1}{2}\Z\setminus \Z$), then a {\it weight $k$ harmonic weak
Maass form on the congruence subgroup
$\Gamma_1(N)$} is any
smooth function $M:\H\to \C$ satisfying the following:
\begin{enumerate}
\item For all $A= \left(\begin{smallmatrix}a&b\\c&d
\end{smallmatrix} \right)\in \Gamma_1(N)$ and all $\tau\in \H$, we
have
\begin{displaymath}
M\left(\frac{a\tau +b}{c\tau +d}\right)
=\begin{cases} (c\tau +d)^k M(\tau) \ \ \ \ \ &{\text {\rm if}}\
k\in \Z,\\
\leg{c}{d}^{2k}\epsilon_d^{-2k}(cz+d)^{k}\ M(\tau) \ \ \ \ \
&{\text {\rm if}}\ k\in \frac{1}{2}\Z\setminus \Z.
\end{cases}
\end{displaymath}
Here $\leg{c}{d}$ denotes the extended Legendre symbol, and
$\sqrt{\tau}$ is the principal branch of the holomorphic square
root.
\item We have that $\Delta_k
M=0$.
\item There is a polynomial $P_M=\sum_{n\leq 0}c^{+}(n)q^n
\in \C[q^{-1}]$ such that
\begin{displaymath}
 M(\tau)-P_M(\exp(2\pi i \tau))=O(e^{-\epsilon v})
\end{displaymath}
as $v\rightarrow+\infty$ for some $\epsilon>0$. Analogous conditions
are required at all cusps.
\end{enumerate}
\end{definition}
\medskip
\noindent
{\it Four remarks.}
\smallskip
\noindent
1) We refer to the polynomial $P_M$ as the {\it
principal part} of $M(\tau)$ at $\infty$. Obviously, if $P_M$ is
non-constant, then $M(\tau)$ has exponential growth at $\infty$.
Similar remarks apply at all cusps.
\smallskip
\noindent 2) Bruinier and Funke \cite{BF} give a
slightly different definition for a harmonic weak Maass form. In
place of Definition~\ref{defhwmf} (3), they require that $M(\tau)$ has
at most linear exponential growth at cusps. Zagier's weight 3/2 Maass-Eisenstein
series is a harmonic weak Maass form in this relaxed sense.
\smallskip
\noindent
3) Since holomorphic functions on $\H$ are harmonic, it follows that
{\it weakly holomorphic modular forms}, those forms whose poles (if any)
are supported at cusps, are harmonic weak Maass forms.
\smallskip
\noindent
4) For $k\in \frac{1}{2}\Z\setminus \Z$,
the transformation law in Definition~\ref{defhwmf} (1) coincides
with those in Shimura's theory of half-integral weight modular forms
\cite{ShimuraHalf}.
\medskip
In this paper we consider weight 1/2 harmonic weak Maass forms.
Harmonic weak Maass forms have distinguished Fourier expansions
which are described in terms of the incomplete Gamma-function
$\Gamma(\alpha;x)$
\begin{equation}\label{Gamma}
\Gamma(\alpha ;x):=\int_{x}^{\infty}e^{-t}t^{\alpha-1}\ dt.
\end{equation}
The following
characterization is straightforward
(for example,
see Section 3 of \cite{BF}).
If $f(\tau)\in H_{1/2}(N)$, the space of weight 1/2 harmonic weak Maass forms
on $\Gamma_1(N)$, for some $N$, then its
Fourier expansion is of the form
\begin{equation}\label{fourier}
f(\tau)=\sum_{n\gg -\infty} c_f^+(n) q^n + \sum_{n<0} c_f^-(n)
\Gamma(1/2,4\pi |n|v) q^n.
\end{equation}
One sees that $f(\tau)$ naturally decomposes
into two summands, its {\it holomorphic part}
\begin{equation}
\label{holomorphic_part}
f^{+}(\tau):=\sum_{n\gg -\infty} c_f^+(n) q^n,
\end{equation}
and its {\it non-holomorphic part}
\begin{displaymath}
 f^{-}(\tau):=\sum_{\substack{n<0}}
c_f^-(n)\Gamma(1/2,4\pi|n|v)q^n.
\end{displaymath}
\begin{remark} A harmonic weak Maass form with trivial
non-holomorphic part is a weakly holomorphic modular form.
\end{remark}
Harmonic weak Maass forms are related to classical modular forms
thanks to the properties of differential operators. One
nontrivial relationship depends on the differential operator
\begin{equation}
\xi_w:=2i v^w\cdot\overline{\frac{\partial}{\partial \overline{\tau}}}.
\end{equation}
This operator relates $H_{1/2}(N)$ with $S_{3/2}(N)$, the space of
weight 3/2 cusp forms on $\Gamma_1(N)$.
The following lemma\footnote{The formula for $\xi_{1/2}(f)$ corrects
a typographical error in \cite{BF}.}, which is a
straightforward refinement of a proposition of Bruinier and Funke
(see Proposition 3.2 of \cite{BF}), plays a central
role.
\begin{lemma}\label{xiprop}
If $f\in H_{1/2}(N)$, then
\begin{displaymath}
\xi_{1/2}: H_{1/2}(N)\longrightarrow S_{3/2}(N).
\end{displaymath}
Moreover, if
\begin{displaymath}
f(\tau)=\sum_{n\gg -\infty} c_f^+(n) q^n + \sum_{n<0} c_f^-(n)
\Gamma(1/2,4\pi |n|v) q^n,
\end{displaymath}
then we have that
\begin{displaymath}
\xi_{1/2}(f)=-(4\pi)^{1/2}\sum_{n=1}^{\infty}\overline{c_f^{-}(-n)}n^{1/2}q^n.
\end{displaymath}
\end{lemma}
Thanks to Lemma~\ref{xiprop}, we are in a position to relate the
non-holomorphic parts of harmonic weak Maass forms, the expansions
\begin{displaymath}
 f^{-}(\tau):=\sum_{\substack{n<0}}
c_f^-(n)\Gamma(1/2,4\pi|n|v)q^n,
\end{displaymath}
with the ``period integral'' of the modular form $\xi_{1/2}(f)\in
S_{3/2}(N)$. This observation was
critical in Zwegers's work on Ramanujan's mock theta functions.
To make this precise, we relate the Fourier expansion of the
cusp form $\xi_{1/2}(f)$ with $f^{-}(\tau)$. This connection is made by
applying the simple integral identity
\begin{equation}\label{int}
\int_{-\overline{\tau}}^{i \infty}\frac{e^{2\pi i n z}}
{\left(-i(\tau+z)\right)^{1/2}}\ dz= i (2 \pi n)^{-1/2} \cdot
\Gamma(1/2,4 \pi n v) q^{-n}.
\end{equation}
This identity follows by the direct calculation
\begin{displaymath}
 \int_{-\overline{\tau}}^{i \infty}\frac{e^{2\pi i n z}}
{\left(-i(\tau+z)\right)^{1/2}}
  \ dz =
  \int_{2iv}^{i\infty}\frac{e^{2\pi i n(z-\tau)}}
  {(-iz)^{1/2}}\ dz
= i (2 \pi n)^{-1/2}\cdot \Gamma(1/2,4 \pi n v)\, q^{-n}.
\end{displaymath}
By interchanging summation with integration, this identity then implies
that the period integral
\begin{displaymath}
\int_{-\overline{\tau}}^{i\infty} \frac{\overline{\xi_{1/2}(f)(-\overline{z})}}
{(-i(\tau+z))^{1/2}} dz
\end{displaymath}
is proportional to
$\overline{f^{-}}$, the non-holomorphic part of $\overline{f}$.
In this way the non-holomorphic parts of weight
$1/2$ harmonic weak Maass forms are period integrals of weight $3/2$
cusp forms.
Zagier refers to $\xi_{1/2}(f)$ as the {\it shadow} \cite{ZagierBourbaki}
of  $f^+$. The holomorphic part of a weight 1/2 harmonic
weak Maass form is called a {\it mock theta function} if $\xi_{1/2}(f)$
is a linear combination of weight 3/2 theta series
$\vartheta(a,b;\delta \tau)$, where
\begin{displaymath}
\vartheta(a,b;\tau):=\sum_{n\equiv a\pmod b}nq^{n^2}.
\end{displaymath}
The mock theta functions which are relevant for Theorem~\ref{MainTheorem}
are holomorphic parts of Maass forms whose shadows
turn out to be proportional to
\begin{displaymath}
\eta(8\tau)^3=\sum_{n=0}^{\infty}(-1)^n(2n+1)q^{(2n+1)^2}.
\end{displaymath}
\section{Work of Zwegers}\label{Section6}
In his Ph.D. thesis on mock theta functions \cite{Z2}, Zwegers constructs weight 1/2
harmonic weak Maass forms by making use of the transformation properties of
functions which were investigated earlier by Appell and Lerch.
Here we briefly recall some of his results.
For $\tau\in \H$, and $u, v\in \C\setminus (\Z\tau+\Z)$, Zwegers defines the function
\begin{equation}\label{mu}
\mu(u,v;\tau):=\frac{a^{1/2}}{\theta(v;\tau)}\cdot
\sum_{n\in \Z}\frac{(-b)^nq^{n(n+1)/2}}{1-aq^n},
\end{equation}
where $a:=e^{2\pi i u}$, $b:=e^{2\pi i v}$ and
\begin{equation}\label{JTP}
\theta(v;\tau):=\sum_{\nu\in \Z+\frac{1}{2}}(-1)^{\nu-\frac{1}{2}}b^{\nu}
q^{\nu^2/2}.
\end{equation}
Zwegers (see Section 1.3 of \cite{Z2}) proves that $\mu(u,v;\tau)$ satisfies the following important properties.
\begin{lemma}\label{muProperties}
Assuming the notation above, we have that
\begin{displaymath}
\begin{split}
\mu(u,v;\tau)&=\mu(v,u;\tau),\\
\mu(u+1,v;\tau)&=-\mu(u,v;\tau),\\
a^{-1}bq^{-\frac{1}{2}}\mu(u+\tau,v;\tau)&=-\mu(u,v;\tau)+a^{-\frac{1}{2}}b^{\frac{1}{2}}q^{-\frac{1}{8}},\\
\mu(u,v;\tau+1)&=\zeta_8^{-1}\mu(u,v;\tau) \ \ \ \ \ (\zeta_N:=e^{2\pi i/N})\\
(\tau/i)^{-\frac{1}{2}}e^{\pi i (u-v)^2/\tau}\mu\left(\frac{u}{\tau},
\frac{v}{\tau};-\frac{1}{\tau}\right)&=-\mu(u,v;\tau)+\frac{1}{2}h(u-v;\tau),\\
\end{split}
\end{displaymath}
where
\begin{displaymath}
h(z;\tau):=\int_{-\infty}^{\infty}\frac{e^{\pi i x^2\tau-2\pi xz} dx}{\cosh \pi x}.
\end{displaymath}
\end{lemma}
\begin{remark}
The integral $h(z;\tau)$ is known as a {\it Mordell integral}.
\end{remark}
Lemma~\ref{muProperties} shows that $\mu(u,v;\tau)$ is nearly a weight 1/2
Jacobi form, where $\tau$ is the modular variable.
Zwegers then uses $\mu$ to construct weight 1/2 harmonic weak Maass forms.
He achieves this by modifying $\mu$ to obtain a function $\widehat{\mu}$
which he then uses as building blocks for such Maass forms.
To make this precise, for $\tau\in \H$ and $z\in \C$, let
\begin{equation}\label{R}
R(z;\tau):=\sum_{\nu\in \Z+\frac{1}{2}} (-1)^{\nu-\frac{1}{2}}
\left\{\sgn(\nu)-E\left( (\nu+\im(z)/\im(\tau))\sqrt{2\im(\tau)}\right) \right\}
e^{-2\pi i \nu z}q^{-\nu^2/2},
\end{equation}
where $E(z)$ is the odd function
\begin{equation}\label{E}
E(z):=2\int_{0}^{z}e^{-\pi u^2} du.
\end{equation}
Using $\mu$ and $R$,  we let
\begin{equation}\label{muHat}
\widehat{\mu}(u,v;\tau):=\mu(u,v;\tau)-\frac{1}{2}R(u-v;\tau).
\end{equation}
Zwegers's construction of weight 1/2 harmonic weak Maass forms
depends on the following theorem (see Section 1.4 of \cite{Z2} and \cite{ZagierBourbaki}).
\begin{theorem}\label{ZwegersThm}
Assuming the notation and hypotheses above, we have that
\begin{displaymath}
\begin{split}
\widehat{\mu}(u,v;\tau)&=\widehat{\mu}(v,u;\tau)\\
\widehat{\mu}(u+1,v;\tau)&=a^{-1}bq^{-\frac{1}{2}}
\widehat{\mu}(u+\tau,v;\tau)=-\widehat{\mu}(u,v;\tau)\\
\zeta_8\widehat{\mu}(u,v;\tau+1)&=-(\tau/i)^{-\frac{1}{2}}e^{\pi i (u-v)^2/\tau}
\widehat{\mu}\left(\frac{u}{\tau},\frac{v}{\tau};-\frac{1}{\tau}\right)=
\widehat{\mu}(u,v;\tau).
\end{split}
\end{displaymath}
\end{theorem}
Theorem~\ref{ZwegersThm} gives the modular transformation properties
for $\widehat{\mu}$. Since $R$ is non-holomorphic, this theorem,
combined with the following lemma, allows us to relate suitable specializations
of $\widehat{\mu}$ to weight 1/2 harmonic weak Maass forms whose
shadows are linear combinations of weight 3/2 theta functions.
\begin{lemma}\label{Rlemma}{\text {\rm  [Lemma 1.8 of \cite{Z2}]}}
The function $R$ is real analytic and satisfies
\begin{displaymath}
\frac{\partial R}{\partial \overline{u}}(u;\tau)=\sqrt{2}i
v^{-\frac{1}{2}}e^{-2\pi a^2 v}\theta(\overline{u};-\overline{\tau}),
\end{displaymath}
where $a:=\im(u)/\im(\tau)$.  Moreover, we have that
\begin{displaymath}
\frac{\partial}{\partial \overline{\tau}} R(a\tau-b;\tau)=-
\frac{i}{\sqrt{2v}}e^{-2\pi a^2 v}\sum_{\nu\in \Z+\frac{1}{2}}
(-1)^{\nu-\frac{1}{2}}(\nu+a)e^{-\pi i \nu^2\overline{\tau}-
2\pi i \nu(a\overline{\tau}-b)}.
\end{displaymath}
\end{lemma}
\section{Some $q$-series identities}\label{qIdentities}
Here we relate some important $q$-series to the $\mu$-function.
For convenience, we first recall our notation for the following
theta-functions which are simple Dedekind-eta quotients:
\begin{equation}\label{theta2}
\Theta_2(\tau):=\frac{\eta(16\tau)^2}{\eta(8\tau)}=\sum_{n=0}^{\infty}q^{(2n+1)^2}=
q+q^9+q^{25}+\cdots,
\end{equation}
\begin{equation}\label{theta3}
\Theta_3(\tau):=\frac{\eta(8\tau)^5}{\eta(4\tau)^2\eta(16\tau)^2}=1+2\sum_{n=1}^{\infty}
q^{4n^2}=1+2q^4+2q^{16}+2q^{36}+\cdots,
\end{equation}
\begin{equation}\label{theta4}
\Theta_4(\tau):=\frac{\eta(4\tau)^2}{\eta(8\tau)}=1+2\sum_{n=1}^{\infty}(-1)^nq^{4n^2}=
1-2q^4+2q^{16}-2q^{36}+\cdots.
\end{equation}
These are related to the theta functions $\vartheta_2(\tau), \vartheta_3(\tau)$ and $\vartheta_4(\tau)$ by
$$\vartheta_2(\tau)=2\Theta_2\left(\frac{\tau}{8}\right),\ \ \
\vartheta_3(\tau)=\Theta_3\left(\frac{\tau}{8}\right),\ \ \
\vartheta_4(\tau)=\Theta_4\left(\frac{\tau}{8}\right).
$$
Thanks to Theorem~\ref{criterion}, the proof of Theorem~\ref{MainTheorem}
is reduced to the properties of certain specific power series.
For non-negative even $t$, we define the series
\begin{equation}\label{Ft_old}
F_t(q)=q^{-\frac{1}{8}} \, \sum_{\beta=0}^{\infty}\sum_{\alpha=\beta+1}^{\infty}
(-1)^{\alpha+\beta}(2\beta+1)^t q^{\frac{\alpha^2-\beta(\beta+1)}{2}}.
\end{equation}
We shall work with the renormalizations
\begin{equation}\label{Ft}
\mathcal{F}_t(q):=F_t(q^8)=\sum_{\beta=0}^{\infty}\sum_{\alpha=\beta+1}^{\infty}
(-1)^{\alpha+\beta}(2\beta+1)^tq^{4\alpha^2-(2\beta+1)^2}.
\end{equation}
We shall also require the renormalization of the series
$Q^+(\tau)$ in (\ref{Qplus}):
\begin{equation}\label{MockQ}
\mathcal{Q}(q):=Q^+(8\tau)=q^{-1}+28q^3+39q^7+196q^{11}+161q^{15}+756q^{19}+\cdots.
\end{equation}
\subsection{The $\mathcal{F}_t(q)$ series}
Here we interpret the series $\mathcal{F}_t(q)$ in terms of
Zwegers's $\mu$-function. To make this precise, we require the differential
operator
\begin{displaymath}
D_{\omega}:=\frac{1}{2\pi i}\cdot \frac{\partial}{\partial \omega}.
\end{displaymath}
The main result in this section is the following theorem
which expresses the $q$-series $\mathcal{F}_t(q)$ in terms of the image
of the iterated $t$-th derivative, with respect to $\omega$ then evaluated
at $\omega=0$,
of a certain specialization of the $\mu$-function.
\begin{theorem}\label{FasMu}
If $t$ is a non-negative even integer, then
\begin{displaymath}
\frac{\mathcal{F}_t(q)}{\Theta_4(\tau)}=
\frac{1}{2}D_{\omega}^t\left(\mu(4\tau+2\omega,4\tau;8\tau)\right) |_{\omega=0}
=\frac{1}{2}D_{\omega}^t \left(\mu(2\omega,4\tau; 8\tau)\right) |_{\omega=2\tau}.
\end{displaymath}
Moreover, $\frac{\mathcal{F}_0(q)}{\Theta_4(\tau)}$
is the holomorphic part
of the weight 1/2 harmonic weak Maass form
$\frac{1}{2}\widehat{\mu}(4\tau,4\tau;8\tau)$
whose non-holomorphic part
is the period integral of $\eta(8\tau)^3$.
\end{theorem}
\begin{remark} That $t$ is even
will be important
in the proof of Theorem~\ref{MainTheorem} due to the close connection
between $D_{\omega}^t$ and the iterated heat operator.
\end{remark}
\begin{proof}
It is not difficult to show that
\begin{displaymath}
\begin{split}
\mathcal{F}_t&(q)=-\sum_{x,y\geq 0}(2x+1)^tq^{16y^2+16y+16xy+4x+3}+
\sum_{x,y\geq 0}(2x+1)^tq^{16y^2+32y+16xy+12x+15}\\
&=-\sum_{y=0}^{\infty}(-1)^yq^{(2y+1)^2}\sum_{x=0}^{\infty}(2x+1)^tq^{2(2x+1)(2y+1)}\\
&=-\sum_{y=0}^{\infty}(-1)^yq^{(2y+1)^2}\sum_{x=0}^{\infty}
\left(\frac{x^t-(-1)^x x^t}{2}\right)q^{2x(2y+1)}\\
&=-\frac{1}{2}\sum_{y=0}^{\infty}(-1)^yq^{(2y+1)^2}\sum_{x=0}^{\infty}x^tq^{2x(2y+1)}
+\frac{1}{2}\sum_{y=0}^{\infty}(-1)^yq^{(2y+1)^2}\sum_{x=0}^{\infty}(-1)^x x^t q^{2x(2y+1)}.
\end{split}
\end{displaymath}
Using this last expression, and by letting $\rho:=e^{2\pi i \omega}$, it follows that
\begin{displaymath}
\begin{split}
\mathcal{F}_t(q)&=
-\frac{1}{2}\sum_{y=0}^{\infty}(-1)^yq^{(2y+1)^2}\cdot D_{\omega}^t\left(
\sum_{x=0}^{\infty}\rho^x q^{x(4y+2)}\right)\ \vline_{\omega=0}\\
&\ \ \ +
\frac{1}{2}\sum_{y=0}^{\infty}(-1)^yq^{(2y+1)^2}\cdot
D_{\omega}^t\left(\sum_{x=0}^{\infty}(-\rho)^x q^{x(4y+2)}\right)\ \vline_{\omega=0}\\
&=-\frac{1}{2} D_{\omega}^t\left(\sum_{y=0}^{\infty} \left (
\frac{(-1)^yq^{(2y+1)^2}}{1-\rho q^{4y+2}}+
\frac{(-1)^{y+1}q^{(2y+1)^2}}{1+\rho q^{4y+2}}\right )\right)\ \vline_{\omega=0}\\
&=-D_{\omega}^t\left( \sum_{y=0}^{\infty}\frac{(-1)^y \rho q^{4y^2+8y+3}}{1-\rho^2 q^{8y+4}}
\right) \ \vline_{\omega=0}.
\end{split}
\end{displaymath}
Since we have that
\begin{displaymath}
\begin{split}
\sum_{y=0}^{\infty}&\frac{(-1)^y \rho q^{4y^2+8y+3}}{1-\rho^2 q^{8y+4}}=
\frac{1}{2}\sum_{y=0}^{\infty}\frac{(-1)^y\rho q^{4y^2+8y+3}}{1-\rho^2 q^{8y+4}}
+\frac{1}{2}\sum_{y=0}^{-\infty}\frac{(-1)^{y}\rho q^{4y^2-8y+3}}{1-\rho^{2}q^{-8y+4}}\\
&=\frac{1}{2}\sum_{y=0}^{\infty}\frac{(-1)^y\rho q^{4y^2+8y+3}}{1-\rho^2 q^{8y+4}}-
\frac{1}{2}\sum_{y=-1}^{-\infty}\frac{(-1)^{y}\rho q^{4y^2-1}}{1-\rho^2 q^{-8y-4}}\\
&=\frac{1}{2}\sum_{y=0}^{\infty}\frac{(-1)^y\rho q^{4y^2+8y+3}}{1-\rho^2 q^{8y+4}}-
\frac{1}{2}\sum_{y=-1}^{-\infty}\frac{(-1)^y \rho^{-1} q^{4y^2+8y+3}}{\rho^{-2}q^{8y+4}-1}\\
&=\frac{1}{2}\sum_{n=0}^{\infty}\frac{(-1)^n \rho q^{4n^2+8n+3}}{1-\rho^2 q^{8n+4}}
+\frac{1}{2}\sum_{n=-1}^{-\infty}\frac{(-1)^n \rho^{-1}q^{4n^2+8n+3}}{1-\rho^{-2} q^{8n+4}},
\end{split}
\end{displaymath}
the fact that $t$ is even then implies that
\begin{equation}\label{Mrho}
\mathcal{F}_t(q)=-\frac{1}{2}D_{\omega}^t\left( \sum_{n\in \Z}\frac{(-1)^n\rho
q^{4n^2+8n+3}}{1-\rho^2q^{8n+4}}\right) \ \vline_{\omega=0}.
\end{equation}
Since we have that
\begin{displaymath}
\frac{1}{2} D_{\omega}^t \left( \mu(4\tau+2\omega,4\tau;8\tau)\right) |_{\omega=0} =
\frac{1}{2} D_{\omega}^t \left ( \mu(2\omega,4\tau;8\tau)\right) |_{\omega=2\tau},
\end{displaymath}
to complete the proof of the claimed identity,
it suffices to compute the expansion of
$\mu(4\tau+2\omega,4\tau;8\tau)$.
By (\ref{mu}), we have that
\begin{displaymath}
\mu(4\tau+2\omega,4\tau;8\tau)=\frac{q^{-1}}{\theta(4\tau;8\tau)}\sum_{n\in \Z}
\frac{(-1)^n\rho q^{4n^2+8n+3}}{1-\rho^2 q^{8n+4}}.
\end{displaymath}
The claimed identity follows immediately from (\ref{JTP}) and (\ref{Mrho}) since
\begin{displaymath}
q\theta(4\tau;8\tau)=\sum_{m\in \Z}(-1)^mq^{4(m+1)^2}=-\sum_{m\in \Z}(-1)^mq^{4m^2}=
-\frac{\eta(4\tau)^2}{\eta(8\tau)}=-\Theta_4(\tau).
\end{displaymath}
That $\mathcal{F}_0(q)/\Theta_4(\tau)$ is the holomorphic part
of the weight 1/2 harmonic weak Maass form $\frac{1}{2}
\widehat{\mu}(4\tau,4\tau;8\tau)$ follows, after a straightforward
calculation, from
Theorem~\ref{ZwegersThm}. That the non-holomorphic part
of this Maass form is a multiple of the period integral of $\eta(8\tau)^3$ follows
from the classical
identity
\begin{displaymath}
\eta(8\tau)^3=\sum_{n=0}^{\infty}(-1)^n(2n+1)q^{(2n+1)^2},
\end{displaymath}
and the explicit Fourier expansion of the non-holomorphic part
of $\frac{1}{2}\widehat{\mu}(4\tau,4\tau;8\tau)$. This expansion
is obtained using (\ref{R}), Lemma~\ref{Rlemma},
and the discussion after Lemma~\ref{xiprop}.
\end{proof}
\subsection{The $\mathcal{Q}(q)$ series}
Here we describe the $q$-series $\mathcal{Q}(q)$
in terms of Zwegers's $\mu$-function. To make this precise,
we require some auxiliary weight 1/2 weakly holomorphic modular forms.
Define modular forms $A(\tau)$ and $B(\tau)$ by
\begin{equation}
\label{B-functions}
\begin{split}
\mathcal{A}(\tau):=A(8\tau)&=\sum_{n=-1}^{\infty}a(n)q^n:=\frac{\eta(4\tau)^8}{\eta(8\tau)^7}=
q^{-1}-8q^3+27q^7-\cdots,\\
\mathcal{B}(\tau):=B(8\tau)&=\sum_{n=-1}^{\infty}b(n)q^n:=\frac{\eta(8\tau)^5}{\eta(16\tau)^4}=
q^{-1}-5q^7+9q^{15}-\cdots.
\end{split}
\end{equation}
We sieve on the Fourier expansion of $\mathcal{A}(\tau)$ to define the modular
forms
\begin{equation}
\label{A-functions}
\begin{split}
\mathcal{A}_{3,8}(\tau)&:=A_{3,8}(8\tau)=\sum_{n\equiv 3\pmod 8}a(n)q^n=-8q^3-56q^{11}+\cdots,\\
\mathcal{A}_{7,8}(\tau)&:=A_{7,8}(8\tau)=\sum_{n\equiv 7\pmod 8}a(n)q^n=q^{-1}+27q^7+105q^{15}+\cdots.
\end{split}
\end{equation}
In terms of these modular forms, we have the following theorem.
\begin{theorem}\label{QasMu}
The following $q$-series identity is true:
\begin{displaymath}
\begin{split}
\mathcal{Q}(q)=-\frac{7}{2}&\mathcal{A}_{3,8}(\tau)+\frac{3}{2}\mathcal{A}_{7,8}(\tau)-
\frac{1}{2}\mathcal{B}(\tau)\\
&+2iq^{-1}\mu\left(-16\tau,-8\tau-\frac{1}{2};32\tau\right)
-2iq^{-1}\mu\left(-16\tau,-24\tau-\frac{1}{2};32\tau\right).
\end{split}
\end{displaymath}
Moreover, $\mathcal{Q}(q)$ is the holomorphic part of a weight 1/2
harmonic weak Maass form whose non-holomorphic part is the period
integral of $\eta(8\tau)^3$.
\end{theorem}
\begin{remark}
One can determine the image of $\mathcal{Q}(q)$ under the inversion
map $\tau\rightarrow -1/\tau$
using the formula for $\mathcal{Q}(q)$ above.
One obtains this expansion by applying the third part of
Theorem~\ref{ZwegersThm} to two $\mu$-functions, and by applying the
classical transformation
\begin{displaymath}
\eta(-1/\tau)=\sqrt{-i \tau}\cdot \eta(\tau)
\end{displaymath}
to the weakly holomorphic modular forms
\begin{displaymath}
\begin{split}
\mathcal{A}_{3,8}(\tau)&=-8\cdot\frac{\eta(16\tau)^8}{\eta(8\tau)^7}=-8q^3-56q^{11}-216q^{19}-\cdots,\\
\mathcal{A}_{7,8}(\tau)&=\frac{\eta(8\tau)^5}{\eta(16\tau)^4}+32\cdot\frac{\eta(32\tau)^8}{\eta(8\tau)^3
\eta(16\tau)^4}=q^{-1}+27q^7+105q^{15}+\cdots,\\
\mathcal{B}(\tau)&=\frac{\eta(8\tau)^5}{\eta(16\tau)^4}=
q^{-1}-5q^7+9q^{15}-\cdots.
\end{split}
\end{displaymath}
\end{remark}
\begin{proof}[Sketch of the proof of Theorem~\ref{QasMu}]
This result follows from Theorem~\ref{ZwegersThm} in an argument
which is analogous to the proof of Theorem~\ref{FasMu}.
Namely, since $\mathcal{A}_{3,8}(\tau), \mathcal{A}_{7,8}(\tau)$ and $B(\tau)$ are
weight 1/2 modular forms, it follows that
$\mathcal{Q}(\tau)$ is the holomorphic part of a weight 1/2 harmonic
weak Maass form whose non-holomorphic part is a multiple of the period integral
of $\eta(8\tau)^3$. Therefore, it must be a linear combination
of certain mock theta functions and weakly holomorphic modular forms
of weight 1/2.
Instead of using Theorem~\ref{FasMu} with $t=0$, we
work with the  mock theta function
\begin{equation}
\label{M-function}
\mathcal{M}(q):=q^{-1}\sum_{n=0}^{\infty}\frac{(-1)^{n+1}q^{8(n+1)^2}\prod_{k=1}^n
(1-q^{16k-8})}{\prod_{k=1}^{n+1}(1+q^{16k-8})^2}=-q^7+2q^{15}-3q^{23}+\cdots,
\end{equation}
and $\mathcal{M}(q)=M(q^8)$.
Using  Watson's $q$-analog of Whipple's theorem, we have that
\begin{displaymath}
\mathcal{M}(q)=-\frac{1}{2\Theta_2(\tau)}\cdot \sum_{n\in \Z}\frac{q^{16n^2-8n}}{1+q^{16n-8}}.
\end{displaymath}
Andrews \cite{Andrews} already interpreted\footnote{
We have reformulated Andrews's work to be consistent with the terminology
in this paper.}
this particular $q$-series
as the holomorphic part of a weight 1/2 harmonic weak Maass form
whose non-holomorphic part is the period integral of $\eta(8\tau)^3$.
Using the fact that
\begin{displaymath}
\sum_{n\in \Z}\frac{q^{16n^2-8n}}{1+q^{16n-8}}=
\sum_{n\in \Z}\frac{q^{16n^2-8n}(1-q^{16n-8})}{1-q^{32n-16}}=
\sum_{n\in \Z}\frac{q^{16n^2-8n}}{1-q^{32n-16}}-
\sum_{n\in \Z}\frac{q^{16n^2+8n-8}}{1-q^{32n-16}},
\end{displaymath}
we find that
\begin{displaymath}
\mathcal{M}(q)=-\frac{iq^{-1}}{2}\cdot\left(
\mu\left(-16\tau,-24\tau-\frac{1}{2};32\tau\right)
-\mu\left(-16\tau,-8\tau-\frac{1}{2};32\tau\right)\right).
\end{displaymath}
Here we have used the fact that
\begin{displaymath}
q\theta\left(-8\tau-\frac{1}{2};32\tau\right)=q^9\theta\left(-24\tau-\frac{1}{2};32\tau\right)=i\Theta_2(\tau).
\end{displaymath}
Since two weight 1/2 harmonic weak Maass forms with the same
non-holomorphic parts differ by a weakly holomorphic modular form,
it then follows that $\mathcal{Q}(q)-4M(q)$ is a weakly holomorphic
modular form. Standard calculations in the theory of modular forms reveals
that
\begin{displaymath}
\mathcal{Q}(q)-4\mathcal{M}(q)=-\frac{7}{2}\mathcal{A}_{3,8}(\tau)+\frac{3}{2}\mathcal{A}_{7,8}(\tau)-
\frac{1}{2}\mathcal{B}(\tau).
\end{displaymath}
\end{proof}
\begin{remark}
One also sees directly that the shadow of $\mathcal{Q}(q)$
is proportional to $\eta(8\tau)^3$ using the proven identity.
One applies Lemma~\ref{Rlemma} to this particular linear
combination $\widehat{\mu}$ functions.
\end{remark}
\section{Proof of Theorem~\ref{MainTheorem} and some numerical examples}
\label{TheProof}
Here we prove Theorem~\ref{MainTheorem} by verifying the criterion
in Theorem~\ref{criterion}. We also conclude this section with some
numerical examples which illustrate the phenomenon which appears in the
proof of the theorem.
\subsection{Proof of Theorem~\ref{MainTheorem}}
Thanks to Theorem 4.26, it suffices to prove that the differences between certain $q$-series have vanishing constant term. We shall derive these conclusions by using differential operators. For brevity we describe the $n = 0$
cases in detail, and then provide general remarks which are required to justify the
remaining cases.
By Theorems 7.1 and 7.2, it follows that both
\[
\mathcal{Q}(q) = q^{-1} + 28q^3 + 39q^7 + 196q^{11} + 161q^{15} + 756q^{19} + \dots ;
\]
\[
-\frac{4\mathcal F_0(q)}{\Theta_4(\tau)} = 4q^3 + 12q^7 + 28q^{11} + \dots
\]
are the holomorphic parts of weight 1/2 harmonic weak Maass forms with equal non-holomorphic parts. Therefore, it follows that
\[\mathcal  Q(q) + \frac{4 \mathcal F_0(q)}{\Theta_4(\tau)} = q^{-1} + 24q^3 + 27q^7 + 168q^{11} + \dots\]
is a weakly holomorphic modular form. Standard calculations using the theory of
modular forms reveals that
\begin{eqnarray}
Z_0(q) := \mathcal Q(q) + \frac{4\mathcal F_0(q)}{\Theta_4(\tau)} = \frac{E^*(4\tau)}{\eta(8\tau)^3},
\end{eqnarray}
where $E^*(4\tau)$ is the weight 2 Eisenstein series
\begin{eqnarray}
E^*(\tau):=-E_2(\tau)+2E_2(2\tau)=1+24\sum_{n=1}^\infty\sigma_{\mathrm{odd}}(n)q^n,
\end{eqnarray}
where $\sigma_{\mathrm{odd}}(n)$ denotes the sum of the positive odd divisors of $n$.
The $n=0$ cases of Theorem 4.26 are equivalent to the claim that the constant terms of
\begin{equation}\label{n0}
\frac{\Theta_4(\tau)^9(16\Theta_2(\tau)^4+\Theta_3(\tau)^4)^mZ_0(q)}{\Theta_2(\tau)^{2m+3}\Theta_3(\tau)^{2m+3}}
\end{equation}
are $0$ for every $m\geq 0$.
In order to verify this claim, we will find if helpful to define
\begin{eqnarray}
\widehat{Z_0}(q) := \frac{E^*(4\tau)}{\Theta_2(\tau)^2\Theta_3(\tau)^2} = \frac{Z_0(q) \eta(8\tau)^3}{\Theta_2(\tau)^2\Theta_3(\tau)^2}.
\end{eqnarray}
A calculation shows that
\begin{eqnarray}
q\frac{d}{dq}\widehat{Z_0}(q) = \frac{\Theta_4(\tau)^9}{\Theta_2(\tau)\Theta_3(\tau)\eta(8\tau)^3}.
\end{eqnarray}
Using this notation, and noting that
\[16\Theta_2(\tau)^4+\Theta_3(\tau)^4=1+24q^4+24q^2+\cdots=E^*(4\tau),\]
equation (\ref{n0}) becomes
\[\widehat{Z_0}(q)^{m+1}\cdot q\frac{d}{dq}\widehat{Z_0}(q),\]
which clearly has a vanishing constant term.
In fact, for each $m,n\geq 0$ we find a similar phenomenon.

To this end, we make use of the combinatorial structure of series $\mathcal G_{k}(q)$ and $H_k(q)$ defined below.
These series are obtained by iteratively applying differential operators to the harmonic Maass forms and modular forms in this paper,
which in turn lead to a nice sequence of modular functions of level 8.
The modular properties of these series follow easily from the results and methods of  (pp. 53-55 of \cite{Zagier123}  and \cite{ChoieLee}).

For every non-negative $k$, define
\begin{equation}
\begin{split}
&\mathcal G _{k}(q):= \\
&\ \ \frac{\eta(8\tau)^3}{\left(\Theta_2(\tau)\Theta_3(\tau)\right)^{2k+2}} \sum_{j=0}^k  \begin{pmatrix}k\\j \end{pmatrix}\frac{\Gamma(\frac12) (-12)^j E_2(8\tau)^{k-j}}{\Gamma\left(\frac 12+j\right)8^j} \left[\frac{(-1)^j 4\mathcal F_{2j}(q)}{\Theta_4(\tau)} +\left(q\frac{d}{dq}\right)^j\mathcal Q(q)  \right].
\end{split}
\end{equation}
Theorem 4.26 is equivalent to the claim that the constant coefficient of
\begin{equation}\label{polyZ}
\left(q\frac{d}{dq}\widehat{Z_0}(q)\right) \widehat{Z_0}(q)^m\sum_{k=0}^n  \begin{pmatrix}n\\k \end{pmatrix}\widehat{Z_0}(q)^{n-k}\mathcal G _{k}(\tau)
\end{equation}
 is zero for each non-negative $m$ and $n$.
 To prove the theorem, it suffices to show that $\mathcal G _{k}(q)$ is a polynomial in $\widehat{Z_0}(q)$. To this end, we define $M_0^*(\Gamma_0(8))$ to be the space of modular function on $\Gamma_0(8)$ which are holomorphic away from infinity, and is a subspace of $\C((q^2))$. One can easily  verify that $M_0^*(\Gamma_0(8))$ is precisely the set of polynomials in $\widehat{Z_0}(q)$.
 In order to show that $\mathcal G _{k}(\tau)$ is in $M_0^*(\Gamma_0(8))$, we first show that a similar function, $\mathcal H_k(q)$ is in $M_0^*(\Gamma_0(8))$. We define the function
\begin{equation}
\mathcal H _{k}(q):= \frac{\eta(8\tau)^3}{\left(\Theta_2(\tau)\Theta_3(\tau)\right)^{2k+2}} \sum_{j=0}^k  \begin{pmatrix}k\\j \end{pmatrix}\frac{\Gamma(\frac12) (-12)^j E_2(8\tau)^{k-j}}{\Gamma\left(\frac 12+j\right)8^j} \left(q\frac{d}{dq}\right)^jZ_0(q).
\end{equation}
We can observe that $\mathcal H _{k}(q)$ is modular on $\Gamma_0(8)$ with weight $0$ by comparing the summation to the expression $\mathcal Z_k (q):= \mathcal E_\frac 12^k[Z_0(q^{1/8})]$ where $\mathcal E_\frac 12^k[*]$ is defined as in Lemma~4.10. A calculation shows that  $\left(\Theta_2(\tau)\Theta_3(\tau)\right)^{-2}$ is holomorphic away from infinity, which, combined with the fact that $\frac{\eta(8\tau)^3 Z_0(q)}{\left(\Theta_2(\tau)\Theta_3(\tau)\right)^{2}}=\widehat{Z_0}(q)$, shows that $\mathcal H_k(\tau)$ is in $M_0^*(\Gamma_0(8))$. Hence we need only show that $\mathcal G_k(q)-\mathcal H_k(q)$ is in $M_0^*(\Gamma_0(8))$ as well. We observe that
\begin{displaymath}
\begin{split}
&\mathcal G _{k}(q)-\mathcal H_k(q) = \\
&\ \frac{\eta(8\tau)^3}{\left(\Theta_2(\tau)\Theta_3(\tau)\right)^{2k+2}} \sum_{j=0}^k  \begin{pmatrix}k\\j \end{pmatrix}\frac{\Gamma(\frac12) (-12)^j E_2(8\tau)^{k-j}}{\Gamma\left(\frac 12+j\right)8^j} \left[ \frac{(-1)^j 4\mathcal F_{2j}(q)}{\Theta_4(\tau)} - \left(q\frac{d}{dq}\right)^j \frac{4\mathcal F_{0}(q)}{\Theta_4(\tau)} \right].
\end{split}
\end{displaymath}
Using Theorem 7.1, this can be written as
\begin{displaymath}
\begin{split}
& \frac{\eta(8\tau)^3}{\left(\Theta_2(\tau)\Theta_3(\tau)\right)^{2k+2}} \sum_{j=0}^k  \begin{pmatrix}k\\j \end{pmatrix}\frac{\Gamma(\frac12) (-12)^j E_2(8\tau)^{k-j}}{\Gamma\left(\frac 12+j\right)8^j}\\ &\ \ \ \ \
 \cdot \left.2\left[ (-1)^j \left(\frac{1}{2\pi i}\frac{d}{d\omega}\right)^{2j} \mu(4\tau+2\omega,4\tau;8\tau) - \left(\frac{1}{2\pi i}\frac{d}{d\tau}\right)^j \mu(4\tau+2\omega,4\tau;8\tau) \right]\right | _{\omega=0}.
 \end{split}
\end{displaymath}
Using the transformation laws for $\mu$ found in Lemma 6.1, we observe that the Mordell integrals that arise as obstructions to the modular transformation of
\[(-1)^j \left(\frac{1}{2\pi i}\frac{d}{d\omega}\right)^{2j} \mu(4\tau+2\omega,4\tau;8\tau) - \left(\frac{1}{2\pi i}\frac{d}{d\tau}\right)^j \mu(4\tau+2\omega,4\tau;8\tau)\] cancel directly. Therefore an argument similar to that used in the proof of Lemma 4.10  shows that $\mathcal G _{k}(q)-\mathcal H_k(q)$ is modular with respect to $\Gamma_0(8)$. It is then straightforward to verify that $\mathcal G _{k}(q)-\mathcal H_k(q)$ is in $M_0^*(\Gamma_0(8))$, and
so is a polynomial in $\widehat{Z_0}(q)$. This then completes the proof.
\subsection{Examples}
Here we give some numerical examples of Theorem~\ref{criterion}.
First we give some examples when $n=0$, and then we give an example
when $n=1$. We conclude with examples of the polynomials in
$\widehat{Z_0}(q)$ which are central to the proof of Theorem~\ref{MainTheorem}.
\begin{example}
Let $Z_0(q)=E^*(4\tau)/\eta(8\tau)^3$ and define $f_{m}(\tau)$
by
$$
f_{m}(\tau):=\frac{\Theta_4(\tau)^9\left(16\Theta_2(\tau)^4+\Theta_3(\tau)^4\right)^m}
{\Theta_2(\tau)^{2m+3}\Theta_3(\tau)^{2m+3}}.
$$
The proof shows that the constant
terms of $Z_0(q)f_m(\tau)$ (see (\ref{n0})) vanish for all $m$.
The initial terms of the expansions of the first few $f_m(\tau)$
are
\begin{displaymath}
\begin{split}
f_0(\tau)&=q^{-3}-24q+273q^5-1976q^9+\cdots,\\
f_1(\tau)&=q^{-5}-4q^{-1}-269q^3+5188q^7+\cdots,\\
f_2(\tau)&=q^{-7}+16q^{-3}-411q+272q^5+\cdots,\\
f_3(\tau)&=q^{-9}+36q^{-5}-153q^{-1}-\cdots,
\end{split}
\end{displaymath}
and the initial terms of $Z_0(q)f_m(\tau)$ are
\begin{displaymath}
\begin{split}
Z_0(q)f_0(\tau)&=q^{-4}-276q^4+4096q^8-33606q^{12}+\cdots,\\
Z_0(q)f_1(\tau)&=q^{-6}+20q^{-2}-338q^2-1208q^6+\cdots,\\
Z_0(q)f_2(\tau)&=q^{-8}+40q^{-4}-8992q^4+65260q^8+\cdots,\\
Z_0(q)f_3(\tau)&=q^{-10}+60q^{-6}+738q^{-2}-11256q^2-\cdots.
\end{split}
\end{displaymath}
One easily sees that the constant terms of these $Z_0(q)f_m(\tau)$
indeed vanish.
\end{example}
\begin{example}
Here we simplify the notation in Theorem~\ref{criterion} by
defining $q$-series $\widehat{\D}^{(i)}_{m,n}(q)$ and
$\Lambda^{(i)}(m,n,k,j;q)$ so that
\begin{equation}
\widehat{\D}^{(i)}_{m,n}(q):=\D_{m,n}^{(i)}(q^8)=\sum_{k=0}^{n}\sum_{j=0}^k
\Lambda^{(i)}(m,n,k,j;q).
\end{equation}
By Theorem~\ref{criterion}, Theorem~\ref{MainTheorem} follows from
the assertion that the constant term of
$\widehat{\D}^{(1)}_{m,n}(q)-\widehat{\D}^{(2)}_{m,n}(q)$
vanishes.
The first few terms of
the series $\Lambda^{(i)}(3,1,k,j;q)$ are:
\begin{displaymath}
\begin{split}
\Lambda^{(1)}(3,1,0,0;q)&=\frac{1}{256}q^{-8}+\frac{43}{256}q^{-4}+\frac{7}{16}
-\cdots,\\
\Lambda^{(1)}(3,1,1,0;q)&=\frac{1}{768}q^{-8}+\frac{35}{768}q^{-4}-\frac{13}{48}
-\cdots,\\
\Lambda^{(1)}(3,1,1,1;q)&=-\frac{1}{768}q^{-8}-\frac{59}{768}q^{-4}-
\frac{85}{96}+\cdots,\\
\Lambda^{(2)}(3,1,0,0;q)&=-\frac{1}{3072}q^{-12}-\frac{7}{256}q^{-8}
-\frac{11}{16}q^{-4}-\frac{85}{96}+\cdots,\\
\Lambda^{(2)}(3,1,1,0;q)&=\frac{1}{3072}q^{-12}+\frac{5}{256}q^{-8}
+\frac{13}{64}q^{-4}-\frac{247}{48}-\cdots,\\
\Lambda^{(2)}(3,1,1,1;q)&=\frac{1}{1024}q^{-12}-\frac{13}{256}q^{-8}
-\frac{203}{64}q^{-4}+\frac{85}{16}+\cdots.
\end{split}
\end{displaymath}
One sees that the constant term of $\widehat{\D}^{(1)}_{3,1}(q)-
\widehat{\D}^{(2)}_{3,1}(q)$ is
\begin{displaymath}
\frac{7}{16}-\frac{13}{48}-\frac{85}{96}+\frac{85}{96}
+\frac{247}{48}-\frac{85}{16}=0.
\end{displaymath}
Notice that the constant terms of $\Lambda^{(1)}(3,1,1,1;q)$ and
$\Lambda^{(2)}(3,1,0,0)$ agree.  This equality is a special case
of a general equality whose proof is equivalent to the $n=0$ case
of the proof of Theorem~\ref{MainTheorem}. In general, the constant
terms of $\Lambda^{(1)}(m,n,0,0;q)$ and
$\Lambda^{(2)}(m,n,n,n;q)$ agree.
\end{example}
\begin{example}
The proof of Theorem~\ref{MainTheorem} relies on the fact that the expressions in
Theorem~4.26 can be written as (see (\ref{polyZ}))
\begin{eqnarray}
\left(q\frac{d}{dq}\widehat{Z_0}(q^{1/8})\right) \widehat{Z_0}(q^{1/8})^mP_n(\widehat{Z_0}(q^{1/8})),
\end{eqnarray}
where $P_n(x)$ is a polynomial.
In the table below we give the first few  polynomials.
\begin{center}
\begin{tabular}{|l|l|c|l|l|}
\hline
$n$& $P_n(x)$\\
\hline
& \ \ \\
0&$x$\\ & \ \ \\
1&$-\frac{1}{2}x^2 + 72$\\ & \ \ \\
2&$\frac{1}{12}x^3 + \frac{32}{3}x$\\ & \ \ \\
3&$-\frac{1}{120}x^4 + \frac{122}{15}x^2 - 64$\\ & \ \ \\
4&$\frac{1}{1680} x^5 + \frac{8}{7} x^3 + \frac{5696}{105}x$\\ & \ \ \\
5&$-\frac{1}{30240} x^6 + \frac{457}{1890} x^4 + \frac{9536}{315} x^2 +
\frac{406016}{315}$\\ & \ \ \\
6&$\frac{1}{665280}x^7 + \frac{4}{135}x^5 + \frac{1712}{315}x^3 -
\frac{2089984}{3465}x$\\ \ \ & \\
\hline
\end{tabular}
\end{center}
\end{example}

\end{document}